\documentclass[letterpaper,12pt,leqno,oneside]{article}
\usepackage{color}
\usepackage{bm}
\usepackage{exscale}
\usepackage{amsmath}
\usepackage{amsfonts}
\usepackage{stmaryrd}
\usepackage{amscd}
\usepackage{graphicx}
\usepackage{amsxtra}
\usepackage{amssymb}
\usepackage{theorem}
\usepackage[final]{epsfig}
\usepackage{eqnarray}
\newtheorem{proposition}{Proposition}[subsection]
\newtheorem{definition}[proposition]{Definition}

\newtheorem{lemma}[proposition]{Lemma}

{\theorembodyfont{\rmfamily}\newtheorem{remark}[proposition]{Remark}}

\newtheorem{corollary}[proposition]{Corollary}

{\theorembodyfont{\rmfamily}\newtheorem{example}[proposition]{Example}}

\newfont{\abc}{cmtt10 scaled 1200}

\def\R{\mathbb{R}}

\def\B{\mathbf{A}}

\def\Z{\mathbb{Z}}

\def\U{\mathbb{U}}
\def\D{\mathbb{D}}

\def\C{\mathbb{C}}

\def\E{\mathbb{E}}
\def\U{\mathbb{U}}

\def\I{\mathbb{I}}
\def\ve{\varepsilon}

\def\ra{\rightarrow}
\def\cs{\symbol{35}}
\def\p{\partial}
\def\qed{\hfill $\Box$ \\}
\def\mm{\mbox}
\def\v{= \emptyset}
\def\n{\neq \emptyset}
\def\D{\mathbf{ID}}
\def\M{\mathbb{A}}
\def\RH{reg^{|A|}_H}
\def\db{d^{\,\flat}}
\def\su{\mathbb{SD}}

\def\llb{\llbracket}
\def\rrb{\rrbracket}

\def\sk{$|A|$-skin}
\def\sks{$|A|$-skins}
\def\bp{\langle A \rangle}

\begin{document}

\vspace*{0cm}

\begin{center}\Large{\bf{Skin Structures on Minimal Hypersurfaces}}\\
\medskip
{\small{by}}\\
\medskip
\large{\bf{Joachim Lohkamp}\\
\medskip}

\end{center}

\noindent Mathematisches Institut, Universit\"at M\"unster, Einsteinstrasse 62, Germany\\
 {\small{\emph{e-mail: j.lohkamp@uni-muenster.de}}}

{\small {\center \tableofcontents}

{\contentsline {subsection}{\numberline {}References}{. .}}}

\bigskip
\setcounter{section}{1}
\renewcommand{\thesubsection}{\thesection}
\subsection{Introduction} \label{int}
\bigskip

Variational methods are broadly used in geometry, analysis and physics. Solutions of variational problems can oftentimes be regarded as new geometric spaces conveying information about the original problem. To understand these spaces and to extract the information they carry, involves analytic problems now formulated on these spaces, most obviously, induced from the original variational problems.\\

However, such spaces may carry delicate singularities. Their geometry and also the analysis on these spaces generally degenerate towards these singular sets in a rather complicated way. This causes difficult issues yet hardly understood.\\

The present paper is the first in a series, followed by [L1] and [L2], where we address these issues for the case of singular area minimizing hypersurfaces. This is one of the principal model cases in this field and it exhibits characteristic peculiarities of singular solutions of elliptic variational problems.  Extensions, to cover more general problems, including almost minimizers, will be discussed in later accounts. \\

\textbf{Singular Area Minimzers} \, The known attempts to get grip on \emph{singular} area minimizing hypersurfaces  largely centered on the structure of the singular set and on perturbation methods to avoid the occurrence of singularities. Examples of such approaches are found in the, in part already classical, work of Federer [F1],[F2] and of Schoen, Uhlenbeck, Caffarelli, Hardt, Simon and White, [SU], [CHS], [HS], [Si1], [W]. Two of the best known results say, the singular set has codimension $\ge 7$ within the given hypersurface, [F2], and, in the analytic case, it is a rectifiable set  [Si1].\\

We propose another strategy and shift the focus from singularities to their regular complements. On these open manifolds we introduce new structural elements, the  \emph{skin structures}, disclosing previously unapproachable and largely unexpected geometric and analytic properties of singular but also of regular hypersurfaces.\\

To give the reader a first impression, we informally describe some results we establish in this triple of papers, even before we properly explain what skin structures are. In the next section, Ch.\ref{gint}, we give a broad overview of these topics, and references, guiding through this series of papers. For now, let $H^n \subset M^{n+1}$ be an area minimizing hypersurface with singular set $\Sigma \subset H$, $M$ a smooth compact Riemannian manifold, and let $scal_X$ denote the scalar curvature of the underlying space $X$.\\

$\bullet$\, We establish a basic boundary regularity of $H \setminus \Sigma$, the \emph{skin uniformity} of $H \setminus \Sigma$, where we view $\Sigma$ as a boundary. This is a refined counterpart to concepts of \emph{non tangentially accessibility} of \emph{Euclidean domains}, like that of uniform or NTA-domains, known e.g. from potential theory or quasi-conformal geometry.\\

$\bullet$\, We get a \emph{hyperbolic unfolding} of $H \setminus \Sigma$, that is,  a canonical \emph{conformal Gromov hyperbolization} of $H \setminus \Sigma$ to a complete space of \emph{bounded geometry}. Its \emph{Gromov boundary}  is homeomorphic to the singular set $\Sigma$ of $H$.\\

$\bullet$\, For a large class of elliptic operators, the \emph{skin adapted operators}, the existence of these unfoldings translates to fine controls for their potential theory and asymptotic analysis, near $\Sigma$, on the degenerating original space $H \setminus \Sigma$. For instance, their  \emph{Martin boundary} on $H \setminus \Sigma$ is  homeomorphic to $\Sigma$ and each boundary point is minimal.\\

$\bullet$\,  The \emph{conformal Laplacian} $L_H=-\Delta  +\frac{n-2}{4 (n-1)} \cdot scal_H$ is an interesting sample case. $L_H$ is skin adapted, if and, in general, only if $scal_M \ge 0$. Also we recall the well-known $scal >0$-heredity principle: when $M$ carries a $scal >0$-metric, then $H$ can be conformally deformed to a $scal >0$-space using the first eigenfunction  of $L_H$.\\

$\bullet$\,  Due to the presence of singularities, the latter $scal >0$-heredity principle runs into intractable problems, when we try to use it inductively, to build a structure theory for $scal>0$-geometries. We resolve this classical problem using skin structures, including the consequences of the skin adpatedness of $L_H$. The outcome is a broader strategy, the \emph{$scal >0$-heredity with surgery}. It allows us to also replace the singularities for regular ends and closures while keeping $scal >0$.\\

$\bullet$\,   In turn, there are many applications of the $scal >0$-heredity with surgery in geometry and general relativity. More specific applications are considered in separate accounts. They include proofs of long conjectured general versions of the positive mass theorems, the Penrose inequality and the fact that so-called enlargeable manifolds cannot admit $scal >0$-metrics, results previously only approachable for spin manifolds or in dimensions $\le 7$.

\subsubsection{General Overview}\label{gint}
\bigskip

We start with a broad overview of the main contents of the whole  series with an emphasis on the pervasive use of skin structures. The results of the present foundational paper are detailed in the next section Ch.\ref{sor}.\\

For this introduction we let $(H^n,g_H) \subset (M^{n+1},g_M)$, $ n \ge 2$ be a connected area minimizing hypersurface.
$M^{n+1}$ denotes a compact smooth Riemannian $n+1$-manifold, $\Sigma \subset H$ the singular set of $H$, which may also be empty. In the hypersurface case one actually knows that $\Sigma \v$, when $n < 7$. Therefore we are primarily interested in the higher dimensional case of dimension $\ge 7$.\\

\textbf{Skin Structures} \quad  The basic skin structural notion is that of a \emph{skin transform}, a particular type of non-negative density function $\bp_H$ naturally assigned to $H$.\\

To describe it intuitively,  we consider the level sets  $\B_c:=\{x \in H \setminus \Sigma \,|\, |A_H|(x)=c\}$, for $c >0$,
of the norm $|A_H|$ of the second fundamental form $A_H$ of $H$. These level sets may have all kinds of peculiar properties which reflect an inhomogeneous wrinkling of
$H$ which occurs as we approach $\Sigma$.\\

The idea is to transform each of the $\B_c$ into an equalized enveloping hypersurface $\M_c \subset H \setminus \Sigma$, an \emph{\sk}.
We assemble the collection of all \sks\ to define a \emph{skin transform} $\bp$.
For this, we set $\bp(x):= c, \mm{ for } x  \in \M_c$ and canonically extend the definition to points which do not belong to any \sk .\\

The skin transform $\bp$ detects inhomogeneities in the underlying space $H \setminus \Sigma$ and gives us a means to uniformly reshape and
unfold this delicate and wrinkled geometry.\\

In Ch.2, we describe ways to define such equalizing procedures transforming the $\M_c$ into some $\B_c$ and, from these, skin transforms. These procedures satisfy a common set of simple axioms for the resulting $\bp$. This leads to an important simplification since we can base our applications on only these axioms. For comparison, there are several distinct constructions of the Brownian motion. They all satisfy the same set of axioms and most applications rely on the exclusive use of these few elementary properties. \\

Such a concise axiomatic description for skin transforms reads as follows: An assignment $H \mapsto \bp_H$, of a non-negative, measurable function to any connected area minimizing hypersurface $H$, is called a \emph{skin transform} provided

\begin{itemize}
\item $\bp_H$ is \emph{naturally} assigned to $H$, in other words, the assignment commutes with the convergence of sequences of underlying area minimizers.
    \item $\bp_H \ge |A_H|$ and for any $f \in C^\infty(H \setminus \Sigma,\R)$
        compactly supported in $H \setminus \Sigma$ we have the following Hardy type inequality*:
\[\int_H|\nabla f|^2  + |A_H|^2 \cdot f^2 dA \ge \tau \cdot \int_H \bp_H^2\cdot f^2 dA, \mm{ for some } \tau = \tau(\bp,H) \in (0,1).\]
    \item $\bp_H \equiv 0$, if $H \subset M$ is totally geodesic, otherwise, $\bp_H$ is  strictly positive.
    \item When $H$ is not totally geodesic, the function $\delta_{\bp}:=1/\bp$  is well-defined and it is
        $L_{\bp}$-Lipschitz regular, for some constant $L_{\bp}=L(\bp,n)>0$:
        \[|\delta_{\bp}(p)- \delta_{\bp}(q)|   \le L_{\bp} \cdot d(p,q), \mm{ for } p,q \in  H \setminus \Sigma.\]
\end{itemize}
For further details and explanations we refer to Ch.\ref{sor},Def.1 and to Ch.\ref{hardy}, Rm.\ref{moti}, for the relation of (*) to classical (metric) Hardy inequalities.\\

 $\delta_{\bp}$ can be seen as a distance function to $\Sigma$ more naturally associated to $H$ than the metric distance function $dist(\cdot,\Sigma)$. Up to some gauging, $\delta_{\bp}(p)$ is a \emph{guessed} distance to $\Sigma$ derived from knowing $|A|$ locally around $p$. This gives $\delta_{\bp}$ a proper meaning even when $\Sigma \v$. We call  $\delta_{\bp}$ the \emph{${\bp}$-distance}.\\

\textbf{Skin Uniformity of $H \setminus \Sigma$} \quad In Ch.4 we establish a basic asymptotic regularity result,
the so-called skin uniformity, for the manifold  $H \setminus \Sigma$.\\

To describe this key property, we recall that in dimension $2$, the Riemann uniformization theorem shows that any domain $D \subset \R^2 = \C$, $D \neq \C$ or $\C^*$, admits a conformally equivalent complete hyperbolic metric.\\

 If we also want to ensure some finer details, like asymptotic homogeneity properties, we need to keep some control over the distortion while we approach $\p D$. It turns out that this precisely happens when $\p D$ is \emph{uniformly perfect}, cf.[K] for a broad discussion. This uniformity also ensures that the metric is Gromov hyperbolic, a hyperbolicity property oftentimes more appropriate for the asymptotic analysis of $D$ than classical hyperbolicity, cf.[An1].\\

Also one has a classification theory for Riemann surfaces due to Ahlfors and others in terms of properties of (positive) harmonic functions on the surface which in turn reflect the asymptotic behavior of the surface towards ideal boundaries, cf.[FK],[AS],[CC].\\

Due to geometric work, in particular of Gehring, Osgood [GO] and  Bonk, Heinonen, Koskela [BHK] and analytic results of Jersion, Kenig [JK] and Aikawa [Ai1]-[Ai3],we have
generalizations of such classically two-dimensional results for \emph{uniform domains} in $\R^n$, for any $n \ge 2$. \\

 A domain $D\varsubsetneq\R^n$ is called uniform provided any two points $p,q \in D$ can be
linked by a
 path $\gamma: [0,1] \ra D$, so that for some $c(D) \ge 1$, independent of $p$ and $q$:
\begin{itemize}
    \item \quad  $l(\gamma)  \le c \cdot d_D(p,q)$, for any two points $p,q \in D$
    \item \quad $l_{min}(\gamma(z)) \le c \cdot dist(z,\p D)$, for any $z \in \gamma_{p,q},$
\end{itemize}
Here $d_D$ is the restriction of the Euclidean distance to $D$. The second condition says that any point in $D$
 can be approached by a twisted cone within $D$, while the first condition says that the soul of this cone is a quasi-geodesic arc.\\

Examples of uniform domains include smoothly bounded and Lipschitz domains or the more general NTA-domains, cf.[JK], which played an important role in the development of the potential theory on irregular domains. What makes the broader concept of uniform domains remarkable is that uniformity of a Euclidean domain is \emph{equivalent} to each of the following two conditions, up to modest technical assumptions we omit here, cf.[BHK],[Ai2],[Ai3]:\\

Firstly, the quasi-hyperbolic metric, we get from conformally deform the Euclidean metric by $dist(\cdot,\p D)^{-2}$, is
 Gromov hyperbolic and its Euclidean boundary is homeomorphic to the Gromov boundary  and, secondly, uniformity is also equivalent to the validity of  boundary Harnack principles, for the Laplacian $\Delta$, relative $\p D$. This, in turn, implies that the Euclidean boundary equals the Martin boundary for $\Delta$.\\

The uniformity concept can also be considered for \emph{non-complete} metric spaces, with mild regularity properties, regarding the points added in their metric completion as their boundary. These \emph{uniform spaces} have a completely similar Gromov hyperbolization, cf.[BHK],[He]. \\

However, for general uniform spaces, there is no proper counterpart to the asymptotic analytic theory on Euclidean uniform domains since uniformity cannot detect fast degenerating geometries, like sharpening wrinkles, occurring while we pass to the boundary.\\

This is our first occasion to appeal to skin structures. We show that $H \setminus \Sigma$ is a \emph{skin uniform} space, a property stronger than uniformity of metric spaces.
It asserts that, for any $p,q \in H \setminus \Sigma$, there is a rectifiable path $\gamma: [a,b] \ra H \setminus \Sigma$, for some $a <b$,
so that for any given skin transform $\bp$, there is some $s_H \ge 1$ with \begin{itemize}
    \item \quad  $l(\gamma)  \le s_H \cdot d_{g_H}(p,q)$, for any two points $p,q \in H \setminus \Sigma$
    \item \quad $l_{min}(\gamma(z)) \le s_H \cdot \delta_{\bp}(z)$, for any $z \in \gamma_{p,q}.$
    \end{itemize}

We shall see that skin uniformity of $H \setminus \Sigma$ gives us the same degree of control over both the hyperbolicity properties and
the asymptotic analysis of elliptic operators on $H \setminus \Sigma$ as we have on \emph{uniform domains} in $\R^n$ and we also notice
\begin{itemize}
    \item \emph{skin uniformity $\Rightarrow$ uniformity}, since the Lipschitz continuity of $\delta_{\bp}$ readily shows that $\delta_{\bp} \le L
        \cdot dist(\cdot,\Sigma)$. Indeed, skin transforms are not only a means to express but also to \emph{prove} both the skin uniformity and the, also previously unknown, uniformity of  $H \setminus \Sigma$, from the isoperimetric inequality for area minimizers.
    \item The optimal \emph{skin uniformity constant} $s_H$ commutes with the convergence of
    sequences of minimal hypersurfaces, whereas the \emph{uniformity constant} $c_H$ does not. In particular, skin uniformity is a reasonable concept even when $\Sigma \v$ that matches the singular cases under degenerations.\\
\end{itemize}

\textbf{Gromov Hyperbolicity of $H \setminus \Sigma$} \quad Now we explain how skin uniformity gives us means to literally unfold the geometry and analysis on $H \setminus \Sigma$.\\

In [L1],Ch.2, we first establish the relation between skin uniformity and hyperbolic structures on $H \setminus \Sigma$. For this we conformally deform $g_H$ by $\bp^2$. The
resulting metric
 \begin{equation}\label{sme}
 d_{\bp}(x,y) := \inf \Bigl  \{\int_\gamma  \bp \, \, \Big| \, \gamma   \subset  H \setminus \Sigma\mbox{ rectifiable arc joining }  x \mbox{ and } y  \Bigr \}
 \end{equation}
for any two points $x,y \in H \setminus \Sigma$, is called the \emph{skin metric} on $H \setminus \Sigma$. This metric is a counterpart to the \emph{quasi-hyperbolic metric} $k_{H \setminus \Sigma}$
 defined by
\begin{equation}\label{qh}
k_{H \setminus \Sigma}(x,y) := \inf \Bigl  \{\int_\gamma 1/dist(\cdot, \Sigma)   \, \, \Big| \, \gamma   \subset H \setminus \Sigma \mbox{ rectifiable arc joining }  x \mbox{ and } y  \Bigr \}
\end{equation}
 which is basic in quasi-conformal geometry, cf.[Ai1],[BHK],[K].\\

The skin uniformity of $H \setminus \Sigma$ reveals a combination of properties of $d_{\bp}$ not available from the uniformity of $H \setminus \Sigma$ and not valid for the quasi-hyperbolic metric $k_{H \setminus \Sigma}$:
\begin{itemize}
    \item $(H \setminus \Sigma,d_{\bp})$  is a complete, \emph{Gromov hyperbolic space with bounded geometry} and it admits quasi-isometric Whitney type regularizations to
        smooth manifolds $(H \setminus \Sigma, d_{\bp^*})$. We call these spaces \emph{hyperbolic unfoldings} of $H \setminus \Sigma$ of the conformally equivalent skin uniform space $(H \setminus \Sigma,g_H)$.
    \item  The Gromov boundary $\p_G (H \setminus \Sigma,d_{\bp})$ of $(H \setminus \Sigma,d_{\bp})$ equals $\Sigma \subset H$.
    \item The skin metrics commute with the convergence of minimal hypersurfaces. In particular, $d_{\bp}$ commutes with blow-ups leading e.g. to tangent cones.
     \item Also when $\Sigma_H \v$, the metric space $(H,d_{\bp})$  is still well-defined. In general, it is homeomorphic to $(H,g_H)$, only when $H$ is totally geodesic it degenerates to a single point, whereas $(H,k_{H \setminus \Sigma})$ always is the one-point space, since $dist(\cdot, \Sigma) = \infty$.\\
\end{itemize}

\textbf{Potential theory and Asymptotic Analysis on $H \setminus \Sigma$} \quad It is the alliance of hyperbolicity and bounded geometry that makes the unfolded versions of
$(H \setminus \Sigma,g_H)$ analytically approachable.\\

 In [L1],Ch.3, we consider the \emph{potential theory} for a large class of elliptic operators on the smooth manifold $(H \setminus \Sigma, d_{\bp^*})$. The point is that the potential theory, for reasonably chosen operators,  on Gromov hyperbolic manifolds with bounded geometry is pleasantly transparent. This had been discovered by Ancona, cf. [An1], [An2].\\

We transfer this analysis from the hyperbolic unfolding $(H \setminus \Sigma, d_{\bp^*})$ back to the original hypersurface $H \setminus \Sigma \subset M$. One of the outcomes is a surprisingly simple asymptotic analysis, towards $\Sigma$, for any \emph{skin adapted operator}. A linear second order elliptic operator $L$ is skin adapted provided
\begin{itemize}
     \item $L$ does not degenerate faster than $\bp^2$ when we approach $\Sigma$.
     \item $L$ satisfies a weak coercivity condition: we require that there is a subsolution $s
         >0$ of $L \, u=0$ and some $\ve >0$, so that: $L \, s \ge \ve \cdot \bp^2 \cdot s.$\end{itemize}

To mention some concrete results, we first consider the associated Martin theory.  In general, Martin boundaries, in particular of wrinkled spaces with irregular boundaries, like $H \setminus \Sigma$, are difficult to understand. But the skin uniformity, in the guise of hyperbolic unfoldings, shows that the \emph{Martin boundary $\p_M (H \setminus \Sigma,L)$ of  $H \setminus \Sigma$ for skin adapted $L$} simply equals the Gromov boundary of its unfolded version and, thus, it can be identified with the singular set $\Sigma$:
\[\p_M (H \setminus \Sigma,L) \cong \p_G (H \setminus \Sigma,d_{\bp}) \cong \Sigma\]
where $ \cong $ means homeomorphisms. Hence, for any positive function $u$ on $H \setminus \Sigma$ we have: $u$ solves $L \, v = 0$ if and only if $u$ admits a integral representation in terms of a (uniquely determined) finite Radon measure $\mu=\mu_u$ on $\Sigma$: \[u(x) = \int_{\Sigma} k_L(x;y) \, d \mu(y),\] where $k_L(x;y)$ denotes the Martin kernel,  generalizing the Poisson kernel in
the integral representation of harmonic functions on the unit disc of Herglotz, cf.[BJ],1.7.2.\\

We also get a fine control for the asymptotic behavior of positive solutions of $L \, v=0$ along $\Sigma$. For instance, the quotient of any two positive solutions  of $L \, v=0$ on $H \setminus \Sigma$, with minimal growth towards a given open subset $A \subset \Sigma$, admits a \emph{continuous extension} to $A$, whereas, the individual solutions usually diverge towards $A$.\\

A simple but important observation is that skin adaptedness persists under blow-ups: induced operators on tangent cones are again skin adapted. This makes the analysis of  skin adapted operators amenable to inductive dimensional reduction schemes by means of blow-ups of singularities. This is a key detail in our applications of skin structural concepts  to scalar curvature geometry, in [L2], we shall illustrate below. \\

\textbf{Hardy Inequalities} \, In the results we described so far we exploited the axioms for skin transforms, actually several times, with the exception of the \emph{Hardy type inequality} \[\int_H|\nabla f|^2  + |A_H|^2 \cdot f^2 dA \ge \tau \cdot \int_H \bp_H^2\cdot f^2 dA, \mm{ for some } \tau_H >0.\]
Its role is to couple this theory to classical problems. To this end, we observe that the class of skin adapted operators is rich but, a priori, it could be unrelated to the relevant geometric analysis on $(H \setminus \Sigma,g_H)$ and it depends on the chosen skin transform $\bp$.\\

However, the Hardy inequality shows that, independent of the chosen skin structure, many interesting elliptic operators are indeed skin adapted.  For now, we mention only the one example of particular importance in scalar curvature geometry:
 \begin{itemize}
 \item The \emph{conformal Laplacian} $L_H=-\Delta  +\frac{n-2}{4 (n-1)} \cdot scal_H,$ is a skin adapted operator on $H \setminus \Sigma$, if and, in general, only if $scal_M \ge 0.$
\end{itemize}
Note that, in turn, the statements of the results in the potential theory for $L_H$ on $H$, like the description of the Martin boundary, do not involve skin structures. In other words, through the use of skin structures we can approach new, although classically expressible results for area minimizers. Other examples, we already mentioned, are  the uniformity of $H \setminus \Sigma$ and, thereby, the hyperbolicity of $k_{H \setminus \Sigma}$.\\

\textbf{Scal$>$0-Heredity with Surgery} \, Another application of skin structures is the incorporation of singular area minimizing hypersurfaces in the study of scalar curvature constraints in [L2]. \\

We recall that when $M^{n+1}$ carries a $scal >0$-metric, any area minimizer $H^n \subset M^{n+1}$ can be conformally deformed to a $scal >0$-space. A simple variational argument reveals that the first eigenvalue $\lambda_1$ of the conformal Laplacian $L_H$ is positive. Then the con´formal defomation by eigenfunction $\phi >0$ gives a $scal >0$-metric $\phi ^{4/n-2} \cdot g_H$:
\begin{equation} \label{y01} scal(\phi ^{4/n-2} \cdot g_H) \cdot \phi ^{\frac{n+2}{n-2}} = - \triangle \phi   + \frac{n-2}{4 (n-1)}  \cdot scal_H \cdot  \phi = \lambda_1 \cdot \phi.\end{equation}

 This remarkable $scal >0$\emph{-heredity property}, observed by Schoen and Yau in the late 70ties, suggested an approach to scalar curvature geometry by means of an inductive dimensional descent along towers of nested hypersurfaces until one reaches a well-understood lower dimensional space, cf.[SY], [GL]. However, in these inductive arguments singularities of $H$ cause accumulating problems: they appear in dimensions $> 7$ and
have made the implementation of this idea intractable in higher dimensions.\\

In [L2], we introduce the broader strategy of a \emph{$scal >0$-heredity with surgery} to replace the singularities for regular ends and closures, without sacrificing the $scal >0$-condition. Thereby, we can study scalar curvature problems in arbitrary dimensions.\\

The underlying core result of [L2] is the existence of conformal deformations of $H \setminus \Sigma$ to a $scal>0$-space with a \emph{spreading open end} along $\Sigma$ retaining area minimizers within $H \setminus \Sigma$ from approaching $\Sigma$. To them, $H \setminus \Sigma$ appears to be a compact smooth $scal>0$-manifold, whereas $\Sigma$ may be thought as lying beyond an \emph{inner horizon}.\\

To indicate how these conformal deformations can be defined we first observe that the skin adpatedness of $L_H$ shows that  $L_{H,\lambda}:= L_H - \lambda  \cdot  \bp^2 \cdot id$ is \emph{again skin adapted}, when $\lambda >0$ is small enough. Thus we have (many) positive solutions for the eigenvalue equation for the \emph{skin conformal Laplacian} $\delta^2_{\bp} \cdot L_H$ with eigenvalue $\lambda$ and get
\begin{equation} \label{y1} scal(v^{4/n-2} \cdot g_H) \cdot v^{\frac{n+2}{n-2}} = - \triangle v  + \frac{n-2}{4 (n-1)}  \cdot scal_H \cdot  v = \lambda  \cdot  \bp^2  \cdot v \, \mm{ on } H \setminus \Sigma.\end{equation}
This suggests to exploit the Martin theory for
 $L_{H,\lambda}$ to find a solution $u>0$ of (\ref{y1}) so that $(H \setminus \Sigma, u^{4/n-2} \cdot g_H)$ has $scal >0$ and it also contains a horizon shielding $\Sigma$. Note that, since $\bp^2 $ is locally Lipschitz, any such solution is $C^{2,\alpha}$-regular, for $\alpha \in (0,1)$.\\

Heuristically, the strategy is to define an \emph{equidistributed} positive Radon measure $\mu$ on $\Sigma$ and to use growth estimates for the Martin kernel to show that
\begin{equation}\label{mic}
u_\mu(x) := \int_{\Sigma} k_{L_{H,\lambda}}(x;y) \, d \mu(y)
\end{equation}
is a solution that induces an inner horizon shielding $\Sigma$. However, inhomogeneities of $H$ along $\Sigma$ make it hard to describe or evaluate equidistributed measures on $\Sigma$ directly.\\

But we can bypass this issue. In simple terms, we reduce the problem to the case of couples $\Sigma^* \subset H^*$,  where $\Sigma^*$ and $H^*$ are simplicial complexes. The simplices of $\Sigma^*$ are pieces of symmetry axes of  iterated tangent cones $C^n$ of the original $\Sigma \subset H$, and we define a compatible system of equidistributed measures for these (pieces of) cones.\\

This use skin structures in several ways: to control the again skin adapted operators $L_{C,\lambda}$, to geometrically exploit the scaling invariance of the $\bp$-weighted eigenvalue $\lambda$ in (\ref{y1}) and to assemble the systems of equidistributed measures for $\Sigma^* \subset H^*$. And, the actual reduction to the approximating couples $\Sigma^* \subset H^*$ uses again the skin uniformity.\\

\textbf{Acknowledgements} \,  The author thanks Misha Gromov, Jan-Mark Iniotakis, Matthias Kemper and Frederik Witt for many helpful suggestions on how to improve the exposition of these papers. Part of this work, mainly the applications in scalar curvature geometry in [L2], had been announced in [L3].\\

\subsubsection{Results in the Present Paper}\label{sor}
\bigskip

In this paper we introduce skin structures and apply them to derive basic structural results for area minimizing hypersurfaces near their singularities. Our main focus will be on
the uniformity properties of the regular portion of such spaces.\\

\textbf{Setup and Notations}\quad 1. $H^n \subset M^{n+1}$ is a connected locally area minimizing hypersurface without boundary in some $C^\infty$-smooth Riemannian
$n+1$-manifold $(M,g_M)$. $g_H$ is the induced metric on the smooth portion of $H$.\\

In technical terms, $H$ is a locally mass minimizing, integer multiplicity rectifiable current of dimension $n$ without boundary.
The partial regularity theory for these minimizers says that $H$ is a smooth hypersurface except for some singular set $\Sigma_H$ of Hausdorff-dimension $\le n-7$.\\

 2. We consider the following classes of complete area minimizing hypersurfaces.
\begin{description}
    \item[${\cal{H}}^c_n$]:    $H \subset M$ is a compact and connected hypersurface without boundary.
    \item[${\cal{H}}^{\R}_n$]:  $(M,g_M) = (\R^{n+1},g_{\R^{n+1}})$, $H$ is an oriented boundary of some open set $A \subset \R^{n+1}$  and $0 \in H$. Thus, $H$ is non-compact
        and complete.
        \item[${\cal{SH}}^{\R}_n$]:      ${\cal{SH}}^{\R}_n \subset  {\cal{H}}^{\R}_n$ is the subset of those hypersurfaces singular, at least, in $0$.
\end{description}
The main class of hypersurfaces we study in this paper is given by
\[{\cal{H}}_n:= {\cal{H}}^c_n \cup {\cal{H}}^{\R}_n \, \mm{ and }\, {\cal{H}} :=\bigcup_{n \ge 1} {\cal{H}}_n.\]
${\cal{H}}_n$ is closed under blow-ups. That is, the limit of converging subsequences under scaling by a diverging sequence of real numbers belongs to ${\cal{H}}^{\R}_n$.\\

3. In dimensional reduction arguments we also consider more particular classes of minimal hypersurfaces:
\begin{description}
    \item[$\mathcal{C}_{n}$]:   $\mathcal{C}_{n} \subset {\cal{H}}^{\R}_n$ is the space of area minimizing $n$-cones in $\R^{n+1}$ with tip in $0$.
    \item[$\mathcal{SC}_{n}$]:    $\mathcal{SC}_n \subset\mathcal{C}_n$ is the subset of cones singular, at least, in $0$
    \item[$\mathcal{K}_{n-1}$]: For any area minimizing cone $C \subset \R^{n+1}$ with tip $0$, we get the non-minimizing minimal hypersurface $S_C$ in the unit sphere
\[S_C:= \p B_1(0) \cap C \subset S^n \subset  \R^{n+1} \mm{ and set } {\cal{K}}_{n-1}:= \{ S_C\,| \, C \in {\mathcal{C}_{n}}\},\] where $\mathcal{C}_{n}$ is the space of area
minimizing $n$-cones in $\R^{n+1}$ with tip in $0$. We write ${\cal{K}}= \bigcup_{n \ge 1} {\cal{K}}_{n-1}$, for the space of all such hypersurfaces $S_C$.
\end{description}

4. $A = A_H$ denotes the second fundamental form of $H \subset M$, $|A|$ is the norm or length of $A$ and $\RH$ is the set of positive regular values of $|A|$. We
    denote $|A|$-level sets, for $c >0$, by  $\B_c := |A|^{-1}(c)  \subset H$.\\

For readers not familiar with geometric measure theory we
have included the necessary background and references to the literature in the appendix.\\\\

\textbf{Skin Transforms} \quad  We start with an axiomatic description of skin transforms. Intuitively, any skin transform is a distinguished averaging procedure for the norm of
second fundamental form defined from equalizing procedures transforming simultaneously all $|A|$-level sets $\B_c$ into better controlled wraps $\M_c$,  the so-called $|A|$-skins.\\

When $|A| \equiv 0$, that is, when $H \subset M$ is totally geodesic, it is not hard to see that $\Sigma \v$, that is, $H$ is a smooth manifold. As usual, saying $H$ is not totally
geodesic, means there is at least one point where $|A| \neq 0$.\\

\textbf{Definition 1} \textbf{(Skin Transforms)} {\itshape \quad An assignment $\bp$ defined on ${\cal{H}}$ that gives us
 for any $H \in {\cal{H}}$
a measurable and non-negative function $\bp_H$ defined on $H \setminus \Sigma_H$ is called a \textbf{skin transform} when the following holds
 \begin{description}
    \item[(S1)] \emph{\textbf{Trivial Gauge}} \,  $\bp_H \equiv 0$, if $H \subset M$ is totally geodesic.
    \item[(S2)] \emph{\textbf{Skin Property}} \,  When $H$ is not totally geodesic, then $\bp_H$ is  strictly positive. Its level sets $\M_c:= \bp_H^{-1}(c)$, called the
        $\mathbf{|A|}$\textbf{-skins}, surround the $\B_c:= |A|^{-1}(c)  \subset H$, for $c>0$:
    $$\bp_H \ge |A_H|\,\, \mm{ and  }\,\, \M_c \cap \Sigma \v, $$ and, like  $|A|$, $\bp$ anticommutes with scalings, that is, $\bp_{\lambda \cdot H} \equiv \lambda^{-1} \cdot  \bp_{H}$, for any  $\lambda >0$.
    \item[(S3)] \emph{\textbf{Hardy Inequality}}  \,  There is a positive constant $\tau = \tau(\bp,H) \in (0,1)$, so that for any $f \in C^\infty(H \setminus \Sigma,\R)$
        compactly supported in $H \setminus \Sigma$
    $$\int_H|\nabla f|^2  + |A|^2 \cdot f^2 dA \ge \tau \cdot \int_H \bp^2\cdot f^2 dA.$$  $\tau$ is called the \textbf{tightness} of $\bp$ on $H$.
    \item[(S4)] \emph{\textbf{Lipschitz regularity}} \, When $H$ is not totally geodesic, we set $\delta_{\bp}:=1/\bp$,  and call this quantity $\mathbf{\bp}$\textbf{-distance}$^*$. It is
        $L_{\bp}$-Lipschitz regular, for some constant $L_{\bp}=L(\bp,n)>0$:
        \[|\delta_{\bp}(p)- \delta_{\bp}(q)|   \le L_{\bp} \cdot d(p,q), \mm{ for } p,q \in  H \setminus \Sigma, \mm{ for any } H \in {\cal{H}}_n. \]
    \item[(S5)]  \emph{\textbf{Naturality}}$^{**}$ \, The assignment $\bp$ commutes with the convergence of sequences of area minimizers $H_i \in {\cal{H}}_n$, $i \ge 1$ to a
        limit space $H \in {\cal{H}}_n$, so that for any $\alpha \in (0,1)$: \, $\bp_{H_i} \overset{C^\alpha}  \longrightarrow {\bp_{H}}.$\\
 \end{description}}
\textbf{Remark 1} \quad  $^*$When $H$ is totally geodesic, and thus $\bp \equiv 0$, we set accordingly $\delta_{\bp}\equiv +\infty$. Also, for totally geodesic hypersurfaces the assertions in this paper are largely empty or may be taken as conventions.\\

$^{**}$Technical details, like the precise meaning of this commutativity or the $C^\alpha$-convergence are discussed in Ch.\ref{rana}.
 \qed

A basic model of skin transforms are interpolations between $|A|$ and $1/dist(x,\Sigma)$.\\

 \textbf{Theorem 1} \textbf{(Metric Skin Transforms)} {\itshape \, There is a family of skin transforms $\bp_\alpha$, for any $\alpha >0$, with the following properties:\\

The \sks\ $\M_c$ of $\bp_\alpha$ bound the outer $\alpha/c$-distance collar of the level sets $\B_c$ of $|A|$ in $H$  and we have:
 \begin{itemize}
    \item \quad   $\bp_{\alpha}(x) \ra |A|(x)\mm{ in } L_{loc}^\infty,\mm{ for }\alpha \ra 0, \mm{ on } H \setminus \Sigma$
        \item \quad   $1/\alpha \cdot \bp_{\alpha}(x) \ra 1/ dist(x,\Sigma)\mm{ in } L_{loc}^\infty,\mm{ for }\alpha \ra \infty, \mm{ on } H \setminus \Sigma$.
  \end{itemize}}
  The outer distance collar of $\B_c$ means the distance neighborhood of $\{x \in H \setminus \Sigma \,|\, |A|(x) \ge c\}$.\\

 \textbf{Remark 2} \quad The two limits $|A|$ and $c/ dist(\cdot,\Sigma)$ are no longer skin transforms. In general $\B_c \cap \Sigma \n$ and there is no Lipschitz bound for
the $|A|$-distance. Thus $|A|$ violates (S2) and (S4). In turn, there is no constant $c>0$, so that $c/ dist(\cdot,\Sigma) \ge |A|$, as required for (S2), and there is no proper
correlation between the singularities of converging sequence of hypersurfaces and their limit. Thus $1/ dist(\cdot,\Sigma)$ also violates the naturality axiom (S5).\qed

The space of all skin transforms is contractible and each skin transform comes with entourage of others e.g. in dimensional reduction processes. Therefore, we focus on results valid for all skin transforms. Later on, we employ skin structures, without reference to their pedigree, exclusively based on the skin axioms of definition 1.\\

For the remainder of this introduction we may use any skin transform $\bp$. The results change only by global constants depending on the chosen $\bp$.\\

\textbf{Accessibility of $\mathbf{\Sigma}$} \quad  Now we turn to the central skin structural concept in this series of papers. The description of the quality of the
 approachability of $\Sigma$ from within $H \setminus \Sigma$. To this end we use the notion of uniform spaces cf.[BHK], [V] and [He]  as a
starting point.\\

\textbf{Definition 2} \textbf{(Uniform Spaces)}\label{ud} \, \emph{A  \textbf{uniform space}, more precisely a c-uniform space, is a locally compact, \textbf{non-complete},
locally complete, rectifiably connected metric space  $(X,d)$  so that there is a constant $c \ge 1$ such that any two points can be joined by a
c-uniform curve.}\\

\emph{\noindent A \textbf{c-uniform} curve joining $p,q \in X$ is a rectifiable path $\gamma: [a,b] \ra X$, for some $a <b$, from $p$ to $q$ so that $\gamma$ satisfies the
following two axioms
\begin{itemize}
    \item \emph{\textbf{Quasi-Geodesic}:} \quad \quad  $l(\gamma)  \le c \cdot  d(p,q),$
    \item  \emph{\textbf{Twisted Double Cones}:}\quad  $l_{min}(\gamma_{p,q}(z)) \le c \cdot dist(z,\p X),$
\end{itemize}
 for any  $z \in \gamma_{p,q}.$ $l_{min}(\gamma_{p,q}(z))$ denotes the minimum of the lengths of the subarcs of $\gamma_{p,q}$ from $p$ to $z$ and from $q$ to $z$.}\\

In this terminology, we think of $\Sigma_H$ as a boundary of a singular area minimizer $H$, that is, we set $X = H \setminus \Sigma$, $\p X = \Sigma$ and  $\overline{X} = H$. We get a stronger version of this purely metric type of uniformity which also naturally extends to the case of regular $H$.\\

 \textbf{Theorem 2} \textbf{(Skin Uniformity of  $\mathbf{H \setminus \Sigma}$)} {\itshape \, For any hypersurface $H \in {\cal{H}}_n$ with singular set $\Sigma_H$,
 which may also be empty, we have

 \begin{enumerate}
 \item$H \setminus \Sigma$ and $H$ are \textbf{rectifiably connected}. In particular, any compact $H \in {\cal{H}}^c_n$ has a finite intrinsic diameter: $diam_{g_H}H   <
     \infty.$
\item $H \setminus \Sigma$  is a $c$-\textbf{skin uniform space}, for some $c >0$. That is, any pair $p,q \in H \setminus \Sigma$ can be joined by a \textbf{c-skin uniform
    curve} in $H \setminus \Sigma$, i.e. a rectifiable path $\gamma_{p,q}: [a,b] \ra H \setminus \Sigma$, for some $a <b$, with $\gamma_{p,q}(a)=p$,
    $\gamma_{p,q}(b)=q$, so that the following  skin uniformity conditions hold
    \begin{itemize}
    \item \emph{\textbf{Quasi-Geodesic}:} \quad \quad  $l(\gamma)  \le c \cdot  d(p,q),$
    \item  \emph{\textbf{Twisted Double Skin Cones}:}\quad $l_{min}(\gamma_{p,q}(z)) \le c \cdot \delta_{\bp}(z)$,
\end{itemize}
for any $z \in \gamma_{p,q}.$
\end{enumerate}}
\smallskip
 \textbf{Remark 3} \quad 1. The parameter $c=c(H)$ is a natural invariant and thus relevant for both, \emph{regular} and \emph{compact}  $H$: for a sequence
$H_i \in {\cal{H}}_n$ flat norm converging to a minimizer $H \in {\cal{H}}_n$, that the $H_i$ and $H$ are c-skin uniform, for the  same $c>0$. For totally geodesic $H$, we have $\delta_{\bp}=+\infty$. Then, the twisted cone condition holds trivially.\\

2. In turn, for  $\Sigma \n$, we will see that the skin uniformity condition still holds for any $p,q \in H$, that is, we can also link singular points by c-skin uniform
    curves, supported in $H \setminus \Sigma$, except for the endpoints.\\

3. Skin uniformity implies uniformity. For $\Sigma \n$, $H \setminus \Sigma$ is also a uniform space. Indeed, $\bp \ge |A|$ in (S2) and the Lipschitz condition in (S4) imply that for any $x \in H \setminus \Sigma$, cf. \ref{locsking}: \, $\delta_{\bp}(x) \le L  \cdot dist(x,\Sigma)$.\qed

For area minimizers in $\R^{n+1}$, the latter result can be sharpened and we have bounds for the skin uniformity parameter depending only on the dimension.
This is useful since this class of hypersurfaces comprises that of all blow-up limits.\\

\textbf{Theorem 3} \textbf{(Area Minimizers in $\mathbf{\R^{n+1}}$ and Blow-Ups)}    \emph{\begin{itemize}
    \item  Any oriented minimal boundary $H \subset \R^{n+1}$  is $c_n$-skin uniform, for some constant $c_n >0$, depending only on the dimension.
      \item For $C \in {\cal{C}}_n$,
 $S_C \setminus \sigma_C=\p B_1(0) \cap C \setminus \sigma_C$  is a $c^*_n$-skin uniform space relative $\bp:=\bp_C|_{S_C}$,
  for a constant $c^*_n$, depending only on the dimension.
     \item Let $H$ be some complete area minimizing hypersurface, so that $H \setminus \Sigma$ is a $c$-skin uniform space, for some $c >0$.
 Then, for any blow-up limit $F$, $F \setminus \Sigma_F$ is also $c$-skin uniform. In particular, this applies to all tangent cones of $H$.
\end{itemize}}

Note that the inheritance result merely assumes the completeness of $H$. Thus it can be applied iteratively in cone reduction arguments.\\

In turn, in general, for \emph{non-compact} $H$ in a non-Euclidean ambient space,  $H \setminus \Sigma$ is not (skin) uniform. We can neither expect a uniform degree of
non-tangential accessibility of points of non-compact singular sets nor for points at infinity. In this case, similar problems cause a failure of the
Hardy inequality (S3). \\

As explained in the overview Ch.\ref{gint},  the skin uniformity is a vital input to establish the conformal  hyperbolization and to understand the potential theory of $H \setminus \Sigma$ in [L1]. There are skin structural substructures which further amplify these links.\\

\textbf{Theorem 4}  \textbf{(Skin Uniform Domains)} \quad  \emph{Let $H \in {\cal{H}}$ be a $c$-skin uniform space for some $c(H) >0$.
Then there are some $\iota(c),\kappa(c) > 1$ so that for sufficiency small $a >0$, there is a \textbf{skin uniform domain} $\su= \su(a) \subset H \setminus \Sigma$:
\begin{itemize}
  \item    \quad  $\E(\iota \cdot a) \subset \su(a) \subset \E(a)$, where $\E(a):= \{x \in H \setminus \Sigma \,|\, \delta_{\bp}(x) \ge a \},$
\end{itemize}
and any two points $p,q \in \su$ can be linked by an arc $\gamma_{p,q} \subset \su$ with
\begin{itemize}
    \item \quad $l(\gamma)  \le \kappa \cdot d_H(p,q)$,
    \item \quad $l_{min}(\gamma(z)) \le \kappa \cdot \min\{L_{\bp} \cdot dist(z,\p \su),\delta_{\bp}(z)\}$, for any $z \in \gamma_{p,q}$.
\end{itemize}
The class of domains $\su(a)$ is naturally associated to $H$ and their existence is (obviously) equivalent to the skin uniformity of $H$.} \\

\textbf{Organization of the paper} \,  In \textbf{Chapter 2} we describe the naturality concept in more detail and construct skin transforms. We  verify the axioms, that is,
we establish Theorem 1  except for the Hardy inequality which is postponed to  Chapter 3.\\
 Then, in \textbf{Chapter 3} we construct skin adapted covers by balls of radius comparable to their $\bp$-distance. We use them in different places: to derive the Hardy inequality,  Whitney type smoothings of skin transforms, we need  in Martin theory in [L1], and finally to define the skin uniform domains.\\
 In \textbf{Chapter 4} we discuss non-tangential accessibility concepts, prove the skin uniformity of regular part of singular area minimizers and get Theorems 2, 3 and 4. \\

In this paper, area minimizing hypersurfaces are locally mass minimizing integer multiplicity rectifiable currents of dimension $n$.  This is a sufficiently large class of
hypersurface to be able to find a mass minimizer in each $n$-dimensional homology class of $M^{n+1}$ when we identify
singular homology and integral current homology. (For smooth manifolds these homology theories are isomorphic).\\

In the above results we always assumed that $H$ is connected. If we need to consider an area minimizer with several components we apply the result componentwise.
This happens, for instance, when we consider the mass minimizer that realizes a general homology class.\\

These currents are sufficiently regular to locally \textbf{decompose} them into (locally disjoint) oriented \textbf{minimal boundaries} of open sets,  cf. \ref{dic} and \ref{dicc}
and [F],4.5.17,  [Si2],Ch.37 and [Si3]. The case of oriented boundaries is typically treated by means of functions of bounded variations, cf.[AFP] and [Gi], and we occasionally refer
to these methods as the \textbf{BV-techniques}.\\

These local decompositions allow us to use extrinsic arguments for oriented boundaries not available in the current setup and also to eliminate formal concerns. We  benefit from this decomposition in the study of intrinsic concepts like that of skin structures.\\

In most places in this paper we may therefore assume, without loss of generality, that $H$ can locally be written as an oriented boundary of an open set in $M$.
We have included an \textbf{appendix} where we briefly outline these and other relevant pieces of geometric measure theory.
Only in chapter 4, where we use more specific results for oriented boundaries, we discuss further details of this reduction.\\

\setcounter{section}{2}
\renewcommand{\thesubsection}{\thesection}
\subsection{Skin Structures and Naturality}
\bigskip

In this chapter we describe approaches to define skin transforms. One important property of these notions, which sets them apart from metric distance concepts, is their naturality. We start with a discussion of what this naturality means.

\subsubsection{Natural Structures}\label{rana}

\bigskip

We define the notions of $\D$-maps and of natural structures. They come into play when we use compactness theorems for area minimizers and
try to extend these compactness argument to structures assigned to hypersurfaces.\\

Natural structures are those continuously depending on the underlying area minimizer which amounts to an \emph{interchangeability}
 of two limits:  for any converging sequences of underlying
spaces the attached structures converge to the one we assigned to the limit space.\\

 Trivial examples are curvatures or classical elliptic operators, but since we restrict to the rather narrow class of area minimizers, such structures can also be \emph{non-locally} determined: their value, in a given point, may depend on details also from remote parts of the underlying space.\\

\textbf{Prototypical Setup} \quad There are two levels of compactness results we use for area minimizers, the flat norm (sub)convergence of sequences with common intrinsic area and extrinsic spatial bounds. Then, from the regularity theory of these minimizers, we get a compact smooth (sub)convergence in those places where the limit minimizer is smooth, cf. the appendix for a summary of this theory and some references.\\

Having this DeGiorgi-Allard theory at our disposal, we now describe the typical situations we encounter later on and explain the recurrent pattern of how to exploit these compactness results for area minimizers.\\

As already announced in the introduction, we generally use the local decomposition of rectifiable currents into oriented minimal boundaries of open sets in the ambient space described in the appendix, part IV, to restrict the presentation to the less technical case of oriented minimal boundaries.  \\

  We start with a sequence of area minimizers $H^n_i$ locally converging along with their ambience $M^{n+1}_i$ to a limit hypersurface within a limit ambient space, where the $M^{n+1}_i$ are complete smooth Riemannian manifolds.
The most important cases are $M=\R^{n+1}$ and $M_i= i\cdot M$, or $M_i=M$, for all $i \ge 0$,  for a fixed manifold $M$.\\

In detail, we consider the mentioned two level convergence relative some fixed base points  $p_i \in H_i \subset M^{n+1}_i$:\\

\textbf{Ambient Level}  \,  The $M^{n+1}_i$ compactly $C^k$-converge, for some $k \ge 5$, to a limit manifold $M^n$, so that $p_i \ra p \in M^n$, for $i \ra \infty$.   This says, for any given $R
         >0$ and sufficiently large $i$,  we have diffeomorphisms
\[\Psi_i : B_R(p_i) \ra B_R(p), \mm{ so that } |(\Psi_i)_*(g_{M_i}) - g_{M}|_{C^k} \ra 0 \mm{ on } B_R(p).\]

\textbf{Minimizer Level}  \,  The $H_i$ subconverge  (that is a
 subsequence converges) to a limit area minimizer $H \subset M$, with $p \in H$. This means for any $R >0$, we get a subconvergence within $M$
\begin{equation}\label{fl}
\Psi_i (B_R(p_i) \cap H_i) \ra B_R(p) \cap H.
\end{equation}
in \emph{flat norm}, defined in appendix III (\ref{fll}), which, in this case, can also be expressed in terms of \emph{Hausdorff distance}.

\begin{example} \label{ex1} \quad  Common examples for the preceding setup appear in inductive cone reductions and when we consider minimal hypersurfaces within symmetric spaces. \\

1.  When we scale around a fixed singular point $p \in \Sigma \subset H_0 \subset M_0$, we set $M_i:= \tau_i  \cdot M_0$ and $H_i:=\tau_i \cdot H_0$, for some sequence $\tau_i \ra \infty$, for $i \ra \infty$. Then $M_i$
converges compactly to $\R^{n+1}$ and we find a local flat norm subconvergence  to a limit space $H \subset \R^{n+1}$. It is a standard fact,  that $H$ is an
area minimizing cone, a so-called tangent cone,
cf. appendix, Ch.10.II and [Gi],9.3.\\

2. More generally, when we do \emph{not} fix the base point, the subconvergence under blow-ups still leads to a limit hypersurface $H \subset
\R^{n+1}$. It is a complete area minimizer which is an oriented boundary in $\R^{n+1}$. But, a priori, it need not to be a cone.\\

3. Another important example are sequences $H_i = S_{C_i}= C_i \cap S^{n+1} \in \mathcal{K}_n$ for some area minimizing cones $C_i \in \mathcal{C}_{n+1}$. Here we can choose  $R= diam \, S^{n+1} $, this means global flat norm convergence.\\

In this case, $H_i \in \mathcal{K}_n$  are critical points of the Area-functional, but they are not area minimizing. But they are so-called almost minimizers cf.\ref{riem} and share the regularity theory of area minimizers, cf.Appendix II. \qed
\end{example}

\begin{remark}\label{tech} \quad
Technically, when $\Psi_i$ is not an isometry or at least a global scaling, the $\Psi_i (B_R(p_i) \cap H_i)$ are \emph{almost minimizers} in $M$, cf. Ch.\ref{app}.II. Using Euclidean charts we can also locally view $H \subset
M$ as an almost minimizer in $\R^n$. The analysis of the convergence of the $H_i$ is properly covered from Allard regularity theory and, by similar arguments, from the regularity theory of almost minimizers in $\R^n$, a natural extension of  the theory for area
minimizers.\qed
\end{remark}

\textbf{Smooth Approximation} \quad Now, on the minimizer level in this setup,  we assume that $B_R(p) \cap H$ is \emph{smooth}. Then the regularity theory says that  the flat norm convergence of $\Psi_i
(B_R(p_i) \cap H_i)$ implies that, for sufficiently large $i$, the $B_R(p_i) \cap H_i$ will also be smooth
(for possibly slightly shrinked radii)  and the flat norm convergence induces a $C^k$-convergence in the following sense.\\

We denote the normal bundle of $B_R(p) \cap H$ by $\nu$. Then, for $i$ large enough, the $\Psi_i (B_R(p_i) \cap H_i)$ become local $C^k$-sections
$\Gamma_i$ of $\nu$ from using Fermi coordinates. In other words the regularity theory gives us canonical identifying
diffeomorphisms
\[\Gamma_i : B_R(p) \cap H \ra \Psi_i (B_R(p_i) \cap H_i).\]
Moreover, the flat norm convergence implies a $C^k$-convergence of the $\Gamma_i$ to the \emph{zero section} $\equiv B_R(p) \cap H$. Actually, when the ambient remains fixed and $\Psi_i =id_{M}$, the DeGiorgi
theory, cf. appendix II, gives a $C^1$ convergence which, by standard elliptic theory, can be upgraded to $C^{m, \alpha}$, for any $m \ge 0, \alpha \in (0,1)$. Composed with the $\Psi_i$, this gives a $C^k$-convergence.\\

This way we get, for large $i$, a canonical local diffeomorphism between $B_R(p) \cap H$ and $B_R(p_i) \cap H_i$. We use the following mnemonic to describe this particular map.

\begin{definition} \label{idmap} \emph{\textbf{($\D$\textbf{-map}) }}\,  For sufficiently large $i$, we call the uniquely determined section $\Gamma_i$ of $\nu$ over
$B_R(p) \cap H$ \[\D := \Gamma_i : B_R(p) \cap H \ra  \Psi_i(B_R(p_i) \cap H_i)\] the asymptotic  \textbf{identification map}, briefly the $\D$\textbf{-map}.
\end{definition}

Obviously we can extend the notion of $\D$ maps to any domain on $H$ with compact closures in $H \setminus \Sigma_{H}$.

\begin{remark} \label{te} \quad We gently suppressed that $\D(\p B_R(p)) \neq \p B_R(p_i)$ but only gradually approaches $\p B_R(p_i)$, for $i \ra\infty$. However, it is only the portion away from these boundaries we are interested in. Thus we may easily adjust the definitions near the boundary to our needs and henceforth ignore these negligible adjustments.\\

When we write $id_{H}$ for the zero section of the normal bundle $\nu$, the mentioned local $C^k$-convergence of the $H_i$ can be expressed in terms of a local convergence
$$|\D - id_{H}|_{C^k} \ra 0.$$

Also, since we use almost minimizers, a way to omit the use of the maps $\Psi$ is to locally embed the converging sequences into $\R^{n+1}$ in the first place. The  DeGiorgi-Allard regularity theory, cf.[Si2],Ch.5,  give us, for any given pair of smooth balls $B_1 \subset \overline{B_1} \subset B_2 \subset H$, bounds $|\D - id_{H}|_{C^k}$ in terms of the flat norm distance and the local distortion of $M_i$ relative $M$.\qed
\end{remark}

Keeping the setup and notations we discussed so far, we define a notion to describe structures \emph{naturally assigned} to hypersurfaces $H \in {\cal{H}}$.

\begin{definition} \label{cr} \emph{\textbf{(Natural Structures)}}\, An assignment $F$ of functions $H \mapsto F_H$ defined on
$H \setminus \Sigma_H$, for any $H \in {\cal{H}}$,  is called \textbf{natural}, when $F_H$ commutes with convergence of the underlying spaces in the following sense:\\

For any compactly flat norm converging sequence $H_i\in {\cal{H}}$ and $p_i \in H_i \setminus \Sigma_{H_i}$, $p_i \ra p \in H \setminus \Sigma_{H}$ as above, there is a neighborhood $U(p) \subset H \setminus \Sigma_{H}$, so that
        \[|F_{H_i} \circ \D - F_{H}|_{C^k(U(p))} \ra 0, \mm{ for } i \ra \infty,\]
for some $k=k(F) >0$.\\

 More generally, when the assignment  $H \mapsto F_H$ maps $H$ to some $F_H$ in a category of  objects that admits an $\D$-map pull-back to $U(p)$ and a topology so that the pull-backs converge to the object assigned to $H$ we call the assignment a \textbf{natural structure}.\\
\end{definition}

\begin{example} \textbf{(Naturality on ${\cal{H}}$)}  \label{ecr} \, The $C^k$ convergence of the minimizers readily shows that algebraic expressions in terms of germs of $g_H$ and $A_H$, like the classical curvatures or the norm $|A|$, are natural functions. The \emph{Laplacian}, the
\emph{Jacobi field operator} $J_H$ or the \emph{conformal Laplacian} $L_H$ are natural operators.\\

However, the \emph{metric distance function} to the singular set $dist(\cdot, \Sigma)$ is \emph{not} natural, since the singular sets
of converging hypersurfaces and their limits are hardly correlated, for instance, when smooth hypersurfaces degenerate to singular ones.\\

On the other hand, since we restrict to $H \in {\cal{H}}$, we can find natural structures not defined from merely local data, and this includes skin transforms.\\
\end{example}

\subsubsection{Metric Skin Transforms}\label{localskin}

\bigskip
In this and the next chapter we describe methods to define skin transforms. In both cases we start with the definition of \sks\ and use them to assemble functions.\\

\noindent \textbf{Basic Concepts and Results} \quad  The easiest way to define skin transforms is to use distance tubes of the $|A|$-level sets, where we rescale the distance
according
to the value of $|A|$ on the level sets, to ensure the cone reducibility of the definition.\\

 More precisely, for any non-totally geodesic $H \in {\cal{H}}$, any given $\alpha >0$ and
$c>0$, we define the \sks\ $\M_c=\M_c(\alpha)$ of the desired skin transform $\bp_{\alpha}$ as
$$\M_c(\alpha):= \mm{\emph{ the boundary of the} } \alpha/c \mm{\emph{-distance tube }} \U^\alpha_c \mm{ \emph{of }} |A|^{-1}([c,\infty)).$$
We have $\overline{\U^\alpha_d} \subset \U^\alpha_c$ and therefore $\M_c \cap \M_d \v$, for $c < d$, since  $|A|^{-1}([d,\infty)) \subset |A|^{-1}([c,\infty))$ and $\alpha/d < \alpha/c$.\\

 Thus we can uniquely define
\begin{equation} \bp_{\alpha}(x):= c, \mm{ for } x  \in \M_c. \end{equation}
However, usually $\bigcup_{c > 0} \M_c \neq H \setminus \Sigma$. There could be open sets where $|A|$ remains constant, when $H$ is non-analytic. (These jump values are
 non-regular and thus the set of such values has measure zero and actually in the present case it is countable.)
Also $|A|$ vanishes identically on any totally geodesic path component $H_0 \subset H$.\\

But we can canonically extend this definition. We set $\bp_{\alpha} \equiv 0$ on path components of $H \setminus \Sigma$ which are totally geodesic. Thus we may assume $H \setminus \Sigma$ is
connected and not totally geodesic. (There is  no need to distinguish between the connectivity properties of $H$ and of $H \setminus \Sigma$, since $H$ is connected $\Leftrightarrow$ $H \setminus \Sigma$ is connected, cf. \ref{coh}.)\\

From this, we notice that $H \setminus \Sigma = \bigcup_{c>0} \U^\alpha_c$.

\begin{definition} \emph{\textbf{(Metric Skin Transforms)}} \label{lsk} \quad For $\alpha > 0$, we define
\begin{equation} \bp_{\alpha}(x):= \sup \{c \,| \, x \in \overline{\U^\alpha_c}\}.\end{equation}
for any $x \in H \setminus \Sigma$. For totally geodesic $H$, we set $\bp_{\alpha}\equiv 0$.
\end{definition}

\begin{proposition}\label{locskin} \emph{\textbf{(Relations between $\mathbf{\bp_\alpha, |A|}$ and distances)}}\,
For any non-totally geodesic $H$ we have the following estimates.
\begin{description}
    \item[A.] \emph{\textbf{Growth Estimates and Lipschitz Properties}}
    \begin{enumerate}
\item For any $x \in H \setminus \Sigma$ we have  \begin{equation} \label{adis} \bp_{\alpha}(x) \ge \alpha / dist(x,\Sigma) \end{equation}
\item $\delta_{\bp_{\alpha}} = 1/\bp_{\alpha}$ is compatible with the usual distance function
\begin{equation}\label{dis}  |\delta_{\bp_{\alpha}}(p)- \delta_{\bp_{\alpha}}(q)|   \le d(p,q)/\alpha \end{equation}
In more analytic terms, $\delta_{\bp_{\alpha}}$ is an $1/\alpha$-Lipschitz function.
\end{enumerate}
    \item[B.] \emph{\textbf{Interpolation Properties}}
    \begin{enumerate}
\item $\bp_{\alpha} \ge \bp_{\beta}, \mm{ for }\, \alpha \ge \beta >0$.
\item $\bp_{\alpha} \ra |A| \mm{ in } L_{loc}^\infty, \mm{ for }\,\alpha \ra 0$.
\item $1/\alpha \cdot \bp_{\alpha}  \ra 1/ dist(\cdot,\Sigma)\mm{ in } L_{loc}^\infty, \mm{ for }\,\alpha \ra \infty$.
\end{enumerate}
\end{description}
\end{proposition}
Here,  $d(x,y)$ denotes the distance between $x$ and $y$, $dist(\cdot,\Sigma)$ the distance to $\Sigma$ on $H$ in terms of $g_H$. The latter
convergence results are on $H \setminus \Sigma$.\\

 {\bf Proof} \quad
For \textbf{A(i)} we note that $c = \bp_{\alpha}(x)$ means that  $|A| \le c$ on $B_{\alpha/c}(x)$. Thus we infer from \ref{tt} that for all balls $B_r(x)$, with $r < \alpha/c$: $B_r(x) \cap \Sigma \v$.\\

 In other words:
$dist(x,\Sigma) \ge \alpha/c$ and thus $\bp_{\alpha}(x) \ge \alpha / dist(x,\Sigma)$.\\

For \textbf{A(ii)},  again we may assume that $p \in \M_c$ and $q \in \M_d$ for some $d > c >0$.
\[|\delta_{\bp_{\alpha}}(p)- \delta_{\bp_{\alpha}}(q)| = |1/\bp_{\alpha}(p)- 1/\bp_{\alpha}(q)| =  \alpha^{-1} \cdot \left|\frac{\alpha}{c} - \frac{\alpha}{d}\right| \le
 d(p,q)/\alpha.\]
The latter inequality follows from the inclusion  $\overline{\U^\alpha_d} \subset \U^\alpha_c$.\\

\textbf{B(i)} and \textbf{B(ii)} follow easily from the definitions. Note that $\U^\beta_c \subset \U^\alpha_c$, for $\alpha \ge \beta >0$.\\

For \textbf{B(iii)}, we observe that  the definition $d:= \bp_\alpha(x) = \sup \{c \,| \, x \in \overline{\U^\alpha_c}\}$ means that $dist(x, |A|^{-1}([d,\infty))) \le \alpha/d$. Also
we recall that $|A|^{-1}([d,\infty))$ shrinks to $\Sigma$, while $ d \ra \infty$ and we conclude that for fixed $x \in H \setminus \Sigma$ and any $\ve > 0$ there is some
$d_0(\ve)
>0$ so that for any $d \ge d_0(\ve)$
\begin{equation}\label{vee} dist(z,\Sigma) \le (1 + \ve) \cdot \alpha /\bp_{\alpha}(z), \mm{ for any } z \in B_{dist(x,\Sigma)/4}(x) \end{equation}
Now A(i) tells us that $\alpha \ra \infty$ implies $\bp_{\alpha}(z)  \ra \infty$ in a locally uniform manner.\\

Therefore, for any given $\ve >0$, we can find an $\alpha_0(\ve)>0$ so that for $\alpha > \alpha _0$, (\ref{vee}) holds. Using A(i) again,  gives us the reversed inequality and
we conclude the $L_{loc}^\infty$-convergence of $1/\alpha \cdot \bp_{\alpha}(x)$ to $1/ dist(x,\Sigma)$. \qed

 \begin{proposition}\label{locskintr} \quad
 $\bp_\alpha$ is a skin transform, for any $\alpha > 0$.
\end{proposition}

 {\bf Proof} \quad With the preparation above we can now check the validity of the axioms from Definition  $1$ in Ch.2.

\begin{description}
  \item[(S1) and (S2)]  \,  The properties that $\bp_H \equiv 0$, if $H \subset M$ is totally geodesic and when $H$ is not totally geodesic, then $\bp_H$ is strictly positive
 with $\bp_H \ge |A_H|$. Since distances on $H$ resp. $|A|$ scale by $\lambda$ resp. $1/\lambda$ under scalings of $H$ by $\lambda$, we get $\bp_{\lambda \cdot H} \equiv \lambda^{-1} \cdot  \bp_{H}$ right from the definition.
    \item[(S3)] \, The proof of the Hardy inequality uses some additional  techniques we yet have to develop. We postpone the proof  to Ch.\ref{hardy} below.
    \item[(S4)]  \, The Lipschitz regularity of $\delta_{\bp_\alpha}:=1/\bp_{\alpha}$ with Lipschitz constant $1/\alpha$ is \ref{locskin}A(ii).
    \item[(S5)]  \, The naturality of $\bp_{\alpha}$ follows from that of $|A|$ and the regularity theory.  Namely,  the convergence of a sequence $H_i \in {\cal{H}}$ to a limit $H \in {\cal{H}}$, expressed in terms of $\D$-maps,  in an arbitrarily small ball of $H$, extends to compact smooth convergence on the entire space
        $H \setminus \Sigma_{H}$, in fact, for connected $H$, this complement is rectifiably connected cf. \ref{coh}.\\

From a compact  convergence of minimizers, we get compact convergence of
$|A|$ and thus compact  $L^{\infty}$-convergence of $\bp_{\alpha}$: either the value in a converging sequence of $p_i \in H_i \setminus \Sigma_{H_i}$
remains lower positively bounded by some $d >0$ and then $p_i \in \U^\alpha_d$ and thus $p \in \U^\alpha_d$  so that the convergence of $|A|$ even gives
a convergence of the values $\bp_{\alpha}(p_i)$. Or $\bp_{\alpha}(p_i) \ra 0$, but then the definition shows that the $H_i$ converge to a totally geodesic limit.\\

Next, we have  a uniform Lipschitz estimate for $\delta_{\bp_\alpha}$  on all minimizers,  therefore the Rellich compactness gives us
the H\"older convergence under smooth convergence of the underlying hypersurfaces.\qed
\end{description}

\textbf{Alternative Approaches} \, Arguably, this is the simplest construction of a skin transform.  Alternative definitions of \sks\ use hypersurfaces within $H \setminus \Sigma$ spanned over the sets $|A|^{-1}([c,\infty))$ by other natural processes. We mention two of them.\\

 1. We can choose area minimizing hypersurfaces within $H \setminus \Sigma$  with obstacles $|A|^{-1}([c,\infty))$. Then the \sks\ have \emph{non-negative mean curvature} relative $H$ and suggest a decomposition of $H$ into pieces with area minimizing boundaries, cf. [G] for a discussion of the closely related concept of Plateau-Stein manifolds due to Gromov.\\

Due to results of David and Semmes [DS], we explain in \ref{ds} below, these \sks\ bound domains with regularity properties allowing a particularly simple definition of skin uniform domains, cf.Ch.\ref{ch}.\\

Also the local distances between any two adjacent \sks\ $\M_a,\M_b$ are mutually controlled through Harnack inequalities since $\M_a$ can locally be regarded as positive solutions of elliptic equations on $\M_b$ and vice versa. This is the source for the Lipschitz-continuity of the associated $\bp$-distance, based on arguments of Solomon [So]. \\

2. Alternatively, we observe that being an area minimizer is a scaling invariant property. The compactness result for these minimizers cause some coarse self-similarities of area minimizers and their singular sets when we zoom into the singular set. We use this observation, for instance, in the assembly of surgeries along $\Sigma$ in [L2].\\

 Now, the treatment of multifractals, cf.[Fa],Ch.17, suggests methods to interpolate between fractals of different degree of complexity. In our case this would be $\Sigma$ versus $H \setminus \Sigma$ and can be adapted to define skin transforms.\\

 We discuss particularities and applications of such alternative approaches to skin transform elsewhere.\qed

\textbf{Derived Skin Structures} \quad   By definition, we consider skin structures on area minimizing hypersurfaces in the class ${\cal{H}} =\bigcup_{n \ge 2} {\cal{H}}_n$
described in the introduction. This class is closed under blow-ups around singular points. To incorporate also dimensional reduction arguments it is important to also include the following natural class ${\cal{K}}$  of non-minimizing hypersurfaces:\\

For any area minimizing cone $C \subset \R^{n+1}$ with tip $0$, we define $S_C:= \p B_1(0) \cap C \subset S^n \subset  \R^{n+1}$ in the
unit sphere $S^n$. We set ${\cal{K}}_{n-1}:= \{ S_C\,| \, C \in {\mathcal{C}_{n}}\}$ and write ${\cal{K}}= \bigcup_{n \ge 2} {\cal{K}}_{n-1}$, for the space of all such hypersurfaces.\\

As mentioned earlier, the $S_C \in {\cal{K}}$  are almost minimizers, cf.Appendix II, \ref{riem}. They share crucial properties including the regularity theory
with proper area minimizers, since the cone $C$ over $S_C$ is area minimizing. This also applies to skin structures.

\begin{lemma}\label{skino}  Let $\bp$ be a skin transform on ${\cal{H}}$, then the skin axioms still hold on ${\cal{K}}$, for the restriction $\bp_{S_C}:= \bp_C|_{S_C}$, for any
$C \in {\mathcal{C}_{n}}$, $n \ge 2$.
\end{lemma}

 The proof  is largely straightforward. We  notice $\bp_{S_C}= \bp_C\big|_{S_C} \ge |A_C|\big|_{S_C}=|A_{S_C}|$. Also the skin transform induced on blow-ups around any singular point $p \in S_C$
 equals the restriction of the skin transform on the product tangent cone $\R \times C^*$ in $(1,p) \in C$ to $\{0\} \times C^*$.
  The remaining
 details are left to the reader.\qed

 Alternatively to the given definition $\bp_{S_C}:= \bp_C|_{S_C}$, we could go back to the construction of the given $\bp$ and rebuild it
also for the $S_C$, $C \in {\mathcal{C}_{n}}$. Let us call the result $\bp^*_{S_C}$. Then we experience that  $\bp_{S_C} \neq \bp^*_{S_C}$.
However, the choice we made is the preferable one, since it is directly linked to that on $C$. We actually use this in the other direction, to study $\bp$ on $C$ from $\bp_{S_C}$. \\

\textbf{Metric Skin Transforms and Regularity Theory} \quad There is an another interpretation of $\bp_1$ in the regularity theory of area minimizers which widens the scope of this particular approach in a rather different direction. In recent and complementary work of Cheeger and Naber [CN] we find the notion of the regularity scale $r_H$ defined for any $p \in H \setminus \Sigma$ used to refine classical regularity estimates  for area minimizers and estimates for the dimension of the singular set.\\

It is defined as $r_H(p):= \sup \{r > 0 \, | \, r^2 \cdot |A|^2 \le 1 \mm{ on } B_r(p)\}$. Using this type of invariants can be traced back to the 50s when Heinz [Hz] revisited the
Bernstein theorem. He showed that for minimal graphs over  $B \subset \R^2$ the \emph{scaling invariant} quantity $r^2 \cdot |A|^2(p)$ on disks $B$  with center $p$ and
radius $r$, can be a priori estimated from above. (Thus when $r \ra \infty$, we get $|A| \ra 0$ and the graph must be a plane, proving the classical Bernstein theorem.)\\

Now, there is a pleasingly simple relation to $\bp_1$:  $\delta_{\bp_1} = r_H$. Namely, we notice that $\rho:= r_H(p) = 1/\max_{\overline{B_\rho}(x)} |A| =: 1/m$. Now take
some $y \in {\overline{B_{1/m}}(x)}$ with $|A|(y) =m$. Then $x$ has distance $\le 1/m$ from $y$. Thus $m \le m^* := \sup \{c \,| \, x \in \overline{\U^1_c}\}$. In fact, $m =
m^*$, otherwise we had a point $z \in B_{1/m}(x)$ with $|A|(z) =m^* >m$.\\

In turn, the regularity theory in [CN] involves versions of $r_H$ containing higher derivatives of $|A|$, made of the terms $r^{k+1} \cdot |\nabla^k A|$, which could also be chosen as a basis to define further examples of skin transforms.\qed

 \textbf{The Space of Skin Transforms} \label{all} \quad   The previous discussion shows that besides the different methods to define  skin transforms, any given skin transform induces others we also need to take into account.  We mention two further operations on the space of skin transforms.
\begin{itemize}
    \item \textbf{Convexity} \, With any two skin transforms $\bp$ and $\bp^*$ any convex combination $c \cdot \bp + (1-c) \cdot \bp^*$, $c \in (0,1)$ is again a skin
        transform. Thus the space of all skin transforms is convex and thus contractible.
    \item  \textbf{Natural Averages} \, Mollifiers applied to $|A|$, defined by convolution integrals with smoothing kernels, so that the smoothing kernel is sized by
        $\delta_{\bp}$, in place of the typically chosen $1/dist(\cdot,\Sigma)$, give new skin transforms.
One version of this idea is the Whitney type smoothing we discuss in \ref{skinada1} below.
 \end{itemize}
Summarizing, there is no use to stick with one particular model of a skin transform. Instead, we
always consider the class of all skin transforms, this automatically brings us to focus on the axioms for $\bp$.\\

\bigskip

 \setcounter{section}{3}
\renewcommand{\thesubsection}{\thesection}
\subsection{Localizations and Controlled Covers}
\bigskip

We first describe a technique to find locally finite covers of area minimizers by  balls with well-controlled geometry and intersection numbers, in  \ref{skinada}. These \emph{skin
adapted covers} are used to localize analytic problems on area minimizers.\\

In \ref{hardy} we use such a cover to prove the remaining axiom (S3) for the $\bp_\alpha$ of \ref{localskin}, the Hardy inequality. Thereby we complete the certification
process for one working model of a skin transform. We keep it as a gold reserve, since, from that point on,  we no longer refer to a particular model of a skin transform but use an arbitrary one.\\

 \textbf{We derive all further results directly for the axioms (S1)-(S5).} This said,  the axiom (S3) will not be exploited in this paper, but only in [L1] and [L2], where it is used to couple the developed potential theoretical results to classical operators.\\

Finally, in \ref{skinada1}, we use the covering techniques of \ref{skinada} to derive smoothing techniques in the style of Whitney [Wh1] and Stein [St]. We apply them to the a priori merely locally Lipschitz skin transforms.

\subsubsection{Skin Adapted Covers}\label{skinada}
\bigskip

For the purposes of this section we use an assignment $\bp$ that satisfies all axioms of a skin transform except for (S3), we call it a \emph{pre-skin transform}.
That is, we do not presume that the Hardy inequality holds since we use the argument of this section to establish the Hardy inequality for the $\bp_{\alpha}$.
Actually all results in this paper hold for pre-skin transforms. The validity of axiom (S3) is only used in the analytic and geometric applications we discuss in [L1] and [L2], \\

Here we prove the existence of covers of $H \setminus \Sigma$ particularly adapted to $\bp$. To this end we recall two classical covering methods for sets in $\R^n$:
\begin{itemize}
    \item  Any subset $P \subset \R^n$ admits covers by balls of various radii but with bounded covering number, the so-called Besicovitch covers, cf.[Di], Ch.18.
    \item For closed $P \subset \R^n$, $\R^n \setminus P$ admits Whitney covers by closed cubes with disjoint interior with the particular property that the diameters of
        the cubes are comparable to their distance to $P$, cf.[Wh1], and [St],Ch.6.\\
\end{itemize}

On $H \setminus \Sigma$ the geometry degenerates while we approach $\Sigma$ and to mimic the classical covering methods we also need to use gradually smaller balls the
closer we approach $\Sigma$. However, to keep the cover useful for analytic estimates we need that  - after scaling to unit size - all these balls belong to some compact class of
smooth geometries with a positive \emph{lower} bounded Hausdorff distance to the flat metric. Thus we need to balance the decay of the radii towards $\Sigma$.\\

We accomplish these somewhat opposing goals through a blend of the Besicovitch and the Whitney covering methods with the additional twist that the radii are controlled from a
pre-skin transform. For the sake of a consistent statement, we include the case where $H$ is totally geodesic. Then this result merely expresses the surjectivity of the exponential
map with an infinitely sheeted covering in the case where $H$ is compact. In the proof we focus on the non-totally geodesic case.

\begin{proposition}\label{besi} \emph{\textbf{(Skin Adapted Covers)}}  \, For any $H \in {\cal{H}}$ and any size parameter $\xi \in (0,\xi_0)$,
for some  $\xi_0(n,L_{\bp}) \in (0,1/(10^3 \cdot L_{\bp}))$, we get:\\

A locally finite cover ${\cal{A}}$ of $H \setminus \Sigma$ by closed balls
\[{\cal{A}} = \{\overline{B_{\Theta(p)}}\,|\, p \in A \},\mm{with radius } \Theta(p):=   \xi / \bp(p) = \xi \cdot \delta_{\bp}(p),\] for some discrete set
$A \subset H \setminus \Sigma,$ so that for a suitably small neighborhood $Q$ of $\Sigma$:
\begin{description}
    \item[(C1)] For  $p \in Q$ the exponential map $\exp_p|_{B_{100 \cdot \Theta(p)}(0)}$ is  bi-Lipschitz onto its image, for some bi-Lipschitz constant $l(n) \ge 1$
    \item[(C2)] $A^Q := A \cap Q$ splits into $c(n)$ disjoint families $A^Q(1),...,A^Q(c)$ with
\begin{enumerate}
\item $B_{10 \cdot \Theta(p)}(p) \cap B_{10 \cdot\Theta(q)}(q) \v$, for $p$ and $q$ in the same $A^Q(k)$
\item $q \notin \overline{B_{\Theta(p)}(p)}$, for any two $p, q \in A^Q$.
\end{enumerate}
\end{description}
In particular, the covering number for $z \in Q$, for any $\rho \in (0,10)$, by balls centered in $A \cap Q$, is uniformly bounded:
 \begin{equation} \label{conn}\cs(A \cap Q,z,\rho):= \cs\{x \in A \cap Q \,|\, z \in B_{\rho \cdot \Theta(x)}(x)\} \le c(n). \end{equation}
We call such a cover ${\cal{A}}$, a \textbf{skin adapted cover}. The covering method also ensures the following details.
\begin{enumerate}
    \item Let $H \subset \R^n$ be a complete, non-totally geodesic area minimizing hypersurface. Then we may choose $Q = H$.
    \item For any $\ve  >0$,  we can find some $\xi_{\ve} \in (0,\xi_0(n,L_{\bp}))$ so that, for  every $p \in H \setminus \Sigma$ the exponential map
        $\exp_p|_{B_{100 \cdot \xi_{\ve}/ \bp(p)}(0)}$ is bi-Lipschitz onto its image with bi-Lipschitz constant $1 + \ve$.\\
\end{enumerate}
\end{proposition}

In the proof of \ref{besi}, and of many other results below, we use the following variations of the Lipschitz continuity property of the $\bp$-distance $\delta_{\bp}$

\begin{lemma}\label{locsking}  \quad Let $\bp$ be a pre-skin transform, so that the $\bp$-distance is Lipschitz continuous
with Lipschitz constant $L >0$. Then we have
\begin{enumerate}
\item For any $x \in H \setminus \Sigma$   \begin{equation} \label{adiss} L \cdot \bp(x) \ge 1/ dist(x,\Sigma) \mm{ or equivalently }  \delta_{\bp}(x) \le L \cdot
    dist(x,\Sigma)\end{equation}
\item $\bp$ is locally Lipschitz \begin{equation}\label{liph} \left|\bp(x)/\bp(p)  - 1\right|  \le 2 \cdot L \cdot \bp(p) \cdot d(x,p)\end{equation}  for any $q \in B_{1/(2 \cdot L \cdot \bp(p))}(p)$\\
\end{enumerate}
\end{lemma}

 {\bf Proof} \quad  For (i), we use the particularity of area minimizing \emph{hypersurfaces} that for any $p \in \Sigma$ there is a sequence $p_n \ra p$, with
 $|A|(p_n) \ra \infty$, for $n \ra \infty$, cf.\ref{tt}. $\bp \ge |A|$ implies that also  $\bp(p_n) \ra \infty$.\\

 Thus, for $k >0$, we have some $n_k$, so that for $n \ge n_k$:  $\bp(p_n) \ge k$. Then we have for any $x \in H  \setminus \Sigma$ with $d(x,p_n) \le 2 \cdot (k \cdot L)^{-1}$
\[\delta_{\bp}(x) \le |\delta_{\bp}(x) - \delta_{\bp}(p_n)| + 1/k \le  L \cdot d(x,p_n) + 1/k \le 3/k \]
We may assume that $d(p,p_n) <  (k \cdot L)^{-1}$ and thus $\delta_{\bp}(x)\le 3/k$ on $B_{(k \cdot L)^{-1}}(p)$. In other words, we get the coarse estimate that
$\bp \ge k$ on $U_{(3 \cdot k \cdot L)^{-1}}(\Sigma_H)$.\\

We use this when we start from $|\delta_{\bp}(x) - \delta_{\bp}(p)|   \le L \cdot d(x,p)$, for $x,p \in H  \setminus \Sigma$. Now we consider a
 sequence $p_n \in H  \setminus \Sigma$ with  $p_n \ra p_\infty \in \Sigma$, $n \ra \infty$,
 so that $d(x,p_n) \le dist(x,\Sigma)$. Since $\delta_{\bp}(p_n) \ra 0$, the claim follows.\\

  (At this point we not yet know that $dist(x,\Sigma) < \infty$ and thus that the assertion is non-trivial. This is only proved in \ref{fini}.)\\

  For (ii), we choose any two points $p, q \in H \setminus \Sigma$. Then the Lipschitz inequality $|\delta_{\bp}(p) - \delta_{\bp}(q)|   \le L \cdot d(p,q)$
  can be rewritten as
\[|\bp(p)-\bp(q)|  \le L \cdot \bp(p) \cdot \bp(q) \cdot d(p,q)\]
From this we see
\begin{equation}\label{2lip}
 \bp(x) \le 2 \cdot \bp(p),\mm{ for all } x \in  B_{1/(2 \cdot L \cdot \bp(p))}(p)
\end{equation}
Thus we have
\[|\bp(p)-\bp(q)| \le 2 \cdot L \cdot  \bp^2(p)  \cdot d(p,q)\]  for any $q \in B_{1/(2 \cdot L \cdot \bp(p))}(p)$.\qed\\

 {\bf Proof of \ref{besi}} \quad There are three main steps. In  the first two steps we introduce two scalings to derive pointwise estimates for
 the volume of balls within $B_{\Theta(p)}(p)$.
We use the Harnack type inequality (\ref{liph}) for $\bp$ to turn these estimates into locally uniform estimates.
 Finally we show how this fits into the typical combinatorics of Besicovitch type covers
to derive the claimed properties of ${\cal{A}}$.\\

In each substep of the argument we choose the most appropriate local scaling of $H$ to derive the needed estimates for the
subsequent steps. These are only \emph{temporary scalings} and we undo them once we derived the result. \\

\noindent \textbf{Step 1 (Scaling of $\mathbf{\Sigma \subset H \subset M}$)} \\

As before $L = L(\bp) >0$ denotes the Lipschitz constant for the $\bp$-distance.
 \[|\delta_{\bp}(p)- \delta_{\bp}(q)|   \le L \cdot d(p,q), \mm{ for } p,q \in  H \setminus \Sigma.\]
 Scaling of any ball $B_{1 /(L \cdot \bp(p))}(p)$, for $p \in H \setminus \Sigma$ by
$L \cdot \bp(p)$
produces a ball of radius $1$ in the scaled space $L \cdot \bp(p) \cdot H$.\\
\[B_1(p) \subset L \cdot \bp(p) \cdot H \subset L \cdot \bp(p) \cdot M  \]

Since $L \cdot \bp(x) \ge 1/dist(x, \Sigma)$, we observe that the closer $p$ approaches $\Sigma$, the stronger the scaling of the smooth compact ambient manifold
$M$ of $H$ becomes, that is, $L \cdot \bp(p) \cdot M$ will look nearly flat, for $p$ very close to $\Sigma$.\\

Formally, we denote the exponential map of $s \cdot M$ in $p$, scaled by $s \ge 1$, by $\exp_p[s \cdot M]: T_p M \ra s \cdot M$. $T_p M$ carries the
 flat metric $g_{T_p M}$.\\

Then, for any $\ve >0$, we find a neighborhood $W(\ve) \subset H$ of $\Sigma$ so that for any $p \in W(\ve)$:
 \begin{equation} \label{scc} |\exp^*_p[\bp(p) \cdot M](L^2 \cdot \bp(p)^2 \cdot g_M) - g_{T_p M}|_{C^5(B_{100}(0))} \le \ve \end{equation}
where the $C^5$-norm and the radius are measured relative $\bp(p) \cdot M$.\\\\

\noindent  \textbf{Step 2 (Locally Uniform Estimates on $\mathbf{H \setminus \Sigma}$)} \label{co}\\

1. The first step says that, when our region of interest $W(\ve)$ shrinks to $\Sigma$, that is, we choose some small $\ve$,  and we scale $M$ by $\bp(p)$, for $p \in W(\ve)$,
then $M$ becomes virtually flat since we have a lower bound for $\bp|_{W(\ve)}$ that diverges, when $\ve \ra 0$.\\

$H$ behaves differently, since
approaching $\Sigma$ also means that $|A|$ diverges. This time $\bp \ge |A|$ shows that
\[|A|(p) \le 1 \mm{ after scaling by }  \bp(p).\]

This pointwise estimate implies also local estimates when we use the Harnack type inequality  (\ref{liph}) below, for $\bp$:
 \begin{equation} \label{aest} \left|\bp(x)/\bp(p)  - 1\right| \le  2 \cdot L \cdot \bp(p) \cdot d(x,p) \le 1,\end{equation}
 for any $p \in H \setminus \Sigma$, $x \in B_{1/(2 \cdot L \cdot  \bp(p))}(p) \subset H$.\\

$\bp \ge |A|$ and (\ref{aest}) imply that there is a constant $A_n \ge 1$ depending only on $n$ so that \[|A| \le A_n \mm{ on } B_1(p) \subset L \cdot \bp(p) \cdot H.\]\

2. This locally uniform bound on $|A|$ implies some locally uniform volume estimates for balls in $H$: Gauss equations relating the curvature tensors of $H$ and $M$ and
Rauch's comparison theorem show that for some $\ve
>0$ small enough, we can find positive functions $\Lambda(\zeta),\eta(\zeta)$ with
\[\Lambda(\zeta) \ra \infty, \,  \eta(\zeta) \ra 0,\mm{ for } \zeta \ra 0,\] so that for any $p \in W(\ve)$, measured relative to $\Lambda(\zeta) \cdot L \cdot \bp(p) \cdot H$,
$exp_p[\Lambda(\zeta) \cdot  L \cdot \bp(p)
\cdot H ]$ is a local diffeomorphism from $B_{10^3}(0)$ onto its image in $H$, so that
\begin{equation}\label{hbg0} |\exp^*_p[\Lambda(\zeta) \cdot  L \cdot \bp(p)
\cdot H ](\Lambda^2(\zeta) \cdot  L^2 \cdot \bp(p)^2 \cdot g_H) - g_{T_p H}|_{L^\infty(B_{100}(0))} \le \eta(\zeta).\end{equation} The regularity theory of $H$ allows us to
upgrade to $C^k$-estimates for any given $k \ge 0$ and we get, keeping the same constants for simplicity,
\begin{equation}\label{hbg}  |\exp^*_p[\Lambda(\zeta) \cdot  L \cdot \bp(p) \cdot H ](\Lambda^2(\zeta)    \cdot L^2 \cdot \bp(p)^2 \cdot g_H) - g_{T_p H}|_{C^5(B_{100}(0))} \le
\eta(\zeta).\end{equation}\

3. Thus, once $\ve >0$ has been chosen small enough, we have a uniform control, for any $p \in W(\ve)$:\\

3.A.\,  $\exp_p[\Lambda(\zeta) \cdot  L \cdot \bp(p) \cdot H]$ is a bi-Lipschitz map from $B_{100}(p)$ to its image for a bi-Lipschitz constant $l(\zeta) \ge 1$, with
    $l(\zeta) \ra 1$ for $\zeta \ra 0$.\\

3.B. When we choose $\zeta >0$, so that $l(\zeta) \in [1,2]$,  the following set of volume estimates hold for $z \in B_{50}(p)$ with volumes and radii measured relative
$\Lambda(\zeta) \cdot   L \cdot \bp(p) \cdot H$
\begin{equation}\ \label{vol} k_1 \le Vol(B_{1/3}(z))\mm{ and } Vol(B_{30}(z))  \le k_2,\end{equation}
for some constants $k_i(n,\bp) >0,  i=1,2$.\\

4. We choose some $\zeta \ll 1$, so that $\Lambda(\zeta) \gg 1$, in particular, with $L \cdot \Lambda(\zeta) >100$,   we set $\xi(\zeta):=1/(L \cdot \Lambda(\zeta))$ and
 \begin{equation}\label{oic}\Theta(p):= \xi(\zeta) / \bp(p) = 1/(L \cdot \Lambda(\zeta) \cdot \bp(p)).\end{equation}\

\noindent \textbf{Step 3 (Combinatorics)} \\

We choose a countable dense set $S \subset H \setminus \Sigma$, $S= \{a_m\,|\,m \in \Z^{\ge 0}\}$ and start with the cover
 \[{\cal{B}} = \{\overline{B_{\Theta(p)}(p)}\,|\, p \in S\}\]

Now we define a map $i: S \ra\Z^{\ge 0}$ by induction: we set $i(a_0):=1$ and {\small\[ i(a_{k+1}):=0, \mm{ if } a_{k+1} \in \bigcup_{i \le k} \overline{B_{\Theta(a_i)}(a_i)},
\mm{ and, otherwise,}\]
 \[ i(a_{k+1}):= \min(\{m \le k \,|\,d(a_m,a_{k+1}) >  10 \cdot \Theta(a_m) + 10 \cdot \Theta(a_{k+1})\} \cup \{k+1\}).\]} We define
 \[{\cal{A}}(j) := \{\overline{B_{\Theta(p)}(p)}\,|\, p \in S, i(p)=j\},\, {\cal{A}} := \bigcup_{j\ge 1} {\cal{A}}(j) \mm{ and } A:= \{a \in S\,|\, i(a) \ge 1\} \]
 We observe that ${\cal{A}}(i) \cap {\cal{A}}(j) \v$ for $i \neq j$. Also these families satisfy the asserted properties (i) and (ii).\\

Now we prove that there is a neighborhood $ Q$ of $\Sigma$  and  a constant $c(n,\bp)$ so that $A(i) \v$ for $i > c$. This also implies the local finiteness of ${\cal{A}}$ since
there are only finitely many balls in ${\cal{A}}$ with center point in $H \setminus Q$.\\

From (\ref{aest}) and (\ref{oic}) we may assume that
\begin{equation}\label{estt} 9/10 \cdot \Theta(p) \le  \Theta(x) \le 11/10 \cdot \Theta(p),\,\mm{ for any }x \in \overline{B_{100 \cdot \Theta(p)}(p)}.\end{equation}
We claim that for \begin{equation}\label{defc}c(n,\bp):= (\mm{\emph{the smallest integer}}  \ge k_2/k_1), \mm{ with $k_i$ as in (\ref{vol}): }\end{equation} ${\cal{A}}(i) \v,$
for $i > c.$ Otherwise, we had some $\overline{B_{\Theta(p)}(p)} \in {\cal{A}}(c+1)$ and thus
\[B_{10 \cdot \Theta(p)}(p) \cap B_{10 \cdot \Theta(x_i)}(x_i) \n\]
for at least $c$ different $x_i \in A$. But then we get a contradiction from the following estimates which use (\ref{estt}) to see that the $(c+1)$ balls $B_{1/3}(x_i)$ and
$B_{1/3}(p)$ are pairwise disjoint. We use (\ref{vol}) relative  $\Lambda(\zeta) \cdot   L \cdot \bp(p) \cdot H$ to see
\[(c+1) \cdot k_1 \le \sum Vol(B_{1/3}(x_i))+ Vol(B_{1/3}(p))=\]
\[Vol \big(\bigcup B_{1/3}(x_i) \cup B_{1/3}(p)\big) \le Vol(B_{30}(p)) \le k_2.\] Finally, we observe that
${\cal{A}}$ is a cover: assume there is some point $q \in U^c=U \setminus \bigcup_{p \in A} \overline{B_{\Theta(p)}(p)}$. Since ${\cal{A}}$ is locally finite, this complement is
open and thus there is a point $z \in U^c \cap S$ with $\overline{B_{\Theta(z)}(z)} \in {\cal{A}}$ and $q \in \overline{B_{\Theta(z)}(z)}$,  a contradiction.  \qed

\bigskip

\subsubsection{Hardy  Inequalities on $H \setminus \Sigma$}\label{hardy}

\bigskip

Here we resume our discussion of the $\bp_\alpha$ of  \ref{lsk} and notice that we already verified all axioms for a skin transform except for (S3)
requiring the validity of the Hardy type inequality for $\bp$, which actually  is a Hardy inequality for the operator
$-\Delta + |A|^2$ relative to the $\bp$-distance $\delta_{\bp}$.\\

In this section we prove this Hardy inequality and apply it to deduce the Hardy inequality also relative to the distance $dist(\cdot,\Sigma)$.\\

To this end we observe that the construction of skin adapted covers of \ref{besi} did not use axiom (S3). Thus we may use \ref{besi} in the
proof of (S3) for $\bp_\alpha$. Thereby, we complete the proof of Theorem 1 asserting that $\bp_{\alpha}$ is a skin transform.\\

\begin{proposition}\label{pi}  \quad For any $\alpha > 0$ and any compact area minimizing  hypersurface
$H^n \subset M^{n+1}$, there is a constant $k_{\alpha ,H} > 0$, so that
\begin{equation} \label{pi1}
\int_H|\nabla f|^2 + |A|^2 \cdot f^2 dA \ge k_{\alpha, H}  \cdot \int_H \bp_{\alpha}^2 \cdot f^2 dA
\end{equation}
for any smooth $f$ compactly supported in $H \setminus \Sigma$.\\
\end{proposition}

The argument for \ref{pi} uses the fact that $H$  is compact. However, in the case where the ambient space is the Euclidean space the result also holds, with an actually
uniform estimate for the Hardy constant:

\begin{proposition}\label{pia}  \quad For any $\alpha > 0$, there is a constant $k_{\alpha ,n} > 0$, so that for any oriented minimal boundary
$H^n \subset \R^{n+1}$,
\begin{equation} \label{pia1}
\int_H|\nabla f|^2 + |A|^2 \cdot f^2 dA \ge k_{\alpha, n}  \cdot \int_H \bp_{\alpha}^2 \cdot f^2 dA
\end{equation}
for any smooth $f$ compactly supported in $H \setminus \Sigma$.\\
\end{proposition}

Both \ref{pi} and \ref{pia} also trivially hold for totally geodesic $H$. We henceforth assume that $H$ is \emph{non-totally geodesic}.

\begin{remark} \label{lao} \quad In the proof that follows we will see that   \ref{pi} is largely a consequence of \ref{pia} from some blow-up argument.\\

This is the reason why we consider complete area minimizing  hypersurfaces in
$\R^{n+1}$ and not only minimal cones. Blow-ups around a fixed singular point lead to cones, but here we use blow-ups with possibly varying basepoints. This
may lead to more general oriented boundaries than cones.\qed
\end{remark}

 We already know that $\bp_\alpha(x) \ge \alpha/ dist(x, \Sigma)$, for $x \in H \setminus \Sigma$, cf. \ref{locskin}(ii). Hence, for $\alpha = 1$, we also get a Hardy
inequality for $1/dist(x, \Sigma)$.

\begin{corollary}\label{piw}  \quad For any compact area minimizing  hypersurface
$H^n \subset M^{n+1}$, there is a constant $k^*_H >0$, namely $k^*_H:= k_{1 ,H}$, so that
\begin{equation} \label{pi1w}
\int_H|\nabla f|^2 + |A|^2 \cdot f^2 dA \ge k^*_H \cdot \int_H 1/dist(x, \Sigma)^2 \cdot f^2 dA
\end{equation}
for any smooth $f$ compactly supported in $H \setminus \Sigma$.\\

Similarly, we get a $k^*_n >0$, so that for any oriented minimal boundary  $H^n \subset \R^{n+1}$,
\begin{equation} \label{pi1t}
\int_H|\nabla f|^2 + |A|^2 \cdot f^2 dA \ge k^*_n \cdot \int_H 1/dist(x, \Sigma)^2 \cdot f^2 dA
\end{equation}
and for any smooth $f$ compactly supported in $H \setminus \Sigma$.\\
\end{corollary}

Before we turn to the proof we review these inequalities in a wider context.

\begin{remark} \label{moti} \quad 1. The latter inequalities (\ref{pi1w}) and  (\ref{pi1t}) are weaker than the skin versions (\ref{pi1}) and (\ref{pia1}), since the inequality
$\bp_\alpha(x) \ge \alpha/ dist(x, \Sigma)$ is far from being sharp. For instance, we could have sequences $p_i \in H \setminus \Sigma$, $p_i \ra p_\infty \in \Sigma$ so
that
\[\bp_\alpha(p_i) \cdot dist(p_i, \Sigma) \ra \infty \mm{ when } i \ra \infty,\]
when the $p_i$ belong to increasingly sharp but still smooth wrinkles of $H$. However, although \ref{piw} is weaker than \ref{pi} and \ref{pia}, it is not simpler, but a proper
consequence of \ref{pi}. cf. \ref{motii} below.\\

2. The Hardy inequality for a Lipschitz regularly bounded domain $D \subset \R^n$ says:
\[\int_D |\nabla \phi|^2 dV  \ge c_D \cdot \int_D (1/dist(x,\p D))^2\cdot \phi^2 dV,\]
for some  constant $c_D > 0$ and any smooth function $\phi$  supported in $D$. This is a stronger version of the usual Poincar\'{e}-Sobolev inequality, cf.[An3], [BM] and the recent extensive expositions in [BEL] and [GN].\\

3. The codimension of the singular set $\Sigma$ is larger than $2$, and due to the compactness of $H$, the appropriate coarea formula e.g. [GMS], Vol1, Ch.2.1.5, Th.3, p.103, shows that \[\inf_{f \in C_c^\infty(H \setminus \Sigma)} \int_H
|\nabla f|^2 dA/\int_H \bp_\alpha^2\cdot f^2 dA = 0,\]  where $C_c^\infty(H \setminus \Sigma)$ denotes the smooth functions supported in $H \setminus \Sigma$.\\

In turn, $|A|$ usually vanishes or converges to zero along suitable point series to $\Sigma$. This corresponds to rays in tangent cones so that the cone is totally geodesic
along this ray. But $\bp_\alpha$ remains positively lower bounded and actually converges to $+\infty$. Then, without a term that prevents the support of a series of test
functions to concentrate around these points,  we get  \[ \inf_{f \in C_c^\infty(H \setminus \Sigma)} \int_H  |A|^2\cdot f^2 dA /\int_H \bp_\alpha^2\cdot f^2 dA   = 0.\]

That is, the use of the balanced couple of integrands on the left hand side of (\ref{pi1}) is indispensable to ensure a positive constant $k_H$.\qed \\
\end{remark}

Now we collect some tools and observations for the proof of \ref{pi}. We first note that the optimal constant $k_{\alpha, H}$ in  \ref{pi} is nothing but the first eigenvalue for
the weighted operator
  $$P_\alpha u := \bp_\alpha^{-2} \cdot  \left(-\Delta u + |A|^2 \cdot u \right)$$
  on smooth functions defined on $H \setminus \Sigma$. Since $\bp$ is locally Lipschitz, standard elliptic theory implies that
  eigenfunctions of $P_\alpha$ are  $C^{2,\gamma}$-regular,  for any $\gamma \in (0,1)$, cf. [GT], 6.4.\\

More explicitly, the first eigenvalue $\lambda_{P_\alpha} = k_{\alpha, H}$ of the weighted operator $P_\alpha$ can be obtained from its
 variational characterization (Rayleigh quotient) taken over the set $C^{\infty}_c(H \setminus \Sigma)$ of all smooth $f$ compactly supported  in $H \setminus \Sigma$:\\
\begin{equation} \label{ray}\lambda_{P_\alpha} =
\inf_{f \in C^{\infty}_c(H \setminus \Sigma), f \not\equiv 0} \int_H |\nabla f|^2 dA + |A|^2 \cdot f^2 dA \Big/\int_H \bp_\alpha^2\cdot f^2 dA \end{equation}

To estimate this eigenvalue  we localize the problem to Neumann eigenvalues on regular balls. Then we take covers with controlled covering
numbers (Besicovitch covers) by such balls and use them to derive a positive lower  estimate for the eigenvalue $\lambda_{P_\alpha}$.\\

\textbf{Neumann eigenvalues of $P_\alpha$ on  balls} \quad We observe from the shape of (\ref{ray}) that the (Neumann) eigenvalues of $P_\alpha$ on balls are scaling invariant:\\

For any ball $B_r(p) \subset H \setminus \Sigma$, $r > 0$ and any scaling factor $\mu >0$ we have: \begin{equation} \label{eigsc}\nu_\alpha(B_r(p)) = \nu_\alpha(\mu
\cdot B_r(p))\end{equation} Namely, when we scale $B_r(p)$, by any $\mu > 0$, the integrands of the Rayleigh quotient  $|\nabla f|^2 + |A|^2$ (numerator) and
 $\bp^2$ (denominator) change by the same factor $\mu^{-2}> 0$.\\

Now we want to get lower estimates for $\nu_\alpha(B_r(p))$.  We note that  there is no uniform positive lower bound for $r \ra 0$, when  $|A|(p) = 0$.
Conversely, when the balls are too large, we can hardly control its geometry and we can neither understand the eigenvalues nor the covering numbers.\\

However, when we approach $\Sigma$ we get better and better local approximations by analytic hypersurfaces in $\R^n$ (after rescaling to a unit size). This leads us to the
idea to let $\bp$ determine the radius of the balls when we are close to $\Sigma$.\\

To this end we notice that $\bp_{\alpha}(x) \ge \alpha / dist(x,\Sigma)$ means that for any $p \in  H \setminus \Sigma$: \, $B_{\alpha/\bp_\alpha(p)}(p) \cap \Sigma \v$. For
these balls we have the following estimates

\begin{lemma}\label{rad}  \quad For any $\alpha \in (0,1]$ and any $\mu \in (0,1/2)$, there is a neighborhood $U_{\alpha, \mu}$ of $\Sigma$ and a constant $\zeta(U_{\alpha, \mu}) >0$, so that
 \begin{equation}\label{rade}
 \nu_\alpha(B_{\mu \cdot \alpha /\bp_\alpha(p)}(p))  \ge \zeta
 \end{equation}
for any $p \in U_{\alpha, \mu} \setminus \Sigma$.\\
\end{lemma}

{\bf Proof} \quad To simplify the notations we only consider the case $\alpha =1$. Now we explicitly use the definition of $\bp_1$.\\

We assume there is a sequence of points $p_i \in H \setminus
\Sigma$ with $p_i \ra p_\infty \in \Sigma$ so that $\nu_1(B_{\mu/\bp(p_i)}(p_i)) \ra 0$. When we scale $B_{\mu/\bp_1(p_i)}$ by $\bp_1(p_i)$ to the size $B_{\mu}(p_i) \subset
\bp_1(p_i) \cdot H$ we first notice that (after scaling by $\bp_1(p_i)$)\[ \sup \{|A|(x) \, | \, x \in B_1(p_i) \cap \bp_1(p_i) \cdot H\} = 1.\] And we get a subsequence of the sequence
of pointed spaces $(\bp_1(p_i) \cdot H, p_i)$ which converge compactly in flat norm.
\[(\bp_1(p_i) \cdot H, p_i) \ra (H_\infty,0) \subset (\R^{n+1},0)\]
Since the limit is smooth on $B_1(0) \cap H_\infty$ (which we can see as in  \ref{haa}) the convergence towards this subset can be upgraded to $C^k$-topology for any given
$k \ge
0$.\\

 We notice that
\[ \sup \{|A|(x) \, | \, x \in B_1(p_i) \cap \bp_1(p_i) \cdot H\} = 1 \mm{ and thus }\bp_1(0) = 1 .\] Moreover, we know from \ref{haa}  for any $\mu \in (0,1]$, there is a constant
$c(\mu, n)
>0$, so that
 \[ \sup \{|A|(x) \, | \, x \in B_{\mu}(0) \cap H_\infty\} \ge c.\]

Now we claim that the first Neumann eigenvalue $\nu_1(B_{\mu}(0) \cap H_\infty)$ for $P_1$ on $B_{\mu}(0) \cap H_\infty$ is positive:
\[\inf_{f \in C^\infty (B_{\mu}(0) \cap H_\infty), f \not\equiv 0} \int_{B_{\mu}(0) \cap H_\infty} |\nabla f|^2 dA + |A|^2 \cdot f^2 dA \Big/\int_{B_{\mu}(0) \cap
H_\infty} \bp_1^2\cdot f^2 dA > 0 \]  Since $\bp_1$ is positively upper bounded on  $B_{\mu}(0) \cap H_\infty$, it suffices to consider
the usual non-weighted Neumann eigenvalue
\[\nu(B_{\mu}(0) \cap H_\infty) = \inf_{f \in C^{\infty}(B_{\mu}(0) \cap H_\infty), f \not\equiv 0}
\int_{B_{\mu}(0) \cap H_\infty} |\nabla f|^2 dA + |A|^2 \cdot f^2 dA.\]

Clearly, $\nu(B_{\mu}(0) \cap H_\infty) \ge 0$, and for $\nu(B_{\mu}(0) \cap H_\infty) = 0$, we had a smooth positive function $u$ with $\Delta u = |A|^2 \cdot u$ and
vanishing normal derivative along $\p B_{\mu}(0) \cap H_\infty$. But then Stokes theorem says that $\int_{B_{\mu}(0) \cap H_\infty} \Delta u =0$, whereas
$\int_{B_{\mu}(0) \cap H_\infty} |A|^2 \cdot u > 0$.\\

Thus we infer that for sufficiently large $i$:
 \[\nu_1(B_{\mu/\bp_1(p_i)}(p_i)) \ge \nu_1(B_{\mu}(0) \cap H_\infty)/2 > 0\]
 and this contradicts $\nu_1(B_{\mu/\bp_1(p_i)}(p)) \ra 0$. \qed

Now we combine these estimates for  Neumann eigenvalues using particular covers to derive \ref{pi}.\\

{\bf Proof of \ref{pi} and \ref{pia}} \quad We first consider \ref{pi}. We start with the two easy cases: when $H$ is totally geodesic, then $\bp_{\alpha} \equiv 0$ and the Hardy
inequality becomes trivial. Thus we may assume $H$ is connected and not totally
geodesic. Then we have $\bp_{\alpha} > 0$ on $H$.\\

Secondly, when $H$ is regular, we note that $\bp_{\alpha} > 0$ is upper bounded on  $H$, and it suffices to see the positivity of the usual  eigenvalue of $-\Delta + |A|^2$. Indeed,
$|A|^2 \ge 0$ and in some open set it is positive. Thus, for any smooth positive function $u$  (including the first eigenfunction) we have  $\int_H\Delta u =0$,
whereas $\int_H |A|^2 \cdot u > 0$ and therefore the eigenvalue of $-\Delta + |A|^2$ cannot be zero.\\

Now we turn to the main case where $H$ is singular. We notice that for sufficiently small smoothly bounded neighborhoods $W$ of $\Sigma \subset H$:
\[\nu_{\alpha}(H \setminus \overline{W})
> 0.\] This follows as above: for sufficiently small $W$, $H \setminus \overline{W}$ contains an non-empty open ball where $|A|
>0$. For fixed $W$ we have positive bounds for $\bp_{\alpha}$, and we can consider the standard Neumann eigenfunction and the same
Stokes theorem argument as above shows that $\nu_{\alpha}(H \setminus \overline{W}) > 0$. \\

Also we observe: once we have found one neighborhood $W$ with $\nu_{\alpha}(H \setminus \overline{W}) > 0$, the positivity holds for all such neighborhoods $W^*$ of
$\Sigma$ with $W^* \subset W$. Albeit,  this argument does not yet give a uniform positive lower bound while $W^*$ shrinks to $\Sigma$.\\

Now we build a cover of $H \setminus \Sigma$.  $H$ is merely weakly controlled away from $\Sigma$, since the second fundamental form is not enough to control the
geometry in a curved ambient manifold. Thus to get at least a positive bound we use that $\nu_{\alpha}(H \setminus \overline{W}) > 0$ for sufficiently small open
neighborhoods
$W$ of $\Sigma$.\\

Next we use the \textbf{skin adapted covers}  ${\cal{A}}(V_{\alpha, \mu})$ from \ref{besi}, of small neighborhoods $V_{\alpha, \mu} \subset U_{\alpha, \mu}$ of
$\Sigma$ by balls $B_{\mu \cdot \alpha /\bp_\alpha(p)}(p) \subset H \setminus \Sigma$, for $p \in U_{\alpha, \mu} \setminus \Sigma$ with some upper bounded covering numbers $N$.\\

Then we get the following estimate for any $f \in C^{\infty}_c(H \setminus \Sigma), f \not\equiv 0$, writing briefly $W= V_{\alpha, \mu}$ and $B(p)=B_{\mu \cdot \alpha
/\bp_\alpha(p)}(p)$:
\[\int_H |\nabla f|^2 dA +  |A|^2 \cdot f^2 dA \ge  \]
{\small \[1/(N+1) \cdot \left(\int_{H \setminus \overline{W}} |\nabla f|^2  +  |A|^2 \cdot f^2 \, dA +   \sum_{B(p) \in {\cal{A}}} \int_{B(p)} |\nabla f|^2  +  |A|^2 \cdot f^2 \, dA
\right) \ge
\] \[ 1/(N+1) \cdot \left( \nu_{\alpha}(H \setminus \overline{W}) \cdot \int_{H \setminus \overline{W}}\bp_{\alpha}^2 \cdot  |f|^2 \, dA + \sum_{B(p) \in
{\cal{A}}}\nu_{\alpha}(B(p)) \cdot \int_{B(p)} \bp_{\alpha}^2 \cdot  |f|^2 \, dA  \right)\]}
\[ \ge 1/(N+1) \cdot \min \{\nu_{\alpha}(H \setminus \overline{W}), \inf_{B(p) \in {\cal{A}}} \nu_{\alpha}(B(p))\} \cdot \int_H \bp_{\alpha}^2 \cdot |f|^2 dA\]

Now we use \ref{rad} (\ref{rade}) saying that all the Neumann eigenvalues $\nu_{\alpha}(B(p))$ are uniformly lower bounded by $\nu_\alpha(B_{\mu \cdot \alpha /\bp_\alpha(p)}(p))  \ge \zeta$. Thus,
we get
$$\lambda_{P_\alpha} \ge \frac{\min \{\nu_{\alpha}(H \setminus \overline{W}), \zeta\}}{N+1}  > 0,$$
In other words, the Hardy inequality holds for  $k_{\alpha ,H} :=\lambda_{P_\alpha}$.\\

To verify \ref{pia} we review the previous argument for a complete area minimizing hypersurface $H^n \subset \R^{n+1}$. Then we get skin adapted covers not only of small
neighborhoods of $\Sigma$, but of the entire hypersurface cf.\ref{besi}. Hence, the previous chain of inequalities collapses to
\[\int_H |\nabla f|^2 dA +  |A|^2 \cdot f^2 dA \ge  \]
\[1/N \cdot \sum_{B(p) \in {\cal{A}}} \int_{B(p)} |\nabla f|^2  +  |A|^2 \cdot f^2 \, dA \ge
1/N  \cdot   \sum_{B(p) \in {\cal{A}}}\nu_{\alpha}(B(p)) \cdot \int_{B(p)} \bp_{\alpha}^2 \cdot  |f|^2 \, dA. \] Thus, in this case, the Hardy inequality holds for  $k_{\alpha ,n} :=\zeta/N$.
\qed

\begin{remark} \label{motii} \quad There is \textbf{no} similarly controlled cover with not too small balls  (= with lower bounded Neumann eigenvalues for $\Delta +|A|^2$) when we use
$dist(\cdot,\Sigma)$ in place of $1/\bp_\alpha$, since the result relies on the fact that $\bp_\alpha \ge |A|$ and thus $\bp_\alpha$ also controls the local complexity of the
underlying geometry.\\

Thus the Hardy inequality \ref{piw} for $1/dist(\cdot,\Sigma)$ is not a parallel, but a proper consequence, of the Hardy inequality for $\bp_\alpha$.\qed
\end{remark}

\subsubsection{Whitney Smoothings}\label{skinada1}
\bigskip

Here we discuss methods to improve the regularity properties of skin transforms and of \sks .\\

 The axiom (S4) asserts that any skin transform $\bp$ is locally Lipschitz regular. In some cases we would rather like to have a proper
smoothness. An instance is the Ancona-Martin theory we consider in [L1], where we use smoothed versions of $\bp$.\\

In turn, also the \sks\ have some basic regularity properties. In the case of minimizing cones an appropriate implicit function theorem shows
that the \sks\ are locally Lipschitz submanifolds. Albeit, in general, there is no global estimate for the Lipschitz constant.\\

Also we have seen that the \sks\ $\M_c$ may satisfy particular additional properties like an exterior ball property in the case of the $\bp_\alpha$ or they may have non-negative
mean curvature, as for $\bp_\square$. But we may need a higher regularity to use the \sks\ effectively.\\

However, we can hardly improve both regularity properties, that of $\bp$ and its \sks , at the same time. But we can use skin adapted covers to either improve the regularity
of $\bp$ or of its \sks . This way we get two complementary descendants with different characteristics.\\

The idea really is to use ${\cal{A}}$ to define a \textbf{discretization} of $\bp$,  respectively of $\delta_{\bp}$, which can be reassembled in two complementary ways.

\begin{itemize}
    \item \textbf{$\delta_{\bp} \rightarrow \delta^*_{\bp}$} \quad We merely keep the values of $\delta_{\bp}$ for the center points $A$ of the cover ${\cal{A}}$. Then we extend them over the
        balls $B_{\Theta(p)}(p)$ of a skin adapted cover, cut them off by some smooth cut-off function when we leave $B_{2 \cdot \Theta(p)}(p)$, $p \in A$, and add them all up
        to get a smooth function that approximates $\delta_{\bp}$.
    \item \textbf{$\M_c \rightarrow \M_c^*$} \quad For the complementary problem, the regularity of the \sks\ and the domain they bound, we
    consider a specific subcollections ${\cal{A}}[c]$ of balls $\overline{B_{\Theta(p)}(p)}$  in
        ${\cal{A}}$ so that  $\{x \in H \setminus \Sigma \,|\, \delta_{\bp}(x) \ge c \} \subset \bigcup_{{\cal{A}}[c]}  \overline{B_{\Theta(p)}(p)}$. Then we can define new and better controlled \sks\  by \,
        $\M^*_c := \p \, (\bigcup_{{\overline{B_{\Theta(p)}(p)} \in \cal{A}}[c]} B_{2 \cdot \Theta(p)}(p)).$\\
\end{itemize}

Here we focus on the smooth approximations $\delta_{\bp} \rightarrow \delta^*_{\bp}$. We return to the second way of discretization in Ch.\ref{ch}, when we define skin uniform domains.

\begin{proposition}\label{smsk} \emph{\textbf{(Smoothings)}} \, For any skin transform $\bp$ there is a smoothed version
$\bp^*$ on $H \setminus \Sigma$, so that \begin{itemize}
    \item \quad $c_1 \cdot \delta_{\bp}(x) \le \delta_{\bp^*}(x)  \le c_2 \cdot \delta_{\bp}(x)$
    \item \quad $|\p^\beta \delta_{\bp^*}  / \p x^\beta |(x) \le c_3(\beta) \cdot \delta_{\bp}^{1-|\beta|}(x)$,
\end{itemize}
for any $x \in H \setminus \Sigma$. $\beta$ denotes the usual multi-index for derivatives, with respect to normal coordinates around $x \in H \setminus \Sigma$, and $c_i >0$,
$i=1,2,$ are constants depending on $H$, $c_3 > 0$ depends on $H$ and on $\beta$.\\
\end{proposition}

 {\bf Proof } \quad We choose a smooth non-negative function $\phi$ on $\R^n$ with $\phi \equiv 1$ on $B_1(0)$ and $\phi \equiv 0$ on $\R^n  \setminus B_2(0)$.
 We consider, in \ref{besi} Step 2.2.,  $\exp_p[10 \cdot \Lambda(\zeta) \cdot \bp(p) \cdot H]$, for some tiny $\zeta >0$.\\

For $x \in H \setminus \Sigma$ we set
\[\Phi_{p}(x):=  \phi(\exp^{-1}_p[10 \cdot \Lambda(\zeta) \cdot \bp(p) \cdot H](x))\]
Thus the $exp_p$-preimage of the ball $B_5(p) \subset 10 \cdot \Lambda(\zeta) \cdot \bp(p) \cdot H$ is almost isometric to $B_5(0) \subset T_pH$.\\

We notice that $|\p^\beta \Phi_p/ \p x^\beta |(x)\le k(\beta) \cdot (10 \cdot \Lambda(\zeta) \cdot \bp(x))^{|\beta|-1}$ on $B_5(p) \subset 10 \cdot \Lambda(\zeta) \cdot \bp(p)
\cdot H$. Now we define $\bp^*$ through its $\bp$-distance
\[\delta_{\bp^*}(x):= \delta^*_{\bp}(x):= \sum_{p \in A} \delta_{\bp}(p) \cdot \Phi_{p}(x)\]
This is obviously a smooth function and the claimed properties readily follow from the upper bounded covering number of ${\cal{A}}$.\\

$\bp^*$ satisfies the axioms (S1)-(S4). They are readily verified from those we presume for the skin transform $\bp$. But due to the unavoidable breaking of symmetry to define $\bp^*$  it is only \emph{quasi-natural}, in the sense, that the inherited quantities under convergence, equal the genuine one up to fixed positive multiples. This follows from (S5) for $\bp$ and the inequalities $c_1 \cdot \delta_{\bp}(x) \le \delta_{\bp^*}(x)  \le c_2 \cdot \delta_{\bp}(x)$.\qed

\begin{remark}\textbf{(Whitney smoothings)} \label{sws} \,  This smoothing process adapts the idea of Whitney smoothings for the Lipschitz regular
distance function to closed sets $A \subset \R^n$, like, $\p D$ for a domain $D \subset \R^n$, [St], cf. Ch.6.2. There are also skin counterparts to the extension results by
Whitney in [Wh1]. The idea is to insert the cover ${\cal{A}}$ in the extension operators defined by Stein in [St], Ch.6. We leave details
 to the reader.
 \qed
\end{remark}

 \setcounter{section}{4}
\renewcommand{\thesubsection}{\thesection}
\subsection{Non-Tangential Accessibility of $\Sigma$}
\bigskip

The global analysis on spaces with boundary is significantly influenced from the coupling of their interior shape with that of their boundary. The accessibility of boundary
points from the interior (or the exterior) plays a vital role in the elliptic theory on these spaces. For the classical case of Euclidean domains equipped with the Laplacian this is a
broad topic in potential theory, cf.[AG],Ch.8 or [Do],Ch.XII.\\

Here we study the question whether $\Sigma \subset H$ can be non-tangentially approached from its complement through paths. There are
two interpretations of this question, an extrinsic and an intrinsic one. \\

From a classical viewpoint, to study the analysis on $M \setminus H$, we try to reach $H$, and in particular its most delicate constituent $\Sigma \subset H$, from
$M \setminus H$. We will see that the minimality of $H$ causes an extrinsic non-tangentially accessibility of $H$. \\

However, our main focus are potential theoretic and elliptic problems on $H$. For this, we consider $\Sigma$ as a boundary of $H \setminus \Sigma$ and ask whether $\Sigma$
can be non-tangentially approached from within $H \setminus \Sigma$. This intrinsic question introduces another twist since $H$ degenerates while we approach
$\Sigma$.\\

It is the main goal of this chapter to establish a strong form of such an intrinsic non-tangential accessibility of $\Sigma$, measured in terms of $\bp$-distances.
This is the skin uniformity of $H \setminus \Sigma$. In [L1] this will be the essential basis to establish a natural Gromov hyperbolic geometry on $H \setminus \Sigma$ and  to understand the potential theory of typical elliptic operators on $H \setminus \Sigma$, for instance, in terms of a Martin theory.

\subsubsection{Extrinsic Accessibility}\label{exnta}
\bigskip

We start with the question whether $\Sigma$ can be extrinsically approached in a non-tangential way from $M \setminus H$. A basic way to describe such an approachability
quantitatively is that of inner cone conditions.

\begin{definition} \emph{\textbf{(John Domains)}} \label{jdo} \quad  A domain $D \subset \R^n$ is called a \textbf{John domain}, if, for some basepoint $p_0 \in D$, any $p \in D$ can be linked by a path
$\gamma_p \subset D$ to $p_0$, so that the length of the subarc $\gamma_p(z)$ from $p$ to any $z \in \gamma$ can be estimated by a fixed multiple $c >0$ of the distance to the boundary:
\[l(\gamma_p(z)) \le c \cdot dist(z,\p X), \mm{ for any  } z \in \gamma_p.\]
In other words, there are \textbf{twisted cones} of lower bounded width in $D$ pointing from $p_0$ to any $p \in D$.
\end{definition}

With this concept we can give a first indication for the relevance of non-tangential accessibility problems
in the context of area minimizers. The following is a remarkable result due to David and Semmes [DS](1.8) and (1.10), which, in particular, says that any $p \in \Sigma$ can be approached
by a twisted cone in $M \setminus H$.
\begin{proposition} \label{ds} \quad  Let $H \subset \R^{n+1}$ be a hypersurface that bounds two disjoint domains $H^+ \cup H^-=\R^n \setminus H$.
Then we have: \[ H^+ \mm{ and } H^- \mm{ are John domains } \Leftrightarrow H \mm{ is a quasi-minimizer of the area functional.}\]
\end{proposition}

\begin{remark} \label{dspr} \quad  We heuristically simplified the original statements, cf. [DS], Ch.1 for the proper details. The authors
use a non-standard, slightly stronger, definition of John domains with some additional Lipschitz-constraints on the paths, cf.[DS](1.7).\\

Here, a hypersurface means an Ahlfors-regular set of dimension $n-1$, and a quasi-minimizer is a minimizer of a generalized area functional.
In particular, area minimizing hypersurfaces are also quasi-minimizers. We neither use these concepts nor the result \ref{ds} in this series of papers. Therefore, we refer to [DS] for all further details.\qed
\end{remark}

In general,  John domains show pathologies like the occurrence of several minimal Martin boundary points in one given point in $\p D$. They result from a missing control over the length of curves in $D$ relative to the metric distance between their endpoints.\\

In turn the sharpened concept of uniform domains, which precisely takes care of this length control, resolves these issues.
The notion of uniform domains was actually independently suggested from different sources
 like Sobolev-theory, quasi-conformal geometry or Martin theory. Also there is a satisfactory extension of the concept to uniform metric spaces, cf.[Ai],[He],[V] or [BHK], Ch.1.\\

For our study of $H \setminus \Sigma$ the language of uniform spaces is very natural. Therefore we recall the definition of such spaces.\\

For a metric space $(X,d)$ we recall that the \textbf{length} of a path or arc, meaning a continuous map, $\gamma:[0,1] \ra X$, is defined as {\small \[l(\gamma):= \sup \Big\{
\sum_{i=0,..,N}
        d(\gamma(t_{i-1}),\gamma(t_i))\,\Big|\, \mm{ all partitions } 0 =t_0 \le t_1 ... \le t_N = 1\Big\}.\]}
$\gamma$ is called \textbf{rectifiable}, when $l(\gamma) < \infty$. $X$ is \textbf{rectifiably connected}, when each pair of points  in $X$  can  be  joined  by a
        rectifiable path.\\
\begin{definition} \emph{\textbf{(Uniform Spaces)}} \label{us} \quad  For a locally compact, non-complete, locally complete, rectifiably connected metric space  $X$
we set $\p X:= \overline{X} \setminus X$, where
        $\overline{X}$ is the metric completion of $X$.\\

         Such a metric space $X$ is called a  \textbf{uniform space}, more precisely a c-uniform space,
        if there is a constant $c \ge 1$ such that any two points can be joined by a c-uniform curve in $X$.\\

\noindent A \textbf{c-uniform} curve joining $p,q \in X$ is a rectifiable path $\gamma: [a,b] \ra X$, for some $a <b$, running from $p$ to $q$ so that $\gamma$ satisfies the
following two axioms\begin{itemize}
    \item \emph{\textbf{Quasi-Geodesic}:} \quad \quad  $l(\gamma)  \le c \cdot  d(p,q),$
    \item  \emph{\textbf{Twisted Double Cones}:}\quad  $l_{min}(\gamma_{p,q}(z)) \le c \cdot dist(z,\p X),$
\end{itemize}
 for any  $z \in \gamma_{p,q}.$ $l_{min}(\gamma_{p,q}(z))$ denotes the minimum of the lengths of the subarcs of $\gamma_{p,q}$ from $p$ to $z$ and from $q$ to $z$.\\
\end{definition}

In this framework, a \emph{uniform domain} $D \subset \R^n$ is a uniform space for the  Euclidean metric on $D$. If we choose to use the inner metric, that is the infimum over
 the length of connecting paths in $D$ we get the concept of an \emph{inner uniform domain}.\\

\bigskip

\subsubsection{Intrinsic Accessibility and Skin Uniformity}\label{inta}
\bigskip

\textbf{Skin Uniform Spaces} \quad Now we reach our main objective in this chapter, to understand the \emph{intrinsic accessibility} of $\Sigma$ from within $H \setminus \Sigma$. Thus
we view $\Sigma$ as $\p X$ for the metric space $X = H \setminus \Sigma$ and ask whether the axioms for uniform spaces are satisfied.
Actually, this is the case, but we need to go still one step further and first establish a \emph{skin uniformity} of $X$ to conclude the uniformity.\\

To better understand this strategy, we observe that for $\Sigma$ relative $H \setminus \Sigma$ the problem of accessibility
now also involves the diverging curvature of $H \setminus \Sigma$ when we approach $\Sigma$. This becomes an issue for both uniformity conditions since the complexity and curvature of $H \setminus \Sigma$ are not properly coupled to the metric distance to $\Sigma$.\\

Therefore it seems rather delicate to approach the proof of the desired metric twisted cone condition
\[l_{min}(\gamma_{p,q}(z)) \le a \cdot dist(z,\p X)\] directly. Instead, as in the case of the Hardy inequalities above, we need to prove a sharpened skin version for a given skin
transform $\bp$:
\[l_{min}(\gamma_{p,q}(z)) \le b \cdot \delta_{\bp}(z),\] and, a posteriori, we infer the result also for $dist(z,\p X)$, from the general relation $\delta_{\bp}(z) \le c \cdot dist(z,\p X)$.\\

We point out that the use of a skin version of the twisted cone condition is much more than a technicality. Our intended
analytic applications rely on this stronger skin version (not on the usual metric variant) of uniformity.\\

To formulate our main result we fix some skin transform $\bp$. Also for our notational convenience we assume that the Lipschitz constant for $\delta_{\bp}$ equals $1$.

\begin{proposition} \emph{\textbf{(Skin Uniformity of  $\mathbf{H \setminus \Sigma}$)}} \label{intsuni} \quad  For any hypersurface
$H \in {\cal{H}}^n$ with singular set $\Sigma_H$, which may also be empty, we have

 \begin{enumerate}
 \item$H \setminus \Sigma$ and $H$ are \textbf{rectifiably connected}. In particular, any compact $H \in {\cal{H}}$ has a finite intrinsic diameter: $diam_{g_H}H   < \infty.$
\item $H \setminus \Sigma$  is a $c$-\textbf{skin uniform space}, for some $c >0$. That is, any pair $p,q \in H \setminus \Sigma$ can be joined by a \textbf{c-skin uniform
    curve} in $H \setminus \Sigma$, i.e. a rectifiable path $\gamma_{p,q}: [a,b] \ra H \setminus \Sigma$, for some $a <b$, with $\gamma_{p,q}(a)=p$,
    $\gamma_{p,q}(b)=q$, so that the following  skin uniformity conditions hold
    \begin{itemize}
    \item \emph{\textbf{Quasi-Geodesic}:} \quad \quad  $l(\gamma)  \le c \cdot  d(p,q),$
    \item  \emph{\textbf{Twisted Double Skin Cones}:}\quad $l_{min}(\gamma_{p,q}(z)) \le c \cdot \delta_{\bp}(z)$,
\end{itemize}
for any $z \in \gamma_{p,q}.$
\item Any pair $p,q \in H$, that is, we also allow singular endpoints, can be joined by a c-skin uniform curve supported in $H \setminus \Sigma$, except for its  endpoints.
 \item For $H \in {\cal{H}}^{\R}_n$, that is, for  $H \subset \R^{n+1}$,  we get common constant $c_n$ depending only on the dimension, so that
 $H \setminus \Sigma$  is a $c_n$-skin uniform space.
\end{enumerate}
\end{proposition}

The proof is postponed to the next section. In this section we continue with some remarks, immediate consequences and byproducts of the proof.\\

In simple terms, the connectivity follows from the isoperimetric inequality for area minimizers and we exploit the naturality of $\delta_{\bp}$ to derive the skin uniformity from
 this connectedness of $H \setminus \Sigma$ based on compactness results for area minimizers.\\

In turn $dist(\cdot,\Sigma)$ does not behave naturally under convergence of the underlying spaces and looses its meaning when $\Sigma \v$.
Hence the strategy used in skin uniform case does not apply. Nevertheless,  skin uniformity trivially implies usual uniformity and therefore we get

\begin{corollary} \emph{\textbf{(Uniformity of  $\mathbf{H \setminus \Sigma}$)}} \label{intuni} \quad For any connected compact area minimizing hypersurface
$H \subset M$ with singular set $\Sigma \n$, we have:
\[ \mm{ For }  X= H \setminus \Sigma, \p X=\Sigma, \, X  \mm{ is a }d\mm{-\textbf{uniform space}, for some }d >0.\]

\end{corollary}

We also notice that the same results hold for the minimal, although not area minimizing, hypersurfaces $S_C \in {\cal{K}}$,  where
${\cal{K}}_{n-1}:= \{ S_C\,| \, C \in {\mathcal{C}_{n}}\}$ and
${\cal{K}}= \bigcup_{n \ge 2} {\cal{K}}_{n-1}$. In particular, we get

\begin{corollary} \emph{\textbf{(Extension to Hypersurfaces $\mathbf{S_C \in {\cal{K}}}$)}}
 \label{ext} \quad For any hypersurface $S_C \in {\cal{K}}$ with singular set $\Sigma \n$, there is some $c >0$, so that
\begin{itemize}
\item $S_C \setminus \Sigma$ and $S_C$ are rectifiably connected.
\item $S_C\setminus \Sigma$  is a $c$-skin uniform space and, thus, $c$-uniform space.
\end{itemize}
\end{corollary}

The intrinsic properties of $H \setminus \Sigma$ asserted in \ref{intsuni} and \ref{intuni} clearly rely on the extrinsic property of $H$ as an area minimizer in its ambient
space $M$.\\

 From this, we also get extrinsic estimates when $H$  is an oriented boundary of  some open set $A \subset \R^{n+1}$, that is, $H \in
{\cal{H}}^{\R}_n$ and, completely similarly for $S_C \in {\cal{K}}_{n-1}$ since it  bounds an open and connected set in $S^n$, as follows from the fact that any
$C \in {\mathcal{C}_{n}}$ bounds an open and connected set in $\R^{n+1}$, cf. [BG],Th.1.

\begin{lemma} \emph{\textbf{(Intrinsic versus Extrinsic Metric)}} \label{fini} \quad There are constants $c^{\R}_n \ge 1$ and  $c^{S}_n \ge 1$ depending only on the dimension $n$, so that\begin{itemize}
    \item  For any  $H \in
{\cal{H}}^{\R}_n$ and any two points $p,q \in H \subset \R^{n+1}$ we have
 \[ d_{g_H}(p,q)  \le  c^{\R}_n \cdot  d_{g_{\R^{n+1}}}(p,q) \le c^{\R}_n \cdot d_{g_H}(p,q).\]
    \item  For any  $S_C \in {\cal{K}}_{n-1}$ and any two points $p,q \in S_C \subset S^n$ we have
 \[d_{g_{S_C}}(p,q)  \le c^{S}_n \cdot  d_{g_{S^{n}}}(p,q) \le c^{S}_n  \cdot d_{g_{S_C}}(p,q).\]
\end{itemize}
\end{lemma}

The estimate $d_{g_{S_C}}(p,q) \le c^{S}_n \cdot  d_{g_{S^{n}}}(p,q)$ readily implies the following uniform diameter estimate.

\begin{corollary} \label{finii} \quad  $diam(S_C,g_{S_C}) \le c^{S}_n \cdot diam(S^n) < \infty$, for any $S_C \in {\cal{K}}_{n-1}$.\\
\end{corollary}

Occasionally, it is important to see the relation between the (skin) uniformity constants of a hypersurface and of its blow-up geometries. This is
easily accomplished from the naturality property of $\bp$ and can be stated as follows

\begin{proposition} \emph{\textbf{(Inheritance under Blow-Ups)}} \label{insu} \quad Let $H$ be some complete area minimizing hypersurface, so that $H \setminus \Sigma$ is a $c$-skin uniform space, for some $c >0$.
 Then, for any blow-up limit $F$, $F \setminus \Sigma_F$ is also $c$-skin uniform. In particular, this applies to all tangent cones of $H$.
\end{proposition}

 {\bf Proof} \quad For any two points $p,q \in F \setminus \Sigma_F$, we first find a rectifiable path $\gamma_{p,q}: [a,b] \ra F \setminus \Sigma_F$,
 for some $a <b$, with $\gamma_{p,q}(a)=p$ and  $\gamma_{p,q}(b)=q$, so that \begin{itemize}
    \item \quad  $l(\gamma_{p,q})  \le (c+\ve ) \cdot d_{g_C}(p,q)$,
    \item \quad $l_{min}(\gamma_{p,q}(z)) \le (c+\ve ) \cdot \delta_{\bp}(z)$, for any $z \in \gamma_{p,q}.$
\end{itemize}

This follows from the fact that any compact portion of $F \setminus \Sigma_F$ admits an $\ve$-fine  $C^3$-approximations by suitable pieces of $k \cdot H$ for some
large $k$. As in Step 2 above we can infer the existence of $\gamma_{p,q}$ from that on $k \cdot H$.\\

Namely, when we have scaled $H$ by a sufficiently large constant so that not only the
$\D$-images of $p$ and $q$ can de identified with points in $H$ but also the $c$-skin uniform path that joins them in $H \setminus \Sigma$ is the $\D$-image we simply
take the preimage of this path in $C$ and notice that it is $(c+\ve )$-skin uniform, once we scaled $H$ appropriately, depending on the chosen $p, q$ and $\ve >0$.\\

Finally, we send $\ve \ra 0$ and apply some BV-compactness results,  like Helly's selection principle, cf. [SG],4.5, to a sequence $\gamma_n$ of such $c + 1/n$-skin uniform curves to get a limit curve
that is $c$-skin uniform. \qed\\

We already noted in \ref{skino} and \ref{ext} that for any area minimizing cone $C$, $\bp_{S_C}:= \bp_C|_{S_C}$, with $S_C:= \p B_1(0) \cap C \subset S^n \subset
\R^{n+1}$, satisfies the the skin axioms and $S_C \setminus \Sigma_{S_C}$ is skin uniform.

\begin{corollary}  \label{insus} \quad For a complete singular area minimizing hypersurface
$H \subset M$ with singular set $\Sigma \n$,  let $H \setminus \Sigma$ be a $c$-skin uniform space, for some $c >0$. Then, for any of its tangent cones $C$, \[S_C \setminus
\Sigma_{S_C}\mm{ is } d\mm{-skin uniform, for some } d(c,n) >0.\]
In particular, we get a  constant $c^*_n$ depending only on the dimension, so that
 $S_C \setminus \sigma_C=\p B_1(0) \cap C \setminus \sigma_C$  is a $c^*_n$-skin uniform space relative $\bp:=\bp_C|_{S_C}$,  for any $C \in {\cal{C}}_n$.
\end{corollary}

 {\bf Proof} \quad Given  two points $p,q \in S_C \setminus \Sigma_{S_C}$, we use that $C \setminus \sigma$ is $c$-skin uniform and the fact that the radial projection of
 curves in $(C \setminus \sigma) \setminus B_1(0)$ onto $(C \setminus \sigma)\cap \p B_1(0)$ decreases the length to  find a rectifiable path \[\gamma_{p,q}: [a,b] \ra C \setminus \sigma \cap \overline{B_1(0)}, \mm{
 for some } a <b,\]
  with $\gamma_{p,q}(a)=p$ and  $\gamma_{p,q}(b)=q$, so that \begin{itemize}
    \item \quad  $l(\gamma_{p,q})  \le c \cdot d_{g_C}(p,q)$,
    \item \quad $l_{min}(\gamma_{p,q}(z)) \le c \cdot \delta_{\bp}(z)$, for any $z \in \gamma_{p,q}.$
\end{itemize}

Since \ref{finii} says that $diam(S_C,g_{S_C}) \le s_n:= c^{S}_n \cdot diam(S^n) < \infty$, for any $S_C \in {\cal{K}}_{n-1}$, we have  $l(\gamma_{p,q}) \le s_n$. On the other
hand,  $\delta_{\bp}(p) \le L_{\bp} \cdot dist(p,\sigma_C)$ shows that there is some $\eta(c) \in (0,1)$, so that for any $C$: \[\delta_{\bp}(z) \le c/10 \mm{ on } B_\eta(0) \cap C \mm{ and, hence, }
\gamma_{p,q} \subset (B_{s_n}(0) \setminus B_\eta(0)) \cap C.\] Now we radially project $\gamma_{p,q}$ on $S_C$ and readily verify its $d$-skin uniformity for some $
d(c,n) >0$.
 \qed

\bigskip
\subsubsection{From Isoperimetry to Uniform Networks}\label{intapr}
\bigskip

In this section we present a proof of \ref{intsuni}. The other results will be natural byproducts, consequences or intermediate steps. As usual the totally geodesic case is trivial and
we may assume that our hypersurfaces are non-totally geodesic. In particular, we can use that $\bp >0$ and that $\delta_{\bp}$ is Lipschitz continuous.\\

\textbf{Proof of \ref{intsuni}} \quad We build skin uniform curves in several steps. The strategy is to show that $H \setminus \Sigma$ is rectifiably connected and then to
stepwise
improve this to a quantitative connectivity result. We first get short quasi-geodesics arcs, assemble them in pipelines and finally get the desired skin uniform arcs.\\

The main ingredients are the \textbf{isoperimetric inequality} for area minimizers  and the \textbf{self-similarity} of the pair $(H,\bp_H)$, that is, the scaling invariance of the minimizing condition for $H$ and the naturality of $\bp$.\\

\noindent \textbf{Step 1 (Connectivity of  $H \setminus \Sigma$)} \label{rn}\\

We start from a \textbf{connected} area minimizing hypersurface $H$ which, at this stage, means a integer multiplicity rectifiable current.\\

 Here we want to understand finer connectivity
properties and prove that $H \setminus \Sigma$ is rectifiably connected.\\

Although the codimension of $\Sigma$ is rather large, this is not evident.  The partial regularity theory merely represents $H \setminus \Sigma$ as a union of open $C^1$-regular subsets of $H$ resulting from non-constructive
selection processes using e.g. the theorems of Lusin and Egoroff.  Also $H$ degenerates towards $\Sigma$ and therefore it is by no means  clear whether deleting $\Sigma$ could
harm the connectivity of $H$.

\begin{lemma} \emph{\textbf{(Connectivity  of $\mathbf{H \setminus \Sigma}$)}} \label{coh} \quad $H \setminus \Sigma$ is rectifiably connected.
\end{lemma}

{\bf Proof} \quad Connected Riemannian manifolds are rectifiably connected. They are path connected since they are locally path-connected and moreover each continuous
path in a manifold can be approximated by a rectifiable one. Thus we may assume $\Sigma \n$ and note that the remaining assertion we need to check
is that $H \setminus \Sigma$ is connected.\\

Let us assume that $H \setminus \Sigma$ contains (at least) two open, non-empty and disjoint path components. We consider one such path component $H_1 \subset H
\setminus \Sigma$ and its complement $H_2:= (H \setminus \Sigma) \setminus H_1$. The idea is to think of $\p {H_i} \subset \Sigma$ as the boundary of the minimal current
$H_i$ and to see that even a local disconnectness of $H \setminus \Sigma$ is in conflict with the isoperimetric inequality. \\

Technically, we want to reduce the problem to the case of oriented boundaries. To this end, we use the local decomposition of the rectifiable current $H$ into a locally disjoint collection of oriented minimal boundaries described in the appendix, part IV.  However, since we are about to show that $H \setminus \Sigma$ is connected, this cannot be a global decomposition.\\

 The way out is to prove the following stronger localized claim for the case where $H$ is an oriented boundary. Actually, cf. \ref{dicc} and its discussion, for any $p \in \Sigma$, there is an $r_p >0$ so that $B_{r_p}(p) \cap H \subset M$ is an oriented boundary within $B_{r_p}(p)\subset M$.\\

\textbf{Local Connectedness} \quad \emph{For any $p \in \Sigma$ and $r \in (0,r_p)$, we choose the connected component $B^H_r(p) \subset B_r(p) \cap H$ that contains $p$.}
\begin{center} \emph{Then, for any sufficiently small $r>0$, $B^H_r(p) \setminus \Sigma$ is still connected.} \end{center}

 In turn, since $H$ is connected, this property readily implies that $H \setminus \Sigma$ is connected also when $H$ is a current.  Thus, for the remainder of this argument, which amount to prove this local connectedness, we may assume that $H$ is an oriented boundary.\\

Assume there are arbitrarily small $r \in (0,r_p)$, so that $B^H_r(p) \setminus \Sigma$ contains at least two path components
 $H_1(r) \subset B^H_r(p) \setminus \Sigma$ and its complement $H_2(r):= (B^H_r(p) \setminus \Sigma) \setminus H_1(r)$.\\

Then we can rescale $B_r(p) \subset M$ to unit size and may assume that the resulting geometry on $B_1(p) \subset r^{-1} \cdot M$ is as near (in $C^5$-topology)
to that of a Euclidean ball as we want. \\

For this rescaled ball $B_1(p)$, we choose a narrow tube $U(1)$ of $\p {H_1}(1)$ of small volume.
This is possible since the Hausdorff-dimension of $\Sigma$ is $\le n-7 \le n-2$: we get a ball cover $B[1],..,B[j]$ of $\Sigma$ in $M$ with arbitrarily small total area
 $Vol_n(\bigcup_{k=1,..,j} \p B[k]$. The minimality of $H$ shows that $Vol_n(H \cap \bigcup_{k=1,..,j} B[k]) \le  Vol_n(\bigcup_{k=1,..,j} \p B[k])$.\\

Also we can make sure that $\p U(1) \cap H_2(1)$ is smooth with $Vol(\p U(1) \cap H_2(1)) \ll 1$. For instance, one can apply the coarea formula e.g. [GMS], Vol1, Ch.2.1.5, Th.3, p.103 to some mollified distance function see the existence of such tubes $U(1)$.\\

Then we consider $H^*_1:=  H_1(1) \cup U(1)$ and $H^*_2:= (B_1(p) \setminus \Sigma) \setminus H^*_1$. They are both integral currents and we have a decomposition
\[B_1(p) = H^*_1 \cup H^*_2 \, \mm{ with }\, Vol(B_1(p)) = Vol(H^*_1) + Vol(H^*_2).\]

Now a variant of the \textbf{isoperimetric inequality} for oriented minimal boundaries due to Bombieri and Giusti, cf.[BG], Th.2, p.31 says
\[Vol_{n-1}(\p H^*_1 \cap  B_1(p))   \ge  k_n \cdot \min \{Vol_n(H^*_1 \cap  B_{\beta_n}(p)), Vol_n(H^*_2 \cap  B_{\beta_n}(p))\}^{n-1/n},\]
for some constants $k_n > 0, \beta_n \in (0,1)$ depending only on the dimension.\\

But this leads to a contradiction since the right hand side is positively lower bounded whereas the left hand side becomes arbitrarily small when we shrink
the tube $U(1)$ closer towards $\Sigma$ so that $Vol_{n-1}(\p U(1) \cap H_2(1)) \ra 0$.  \qed

From this we observe a refinement, which will be used in Step 4, in the case where $H$ is compact. Here, and for later use, we introduce the following notation
 \begin{equation}\label{iii}
 \E(\rho):= H \setminus \I(\rho) \mm{ and } \I(\rho):= \{x \in H \setminus \Sigma \,|\, \delta_{\bp}(x) < \rho \},\mm{ for any } \rho >0.
 \end{equation}

\begin{corollary} \label{rnpc} \quad  Let $H \in {\cal{H}}^c_n$ be a compact hypersurface and $K \subset H \setminus \Sigma$ any compact subset. Then there are $l=l(K) >0$, $s = s(K) > 0$ so that, any two points $a,b \in K$, can be linked
through a rectifiable curve in $\E(s)$ of length $\le l$.
\end{corollary}

{\bf Proof} \quad Indeed this is a finite problem: we can cover $K$ by a finite collection of small balls $B_\rho(p_1), ..,B_\rho(p_k)$, $\rho \ll 1$, so that $\overline{B_{2 \cdot
\rho(p_1)}}, ...\overline{B_{2 \cdot \rho(p_k)}} \subset H \setminus \Sigma$.\\

Now we first link both $x$ and $y$ to the center of one the balls they belong to, then we link any two of these centers by a rectifiable curve in $H \setminus \Sigma$. But
these are merely finitely many, namely $k \, !$, compact sets $\Gamma_1, ..\Gamma_{k \, !}$ in $H \setminus \Sigma$.\\

 Therefore we can find some small $s >0$ so that due to
the compactness of the union all the curves and balls $\overline{B_{2 \cdot \rho(p_i)}}$.
 \[\overline{B_{2 \cdot \rho(p_1)}} \cup ...\overline{B_{2 \cdot \rho(p_k)}} \cup \Gamma_1 \cup  ..\Gamma_{k \, !} \subset \E(s)\]
 and we can choose \[l(K) := \max \{l(\Gamma_1), ..l(\Gamma_{k \, !})\} +2.\]
 \qed

\noindent  \textbf{Step 2 (Short Quasi-Geodesic Arcs)} \label{co1}\\

Here we derive the existence of short quasi-geodesic arcs with some controlled $\bp$-distance. We start with the case we get after strong scaling $H \subset M$.

\begin{lemma} \label{copr} \quad Let $H^n  \subset \R^{n+1}$ be a complete area minimizing hypersurface which is an oriented boundary of an open set $A  \subset \R^{n+1}$. Then, for any $t >0$, there is some $\tau(t,n) < t$ and an  $L(t,n)
> 0$,  so that any two points  $p,q \in  \E(t)  \subset H^n$ with $d_{g_{Eucl.}}(p,q) =1$  can be connected by a path:
\begin{equation} \label{eq}\gamma_{p,q} \subset   \E({\tau}) \mm{ of length } l(\gamma_{p,q}) \le L(t,n).\end{equation}
\end{lemma}

 This is a quantitative version of the connectivity result \ref{coh}.  Again, we consider oriented boundaries to derive an intrinsic result also for rectifiable currents in Step 4b below. In the present case, working
 with oriented boundaries allows us to prevent extinctions of minimizing pieces when we apply compactness arguments to derive the asserted estimates.\\

\textbf{Proof} \quad Since $H^n \subset \R^{n+1}$ is an oriented boundary  it is connected and  bounds a connected open set $U_H  \subset \R^{n+1}$, cf. [BG], Ch.2, Th.1 and its Corollary. \\

In particular,  we infer, from Step 1, that $H \setminus \Sigma_H$ is rectifiably connected\\

Now we assume there are \emph{no} estimates as in (\ref{eq}), valid for all hypersurfaces. Then there are a sequence of such hypersurfaces $H_i^n  \subset \R^{n+1}$, bounding
 connected open sets $U_{H_i}  \subset \R^{n+1}$, and points $p_i,q_i \in  \E(t) \subset H_i^n$ with $d_{g_{Eucl.}}(p_i,q_i) =1$, so that:\\

The intrinsic distance $d_{\E({1/i})}(p_i,q_i)$ in $\I({1/i})^c$, that is, the infimum of the lengths of all connecting arcs in $\I({1/i})^c$, determined using $g_{H_i}$,  diverges:
\[d_{\I({1/i})^c }(p_i,q_i) \ge i, \mm{ or it even equals }  +\infty,\] when $\E({1/i})$ is not connected and the points are in different components.\\

We may assume  $p_i= 0= (0,..,0), q_i = e_1 =(1,0,..,0)$. Also we may assume that the $H_i$ converge compactly to a limit area minimizer $H_\infty$, cf. [Gi], 1.19 and 9.1.\\

The Lipschitz estimate $|\delta_{\bp}(p)- \delta_{\bp}(q)|   \le L \cdot d(p,q)$, shows that \[B_r(p_i), B_r(q_i) \subset \E(t/2), \mm{ for any } r \in (0, t/(2 \cdot L)).\]

Now we claim that these balls converge to corresponding balls in the limit. This is not evident from the intrinsic geometry, since the balls $\subset H_i$ could be approached by other sheets of $H_i$
and then be annihilated in the limit. This is another occasion to apply the non-extinction result \ref{nex}. As a consequence we have:\\

 \emph{The limit points of the $p_i$ and $q_i$, $0, e_1$ belong to  $\E(t) \subset H_\infty$.}\\

Now, since,  $H_\infty \setminus \Sigma_{H_\infty}$ is again rectifiably connected and $\bp$ is a proper function on $\overline{B_R(0)} \cap H_\infty \setminus
\Sigma_{H_\infty}$, for any given $R >0$, cf.\ref{locsking}(i), there is a smooth path $\gamma$ in $H_\infty$, that connects $0, e_1$ within $\subset \E(\tau)$ for some suitably small $\tau \in (0, t)$. \\

Again the properness of $\bp$ on $\overline{B_R(0)} \cap H_\infty \setminus
\Sigma_{H_\infty}$ shows that  we can find a smooth tube $V$ around $\gamma$ with $V \subset \E(\tau/2)$ and thus we observe a smooth
convergence of suitably chosen sets $V_i \subset H_i$ to $V$ and rectifiable paths $\gamma_i \subset V_i$, connecting $p_i$ and $q_i$, to the path $\gamma \subset H_\infty$. We infer for large $i$:
$$l(\gamma_i) \le l(\gamma) +1$$ contradicting the assumption.  \qed

\noindent \textbf{Step 3 (Pipelines of Short Quasi-Geodesics)}  \quad \label{netqua}\\

 Now we assemble short quasi-geodesics to pipelines in sufficiently small but uniformly sized balls in $H$, where we deliberately allow center points in $\Sigma_H$.

\begin{lemma} \label{coprq} \quad Let $H^n  \subset \R^{n+1}$ be a complete area minimizing hypersurface which is an oriented boundary of an open set $A  \subset \R^{n+1}$.
Then, for any $\varpi \in (0,1)$ there is some $t_{n,\varpi} \in (0,1)$, so that for any $t \in (0,t_{n,\varpi})$, $k \in \Z$ and any two $x,y \in H$ with $d_{\R^{n+1}}(x,y)=3/2$
{\small \begin{equation}\label{w} \frac{Vol_n(\E(2^{-k} \cdot t) \cap B_{2^{-k+1}} \setminus B_{2^{-k}}(x) \cap H)}{Vol_n(B_{2^{-k+1}}  \setminus B_{2^{-k}}(x) \cap H)} \ge \varpi,
\end{equation}
  \begin{equation}\label{w1}
  \frac{Vol_n(\E(t/2) \cap B_{1}  \setminus B_{1/2}(x) \cap  B_{1}\setminus B_{1/2}(y) \cap H)}{Vol_n(B_{1} \setminus B_{1/2}(x) \cap  B_{1} \setminus B_{1/2}(y) \cap H)} \ge \varpi,
\end{equation}}
where $B_a\setminus B_b(z):=B_a(z)\setminus B_b(z)$ denote, differences of, distance spheres in $\R^{n+1}$.
\end{lemma}

In particular, for any $\ve >0$ and $\varpi(\ve)$ sufficiently close to $1$, there is a point $z_p \in \E(2^{-k} \cdot t) \cap B_{\ve}(p)$ for any given point
$p \in B_{2^{-k+1}} \setminus B_{2^{-k}}(x)$.\\

\textbf{Proof} \quad The isoperimetric inequality cf.[Gi],5.13-5.14, Ineq.(5.16) and a simple comparison with the (larger) volume of $\p B_1(0)$ give
positive constants $c_n^\pm$ depending
only on the dimension $n$, so that for any $x \in H$: \,$c_n^- \le Vol_n(B_{1/2}(x) \cap H) \le c_n^+.$\\

For any given $H$ and $x \in H$ we know that
{\small \[\frac{Vol_n(\E(t/2) \cap B_{1}  \setminus B_{1/2}(x) \cap H)}{Vol_n(B_{1} \setminus B_{1/2}(x) \cap H)} \ra 1,\mm{ for } t \ra 0,\]}
since $\bigcap_{t >0}\I(t/2) = \Sigma$ and $Vol_n(\Sigma) = 0$.\\

Therefore, compactness arguments for area minimizers, completely similarly to those in \ref{copr}, show that there is a $t_{n,\varpi} \in (0,1)$,  so that for any $t \in (0, t_{n,\varpi}]$ and for any point $x \in H$:
{\small \begin{equation} \label{v}Vol_n(\E(t/2) \cap B_{1}  \setminus B_{1/2}(x) \cap H) \ge \varpi \cdot Vol_n(B_{1} \setminus B_{1/2}(x) \cap H) \ge \varpi \cdot c_n^- .\end{equation}}
This implies the assertion (\ref{w}) for $k=0$, the case for $k \neq 0$  follows from scaling by $2^k$ and the way $\bp$ and the volumes transform under these scalings. The second claim, that is  (\ref{w1}), can be derived in the same manner. \qed

Now we return to the general case of $H \in {\cal{H}}^n$, that is,  $H$ could also be a compact area minimizer. Henceforth, we assume, similarly to earlier cases above, that $M$ had been scaled by some large constant,
so that for some $\varpi$ very close to $1$, every ball of radius $r \le 10 \cdot L(t_{n,\varpi},n) + 5$ in $M$ is very close to the ball $B_r(0) \subset \R^n$, in $C^5$-topology.\\

 In particular we may apply \ref{copr} and \ref{coprq} within these balls. \\

\textbf{Assembly of Pipelines}\quad Now we want to assemble quasi-geodesics which join any given point $x$ with points in $B_1 \setminus B_{1/2}(x)$ with controlled
$\bp$-distance.\\

From \ref{coprq},(\ref{w}) we can choose some $p_k \in \E(2^{-k} \cdot t_{n,\varpi}) \cap B_{2^{-k+1}} \setminus B_{2^{-k}}(x) \cap H$, for $k \ge 0$ and may assume that the
extrinsic distance $d_{g_M}(p_k,p_{k+1}) = 2^{-k}$. From \ref{copr} we get
\begin{itemize}
    \item  For some $\tau(t_{n,\varpi},n) < t_{n,\varpi}$,  $p_k$ and $p_{k+1}$ can be connected by a path $\gamma_{p_k,p_{k+1}}$ with:
\[l(\gamma_{p_k,p_{k+1}}) \le L/2^k \mm{ and } \gamma_{p_k,p_{k+1}} \subset \E(2^{-k} \cdot \tau)\subset H\]
    \item Since $L$ may be much larger than $1$ we usually have \[\gamma_{p_k,p_{k+1}} \varsubsetneq \E(2^{-k} \cdot \tau) \cap B_{2^{-k+1}} \setminus B_{2^{-k}}(x) \cap H.\]
However, since $l(\gamma_{p_k,p_{k+1}}) \le L/2^k$, we get, when we write $\tau^*:=  \tau \cdot 2^{-m}$ and
$m:= \mm{\emph{the smallest integer}}\ge \log_2(10 \cdot L)$:
\begin{equation}\label{vk01}  \gamma_{p_k,p_{k+1}} \subset  \E(2^{-(k-m)} \cdot \tau^*) \cap B_{2^{-(k-m)+1}}(x)  \cap H. \end{equation}
\end{itemize}
Now we glue the paths $\gamma_{p_k,p_{k+1}}$, for all $k \ge 0$, that is, we identify the endpoint of $\gamma_{p_k,p_{k+1}}$
with the start point of $\gamma_{p_{k-1},p_{k}}$ and get a path $\Gamma=\Gamma_{p_0,x}$ from $p_0$ to $x$.\\

Each point $z$ on the subarc $\gamma_{p_k,p_{k+1}} \subset \Gamma$  remains within a distance $\le L/2^k$ to
both endpoints $p_k$  and $p_{k+1}$. This way we get estimates for any point $z \in \gamma_{p_k,p_{k+1}} \subset\Gamma$ and for the subarc $\Gamma(z) \subset \Gamma$ from $x$ to $z$:
\begin{equation}\label{vk1}
d_{g_H}(x,z) \le l(\Gamma(z)) \le \sum_{a \ge k} L/2^a = L/2^{k-1} \le  L \cdot d_{g_M}(x,z) \le  L \cdot d_{g_H}(x,z) \end{equation}
\[\mm{ and } \, l(\Gamma(z)) \le 2 \cdot  L/ \tau \cdot \delta_{\bp}(z).\]

For the latter inequality we use  (\ref{vk01}). It shows that $\Gamma(z) \in \E(2^{-k} \cdot \tau)$, that is, $\delta_{\bp}(z) \ge 2^{-k} \cdot \tau \ge \tau/(2 \cdot L) \cdot l(\Gamma(z))$. \qed

\begin{remark} \label{c} \quad Before we proceed with the concluding Step 4, we note some immediate consequences of this construction:  \ref{fini}(i) and \ref{finii}.\\

Any singular point can be linked to regular points in $H$ by a rectifiable path. Since $H \setminus \Sigma$ is rectifiably connected
we get that $H$ is \textbf{rectifiably connected}. Also we observe the one non-trivial inequality in \ref{finii}: \;
$d_{g_H}(p,q)  \le  c^{\R}_n \cdot  d_{g_{\R^{n+1}}}(p,q).$\qed
\end{remark}

\bigskip

\noindent \textbf{Step 4A (Skin Uniformity for Euclidean Hypersurfaces)}  \quad \label{veri}\\

Now we assemble the pipelines of Step 3 to derive the skin uniformity in the case of Euclidean minimal hypersurfaces in ${\cal{H}}^{\R}_n$.\\

 We consider $x,y \in H \setminus \Sigma$. The absolute value of the distance between $x$ and $y$ does not matter since we are in a scaling invariant situation
and also the skin uniformity conditions are scaling invariant.\\

Thus we may assume that $d = 3/2$. Then, we consider the unit balls $B_1(x) \cap H$ and $B_1(y) \cap H$. Then (\ref{w1})
shows that we may choose a common start point  $p_0 \in (B_1(x) \cap B_1(y)) \cap H$, for the pipelines  $\Gamma_{p_0(x),x}$ and $\Gamma_{p_0(y),y}$  described Step 3.\\

Now we use (\ref{vk1}) to see that the composition of these two paths defines a $c$-skin uniform curve linking $x,y$, for $c:=  4 \cdot L/ \tau + 4 \cdot L$.
This $c$, depends only on the dimension, when $\bp$ is given.
 This concludes the proof of \ref {intsuni}(iii).\qed

\bigskip

\noindent \textbf{Step 4B (Skin Uniformity for Compact Hypersurfaces)}  \quad \label{veri2}\\

Finally, for compact hypersurfaces in ${\cal{H}}^c_n$, we combine the pipeline construction with the coarse compactness argument of \ref{rnpc} to verify their skin uniformity.\\

We consider $x,y \in H \setminus \Sigma$ with distance $d := d_M(x,y) >0$. For we assume $M$ has already been scaled as described in the beginning of the previous step,
so that locally it looks virtually like the Euclidean space.\\

We distinguish between different cases depending on the proximity of the points to $\Sigma$ and to each other.
 We note that we can always scale by a constant $>1$ to understand the situation, since this only brings the geometry closer to that of Euclidean hypersurfaces,
 the case we take the estimates from.\\

\textbf{1.}  For $d \le 1$,  we may assume that we have chosen a local decomposition of $H$ into oriented minimal boundaries and that $x,y \in H \setminus \Sigma$ belong to the same boundary.
Now we  argue similarly as in the previous step.
We scale $M$ by $\frac{3}{2\cdot d} \ge 1$ so that we may assume that $d = 3/2$. Then, we consider the unit balls $B_1(x) \cap H$ and $B_1(y) \cap H$. Then (\ref{w1})
shows that we may choose a common start point  $p_0 \in (B_1(x) \cap B_1(y)) \cap H$, for the pipelines  $\Gamma_{p_0(x),x}$ and $\Gamma_{p_0(y),y}$  described Step 3.\\

Now we use (\ref{vk1}) to see that the composition of these two paths defines a $c$-skin uniform curve linking $x,y$, for $c:=  4 \cdot L/ \tau + 4 \cdot L$.\\

\textbf{2.} For $d > 1$, we choose the compact set $K = H \setminus U_{1/4}(\Sigma)$. When $x$ or $y \in U_{1/4}(\Sigma)$, we first choose closest points $x_K, y_K \in K$ and
this time we choose two local decompositions of $H$ into oriented minimal boundaries and that the pairs $x,x_K$ and $y,y_K$, but not necessarily $x$ and $y$,  belong to the same boundary.\\

Now we define the pipelines $\Gamma_x$ and $\Gamma_y$ with endpoints in $x$ and $y$ as in Step 3.
Of course, when  $x$ or $y \notin U_{1/4}(\Sigma)$, this is merely a constant path.\\

Next, for these $x_K,y_K\in K$,   \ref{rnpc} of Step 1 gives us some $l=l(K) >0$, $s = s(K) > 0$ so that they can be linked
through a rectifiable curve $\gamma_{x_K,y_K}$ in $\E(s)$ of length $\le l$. Since $d(x_K,y_K)>1/2$, this can be rewritten as
\begin{equation}\label{pre}
l(\gamma_{x_K,y_K}) \le l = l/d \cdot d \le 2 \cdot l \cdot d_H(x_K,y_K)\end{equation}
\[ \mm{ and }\,\,  l_{min}(\gamma_{{x_K,y_K}}(z)) \le l \le l/s \cdot s \le l/s \cdot \delta_{\bp}(z).\]

Now we append the pipelines $\Gamma_x$ and $\Gamma_y$ to $\gamma_{{x_K,y_K}}$. This defines an arc $\Gamma_{x,y}$ that links $x$ and $y$. We check its skin uniformity properties.\\

The curves $\Gamma_{x,y}$ are quasi-geodesic as is easily seen using (\ref{vk1}), (\ref{pre}) and the inequalities $d(x,x_K),d(y,y_K) \le 1/4 \le d(x,y)/4$:
{\small \[l(\Gamma_{x,y}) \le 2 \cdot L/2 + 2 \cdot l \cdot d_H(x_K,y_K) \le L \cdot d_H(x,y) + 2 \cdot l \cdot (d_H(x,y) + d_H(x,x_K)+d_H(y,y_K))\]
\[\le  (L  + 4 \cdot l) \cdot d_H(x,y)\]}
Now, for the doubled twisted cone condition, we start with an intermediate point $z$ on $\Gamma_x$, and observe that, along $\Gamma_{x,y}$, this point is closer to $x$ than to $y$.
Hence we have from (\ref{vk1})
 \[l_{min}(\Gamma_{x,y}(z)) = l(\Gamma_x(z)) \le 2 \cdot  L/ \tau \cdot \delta_{\bp}(z),\]
and, similarly, we get the counterpart for $z$ on $\Gamma_y$. For $z$ on $\gamma_{x_K,y_K}$, we first consider the subcase $z \in B_{1/8}(x_K) \cap \gamma_{{x_K,y_K}}$.
The Lipschitz continuity of $\delta_{\bp}$, with assumed Lipschitz constant $1$, then shows $\delta_{\bp}(x_K)/2 \le \delta_{\bp}(z)$ and hence we get
 \[l_{min}(\Gamma_{x,y}(z)) = l_{min}(\gamma_{{x_K,y_K}}(z)) + l(\Gamma_x) \le \]
 \[ l/s \cdot \delta_{\bp}(z) + 2 \cdot  L/ \tau \cdot \delta_{\bp}(x_k) \le (l/s + 4 \cdot  L/ \tau) \cdot \delta_{\bp}(z).\]
Again we similarly get its counterpart for $z \in  B_{1/8}(y_K) \cap \gamma_{{x_K,y_K}}$.  To understand the remaining subcase, where
$z \in \big(H \setminus  (B_{1/8}(x_K) \cup B_{1/8}(y_K))\big) \cap \gamma_{{x_K,y_K}}$, we first note from (\ref{pre}) that:
\[1/8 \le  l_{min}(\gamma_{{x_K,y_K}}(z)) \le  l/s \cdot \delta_{\bp}(z).\]
Therefore, we get for these $z$:
 \[l(\Gamma_x) \le  L \cdot d_{g_H}(x,x_K) \le L/4 \le 2 \cdot l/s \cdot L \cdot  \delta_{\bp}(z)\]
and consequently
 \[l_{min}(\Gamma_{x,y}(z)) = l_{min}(\gamma_{{x_K,y_K}}(z)) + l(\Gamma_x) \le  (1+2 \cdot L) \cdot l/s \cdot \delta_{\bp}(z).\]
Summarizing all cases, we get that $H \setminus \Sigma$ is $c$-skin uniform for \[c:= 4 \cdot   (l + l/s + L + L/ \tau + L \cdot l/s).\]
 \qed

\subsubsection{Skin Uniform Domains}\label{ch}
\bigskip

Here we give another characterization of the skin uniformity of $H \setminus \Sigma_H$ in terms substructures, the \emph{skin uniform domains} $\su \subset H \setminus \Sigma_H$. Their uniformity property takes care of both, the metric and the $\bp$-distance. \\

To exclude trivialities, we henceforth assume that our hypersurfaces are not totally geodesic. The main goal of this section the proof of the following existence result. It is again split  in a number of subresults.

\begin{proposition}\label{std} \emph{\textbf{(Skin Uniform Domains)}} \, Let $H \in {\cal{H}}$ be a $c$-skin uniform space for some $c(H) >0$.
Then there are some $\iota(c),\kappa(c) > 1$ so that for sufficiency small $a >0$, there is a domain $\su= \su(a) \subset H \setminus \Sigma$:
\begin{itemize}
  \item    \quad  $\E(\iota \cdot a) \subset \su(a) \subset \E(a)$, where $\E(a):= \{x \in H \setminus \Sigma \,|\, \delta_{\bp}(x) \ge a \},$
\end{itemize}
and for any two points $p,q \in \su$ can be linked by an arc $\gamma_{p,q} \subset \su$ with
\begin{itemize}
    \item \quad $l(\gamma)  \le \kappa \cdot d_H(p,q)$,
    \item \quad $l_{min}(\gamma(z)) \le \kappa \cdot \min\{L_{\bp} \cdot dist(z,\p \su),\delta_{\bp}(z)\}$, for any $z \in \gamma_{p,q}$.
\end{itemize}
For any converging sequence $H_i \ra H_\infty$ in ${\cal{H}}$, skin uniform domains $\su_i(a) \subset H_i$ subconverge to a skin uniform domain $\su_\infty(a) \subset H_i$.
\end{proposition}

\begin{remark}  \label{dspr} \, 1. For a compact $H$, every smoothly bounded domain $D \subset H \setminus \Sigma$ satisfies the skin uniformity conditions for some sufficiently large $\iota_D, \kappa_D >0$.
However, these $\iota_D$ and $\kappa_D$ depend on $D$. Thus, the critical and  non-trivial detail is that the $\su(a)$ are $\kappa$-skin uniform, for some $\kappa >0$ that is \emph{independent} of $a$. \\

2. In turn, for a non-compact $H$, even the definition of one non-compact domain $D \subset H \setminus \Sigma$ satisfying the conditions of \ref{std} becomes non-trivial and uses much of the skin uniformity of $H \setminus \Sigma_H$.\\

3. Indeed, the existence of such a family of domains $\su(a)$, for every small $a>0$, even implies the skin uniformity of $H \setminus \Sigma_H$, as is readily seen from $\delta_{\bp}(x) \le L_{\bp}  \cdot dist(x,\Sigma)$, when we let $a \ra 0$. In other words, the existence of these domains also characterizes skin uniform spaces. \qed
\end{remark}

\textbf{Idea of the Proof} \, To define the domains $\su$ we start from the observation that the domains $\E(a)= \{x \in H \setminus \Sigma \,|\, \delta_{\bp}(x) \ge a \}$ are already nearly skin uniform for pairs of points sufficiently far away from $\p \E(a)$. This is a quite immediate consequence of the skin uniformity of $H \setminus \Sigma$.\\

 But, for a general skin transform $\bp$, we cannot expect the validity of the uniformity conditions for pairs of points while we approach the rather uncontrolled border $\p \E(a)$. (An exception are skin transforms with \sks\ minimizing a generalized area functionals as discussed in Ch.\ref{localskin}. We know from \ref{ds} that for them we automatically get some uniformity control. But we do \emph{not} base our argument on this non-trivial result, since we want to derive a proof valid for any skin transform, from the axioms only.)\\

 Thus we need to modify $\p \E(a)$. The skin uniformity suggests to consider linking paths also outside $\E(a)$. We adjoin these paths to $\E(a)$ and show that a particularly adapted neighborhood of this extension is a domain with the asserted properties.\\

We define these neighborhoods from a discretization: we cover the path extension by balls we take from skin adapted covers, as already indicated in Ch.\ref{skinada1}.  To ensure the uniformity conditions, we now also need some quantitatively controlled transversal intersections, we use the QT-skin adapted covers we define below.\\

Remarkable schemes to treat such transversality problems had been developed by Yomdin, using tools from real algebraic geometry cf.[YC]. However, for our problems we give a self-contained construction without appealing to those methods.\qed\\

\textbf{Quantitative Transversality}  \, We add a further control to skin adapted covers which says that the boundaries of any two  balls in ${\cal{A}}$ intersect transversally with a lower bound for the degree of transversality. More generally, any ball
intersects the intersection of other balls with such a controlled transversality.

\begin{proposition}\label{besi2} \emph{\textbf{(Quantitative Transversality)}}  \, There is a transversality constant $\tau_n \in (0,1)$, depending only on the dimension $n$,
so that a cover ${\cal{A}}$ of $H \setminus \Sigma$, can be chosen so that, in addition to the properties in \ref{besi} we have the following intersection estimates:
\begin{description}
\item[(QT1)] \quad $B_{2 \cdot \Theta(p)}(p) \cap B_{2 \cdot\Theta(q)}(q) \n \, \Rightarrow \, B_{(2-\tau_n) \cdot \Theta(p)}(p) \cap B_{(2-\tau_n) \cdot\Theta(q)}(q) \n$,
\item[(QT2)] \quad $B_{2 \cdot \Theta(p)}(p) \cap B_{2 \cdot\Theta(q)}(q) \v \, \Rightarrow \, B_{(2+\tau_n) \cdot \Theta(p)}(p) \cap B_{(2+\tau_n) \cdot\Theta(q)}(q) \v$,
\end{description} for any two $p,q \in A_{Q}$. We call such a cover, a \textbf{QT-skin adapted cover}.
\end{proposition}

For the proof, we start from an arbitrary skin adapted cover ${\cal{A}}$, for some very small $\xi >0$, as described in \ref{besi}, and show that there is a perturbation
scheme to readjust ${\cal{A}}$ so that, afterwards, it also satisfies the transversality constraints (QT1) and (QT2).\\

 {\bf Proof} \quad We start with the cover $\cal{A}$ of \ref{besi}. Inductively, we slightly move all balls in $A^Q(i)$, $1 \le i \le c(n)$, family by family. We do this
 by uniformly localized arguments which apply simultaneously to the infinitely many balls in each $A^Q(i)$.\\

We note from the argument in \ref{besi}, that for any tiny $\delta >0$, we can choose $\xi$ so small that all balls $B_{10 \cdot \Theta(p)}(p)$, after scaling by $1/(10 \cdot
\Theta(p))$,
 $p \in Q$ are $\delta$-close to the Euclidean unit ball in $C^3$ via exponential maps. Also due to the Lipschitz continuity of $\delta_{\bp}$,
we can assume that $(1-\delta) \cdot \Theta(p) \le \Theta  \le (1+\delta) \cdot \Theta(p)$ on $B_{10 \cdot \Theta(p)}(p)$.\\

This shows that upon choosing $\delta >0$ small enough it suffices to solve the problem for any given cover $\cal{A}$ of the Euclidean space by balls of radius $1$ sorted in
$c(n)$ families $A(1),,,A(c(n))$, so that \[ B_{10}(p) \cap B_{10}(q) \v,\mm{ for } p, q \in A(k) \mm{ and } q \notin \overline{B_1(p)}\mm{, for any two } p, q \in A.\] Our aim is to
move the points in $A$ by at most $1/10$, to get a new set $A^*$, so that,  for any two $p, q \in A^*$ and for any $\tau \in (0,\tau_n]$,  for some suitable $\tau_n \in (0,1/10)$:
\begin{description}
    \item[\textbf{(QT1)}] \quad $B_2(p) \cap B_2(q) \n \, \Rightarrow \, B_{(2-\tau)}(p) \cap B_{(2-\tau)}(q) \n$
    \item[\textbf{(QT2)}] \quad $B_2(p) \cap B_2(q) \v \, \Rightarrow \, B_{(2+\tau)}(p) \cap B_{(2+\tau)}(q) \v$
\end{description}
 We work inductively. We leave the balls in $A(1)$ unchanged, that is, $A^*(1):=A(1)$.\\

 Next, we assume we have constructed all families $A^*(i)$, up to $i \le m$ for some constant $\tau_n[m] \in (0,1/10)$, so that for any pair of points
 $p,q \in \bigcup_{j \le i}A^*(i)$, (QT1) and (QT2) are satisfied for any $\tau \in (0,\tau_n[m]]$,  and so that the distance that any of the points in $\bigcup_{j \le i}A(i)$ has been moved was $\le 1/(100 \cdot c(n))$.\\

Now consider $B_2(p)$, for some $p \in A(m+1)$. We reformulate the condition (QT1) and (QT2): for a pair of points that does \emph{not} satisfy these
conditions, this means: {\small \begin{description}
    \item[(S1)] \quad  $B_2(p) \cap B_2(q) \n$, but $B_{(2-\tau)}(p) \cap B_{(2-\tau)}(q) \v \, \Leftrightarrow \, 2 \cdot (2-\tau) \le  d(p,q) \le 2 \cdot 2$,
    \item[(S2)]  \quad   $B_2(p) \cap B_2(q) \v $, but $B_{(2+\tau)}(p) \cap B_{(2+\tau)}(q) \n \, \Leftrightarrow \, 2  \cdot 2 \le  d(p,q) \le 2 \cdot (2+\tau)$.
\end{description}}

In other words, a failure of (QT1) and (QT2) can be expressed through \emph{shell conditions} (S1) and (S2) for shells around the $q$. We avoid these shells by appropriate moves of $p$.\\

We know from the definition (\ref{defc}) of $c(n)$ that, also after the perturbation of $\bigcup_{j \le i}A(i)$ into $\bigcup_{j \le i}A^*(i)$, (S1) or (S2) are satisfied for at
most $m \le c(n)$ balls $B_2(q)$ with center in $\bigcup_{1 \le i \le m}A(i)^*$, and some given $\tau \in (0,\tau_n[m])$.\\

 For $\tau \in (0,\min\{\tau_n[m],1/(10^3 \cdot c(n)\})$ chosen small enough, the shells of (S1) and (S2) locally look like one-sided $\tau$-distance collars of hyperplanes.
 The two collars on both sides of the hyperplanes combine to the distance tube.\\

For such a $\tau >0$, we consider a ball $B_{\ve}(z) \subset B_2(p)$, for some small $\ve > 0$.
\begin{itemize}
    \item $B_{\ve}(z)$ has volume $Vol_n(B_{\ve}(z))= v_n \cdot \ve^n$, for some constant $v_n > 0$.
    \item For any hyperplane $N$, we have $Vol_{n-1}(N \cap B_{\ve}(z))\le v_{n-1} \cdot \ve^{n-1}$.
    \item  For the intersection the $\tau$-tube $U_\tau(N)$ of $N$ we have \[Vol_n(U_\tau(N) \cap B_{\ve}(z)) \le 2 \cdot \tau \cdot v_{n-1} \cdot \ve^{n-1}\]
\end{itemize}
Thus, for the intersection of $B_{\ve}(z)$ with the union of the relevant hypersurface $\tau$-tubes the volume is $\le 2 \cdot m \cdot \tau \cdot v_{n-1} \cdot \ve^{n-1}$. Therefore, when we take
some $\tau <  \ve \cdot {v_n}/(2 \cdot m \cdot v_{n-1})$, the $\tau$-tubes do \emph{not} cover the entire ball $B_{\ve}(z)$. In other words, we only need to move $p$ for a distance  $\le \ve$ to shift the point away from the shells for any of the
 at most $m \le c(n)$ relevant balls.\\

We choose such an $\ve < 1/(10^4 \cdot c(n))$ and set \[\tau_n[m+1] := \min\{\tau_n[m],1/(10^3 \cdot c(n)),\ve \cdot {v_n}/(2 \cdot m \cdot v_{n-1}\}.\]

 Thus, after such a move for any $p \in A(m+1)$, conditions (QT1) and (QT2) and the induction hypotheses for the next loop are satisfied for this $\tau_n[m+1]$.\\

Since we only need to repeat this process for $c(n)$ times, the resulting transversality constants $\tau_n:=\tau_n[c(n)]$ remains positive and we observe that it only depends on the
dimension $n$. \qed

An interpretation of this QT-enhancement is that any connected union of a subcollection of balls has \emph{quantitative connectivity} properties.\\

To formulate this in detail, we choose again a small size parameter $\xi$. Then we can assume that $(1-\delta) \cdot \Theta(p) \le
\Theta \le (1+\delta) \cdot \Theta(p)$ on $B_{10 \cdot \Theta(p)}(p)$ for some tiny $\delta >0$ and, as in \ref{besi2}, $ \Theta(p)^{-1} \cdot B_{10 \cdot \Theta(p)}(p)$ is
arbitrarily close to the Euclidean ball $B_{10}(0)$ in $C^3$-norm and we can focus on the problem in the Euclidean model case.

\begin{corollary} \label{tra} \quad For any two balls in $\cal{A}$ with $B_{2 \cdot \Theta(p)}(p) \cap B_{2 \cdot\Theta(q)}(q) \n$ there is a radius $\omega_n
>0$ so that, for the geodesic arc $\gamma_{p,q}$ of length $l(\gamma_{p,q}) \le 4 \cdot \Theta(p)$ linking $p$ and $q$, we have
 \begin{equation}\label{tu}
U_{\omega_n \cdot \Theta(p)}(\gamma_{p,q}) \subset B_{2\cdot \Theta(p)}(p) \cup B_{2 \cdot\Theta(q)}(q),
\end{equation} where $U_{\omega_n \cdot \Theta(p)}(\gamma_{p,q})$ denotes the $\omega_n \cdot \Theta(p)$-distance tube around $\gamma_{p,q}$.
\end{corollary}

\textbf{Proof} \quad This follows from  (QT1), i.e. $B_{(2-\tau_n) \cdot \Theta(p)}(p) \cap B_{(2-\tau_n) \cdot\Theta(q)}(q) \n$, since this translates into a lower estimate for the sinus of the largest opening angle
for cones with tip $p$ disjoint to $\p B_{2\cdot \Theta(p)}(p) \setminus B_{2 \cdot\Theta(q)}(q)$, that is, for cones passing through the opening of this spherical complement. \qed

Now we turn to a microscopic variant of (\ref{tu}) to get locally uniform boundary controls for the union $V$ of
subcollections of balls  in some QT-skin adapted cover.\\

 Since we only derive local estimates from given configurations it is enough
to consider some \emph{connected} $V$ which is a union of a finite collection $B_{2 \cdot \Theta(p_1)}(p_1)$,...,
$B_{2 \cdot \Theta(p_m)}(p_m)$, $m \le c(n)$, of such balls
and, as in \ref{besi2}, due to the Lipschitz continuity of $\delta_{\bp}$ we may choose $\xi$  so small, that we can assume that $\Theta_V:=\Theta(p_1)=...=\Theta(p_m)$.
For some $p \in \p V_{\Theta_V}$ and a small $\rho  \in (0,\omega_n/10^5)$ we observe, again as in \ref{besi2}, up to some arbitrarily small distortion:
\begin{enumerate}
    \item $(\rho \cdot \Theta_V)^{-1} \cdot (B_{\rho \cdot \Theta_V}(p) \cap V)$ is the intersection of the unit ball in $\R^n$ with the union of (some of) the halfspaces
        $H_k \subset \R^n$, each corresponding to some balls $B_{2 \cdot \Theta_V}(p_k) \subset V$.  At least one of the $\p H_k$, say $\p H_1$, passes through $p$.
    \item Due to (QT1) there is a uniform lower bound for the angles $\measuredangle(\p H_k,\p H_j)$ along intersections of any two balls
    $B_{2 \cdot \Theta_V}(p_k), B_{2 \cdot \Theta_V}(p_j)$. There are some $\beta^*_n \ge \beta_n >0$, with $\beta^*_n < \pi$, so that
     \begin{equation}\label{ang} \beta^*_n \ge\measuredangle(\p H_k,\p H_j) \ge \beta_n >0.
\end{equation}
\item From (QT2)  we infer that any two $H_j \neq H_k$, both with non-trivial intersection with the unit ball $(\rho \cdot \Theta_V)^{-1} \cdot B_{\rho \cdot \Theta_V}(p)$ intersect.
(ii) shows that this happens within a ball of radius $ct_n:=cotangens(\beta_n/2)$ centered in $p$.
\end{enumerate}
We use (iii) to see that points $x,y \in V$ sufficiently close to the boundary \emph{and} to each other can be linked by some effective path $\subset V$. In the terminology of
uniform spaces this is a $c_n$-\emph{uniform arc}.

\begin{lemma} \emph{\textbf{(Micro-Uniformity)}}\label{pu}  \, For some  $\rho_n \in (0,1)$ and some $c_n >0$, both depending only
on the dimension $n$, we have for any two points $x,y \in B_{\rho_n \cdot \Theta_V}(p) \cap V$, for some $p \in \p V$
\begin{equation}\label{ain}
 l(\gamma)  \le c_n \cdot  d_{\R^n}(x,y)\, \mm{ and }\, l_{min}(\gamma_{x,y}(z)) \le c_n \cdot dist(z,\p V),  \mm{ for any } z \in \gamma_{x,y}.
\end{equation}
where $l_{min}(\gamma_{x,y}(z))$ denotes the minimum of the lengths of the subarcs of $\gamma_{x,y}$ from $x$ to $z$ and from $y$ to $z$.
\end{lemma}

 {\bf Proof} \quad  For $d:= d_{\R^n}(x,y) < 2$, we choose $\rho_n \in (0,\omega_n/10^5)$ not only small enough to ensure that observation (i)-(iii) apply, but that $10 \cdot ct_n \cdot d \ll \omega_n$.
 That is, we act only in regions of a size much smaller than that of the openings of \ref{tra}.\\

We find, under the chosen halfspaces, two ones, labelled $H_x$ and $H_y$, with $x \in H_x$ and $y \in H_y$, so that they either equal or intersect within
 the ball of radius $ct_n \cdot d$ around $p$.\\

  In the first case, we can find a plane $P$ perpendicular to  $\p H_x= \p H_y$ so that $x,y \in P$. Now $P \cap H_x$ is a halfplane and
the Poincar\'{e} model of the hyperbolic plane shows that we can link $x$ and $y$ by a unique  circle segment that extends to a halfcircle that hits $\p (P \cap H_x)$ perpendicularly. These arcs satisfy
(\ref{ain}) for some $c_n(1) >0$, since we locally always have $dist(\cdot, \p V) \ge dist(\cdot, \p H_x )$. Again, this is valid up to arbitrarily small distortion, when $\xi$ has only
been chosen small enough.\\

For the second case, we use observation (iii) above to find a point $p(x,y) \in H_x \cap H_y$ so that
\[d(x,p(x,y)), d(y,p(x,y)) \le 2 \cdot ct_n \cdot d \,\, \mm{ and }\,\, dist(p(x,y), \p (H_x \cup H_y))) \ge ct_n \cdot d.\]
Then we repeat the construction of the first case two times: we choose a circle segment that links $x$ and $p(x,y)$ in $H_x$ and append another one from $p(x,y)$  and $y$ in $H_y$ and we verify the
inequalities (\ref{ain}) for $c_n(2):= 5 \cdot (1 +  ct_n) \cdot c_n(1) > c_n(1) > 0$. Thus, we can choose $c_n :=c_n(2)$ and notice that the uniformity estimates only improve when we take
 further half spaces into account.\qed\\

\textbf{Assembly of the} $\mathbf{\su}$ \quad  The domains $\su$ result from some dual Stein-Whitney regularization of the sets $\E(a)$ where we exploit the skin uniformity
as already indicated in \ref{skinada1}. For starters, we describe the spaces we use to assemble the $\su$.

\begin{definition} \label{arc} \quad  For any skin transform $\bp$ and any $c$-skin uniform area minimizing hypersurface $H \in {\cal{H}}$ and $a >0$,
we define the following spaces \begin{itemize}
    \item The \textbf{maximal link space} $\Lambda_{\max}(H,a)$ is the space of all $c$-skin uniform arcs in $H \setminus \Sigma$  linking points in $\E(a)$, \quad
        $\Lambda_{\max}(H,a) := $
\[\{\gamma \in H \setminus \Sigma \mm{ is a $c$-skin uniform arc and joins two points }p,q \in \E(a) \}.\]
    \item  A subset $\Lambda \subset \Lambda_{\max}$ is called a \textbf{link space} if still any two points $p,q \in \E(a)$ can be linked by some arc in $\Lambda$.
    \item The \textbf{arc hull} $arc_\Lambda$ is the union of all traces of arcs in a link space $\Lambda$,
\[arc_{\Lambda}:= arc_{\Lambda}(\I_H^c(a)):=\{tr_\gamma\,|\, \gamma \in \Lambda\} \subset H \setminus \Sigma.\]
 \item  For a QT-skin adapted cover $\cal{A}$ of $H \setminus \Sigma$, we consider the subcover ${\cal{A}}_\Lambda$ which contains those balls
 $B_{\Theta(p)}(p) \in {\cal{A}}_\Lambda$ with $B_{\Theta(p)}(p) \cap arc_\Lambda \n$.  We call the union of the
 balls in $B_{\Theta(p)}(p) \in {\cal{A}}_\Lambda$ of doubled radius
 \[\mathfrak{B}^H_{{\cal{A}},\Lambda,a}:= \bigcup_{B_{\Theta(p)}(p) \in {\cal{A}}_\Lambda}B_{2 \cdot \Theta(p)}(p) \subset H \setminus \Sigma, \mm{ a
 \textbf{bubbled arc hull}, briefly a $\mathfrak{B}$-hull.}\]
\end{itemize}
\end{definition}

The arc hulls $arc_\Lambda$ can be thought as weak forms of convex hulls which integrate the skin uniformity of $H$ into the definition. With the
$\mathfrak{B}^H_{{\cal{A}},\Lambda}$ we additionally regularize the boundary. The $\mathfrak{B}^H_{{\cal{A}},\Lambda}$  are our candidates for the desired domains $\su$.
We start with some basic properties of these hull concepts.

\begin{proposition}\label{cch} \emph{\textbf{(Arc Hulls)}}  \, Let $\bp$ be a skin transform with Lipschitz constant $L_{\bp}$, $H \in {\cal{H}}$  $c$-skin uniform and $\alpha:= 1/(L \cdot c + 1)<1$.
Then, for any link space $\Lambda$ of $\E(a)$ and $a >0$,  we have
\begin{enumerate}
    \item $\E(a) \subset arc_{\Lambda}(\E(a)) \subset \E(\alpha \cdot a)$.
    \item  $\mathfrak{B}_{{\cal{A}},\Lambda,a}$ is rectifiably connected and for $\xi \in (0,1/(10^3 \cdot L_{\bp}))$ we have \[U_{\xi \cdot \alpha \cdot a/4}\left(arc_{\Lambda}(\E(a)) \right) \subset
        \mathfrak{B}_{{\cal{A}},\Lambda,a} \subset \E(\alpha \cdot a/4),\] where $U_{d}\left(arc_{\Lambda}(\E(a)) \right)$ denotes the $d$-distance neighborhood of
        $arc_\Lambda$.
    \item For any converging sequence $H_i \ra H_\infty$, $\mathfrak{B}^{H_i}_{{\cal{A}},\Lambda,a} \subset H_i$ compactly subconverges to some bubbled arc hull
    $\mathfrak{B}^{H_\infty}_{{\cal{A}},\Lambda_{\infty},a} \subset H_\infty$, for some link space $\Lambda_\infty \subset \Lambda_{\max}$.
\end{enumerate}
\end{proposition}

In general, we find $\Lambda_\infty \subsetneqq\Lambda_{\max}$, thus it does not make sense to focus on $\Lambda_{\max}$ only.\\

{\bf Proof} \quad For (i), $\E(a) \subset arc_{\Lambda}(\E(a))$ is obvious, since we link \emph{any point} in $\E(a)$ with some other one.\\

 Now let $\gamma$ be a c-skin uniform arc in $(H \setminus \Sigma, g_H)$
that links two points $p_,p_2 \in \M_{1/a} = \p \I(a)$. We get, as long as $z \in \gamma$ is closer to $p_1$ than to $p_2$:
\[ |\delta_{\bp}(p_1)- \delta_{\bp}(z)|   \le L \cdot d(p_1,z)\le L \cdot l_{min}(\gamma(z)) \le L \cdot c \cdot \delta_{\bp}(z).\]
When $z$ is closer to $p_2$ we get the corresponding inequality, and thus $\delta_{\bp}(z)$ satisfies one of the following four inequalities
\[\delta_{\bp}(z) \ge \delta_{\bp}(p_i)= a\, \mm{ or } \, \delta_{\bp}(z) \ge (L \cdot c + 1 )^{-1} \cdot \delta_{\bp}(p_i)=: \alpha \cdot a, \mm{ for } i=1,2,\]
and, since $\alpha < 1$, we see that in all cases $\delta_{\bp}(z) \ge \alpha \cdot a$.\\

For (ii), we have $\Theta(p)=\xi / \bp(p)$, for some $\xi \in (0,1/(10^3 \cdot  L_{\bp}))$ and thus the Lipschitz inequality for $\delta_{\bp}$ shows for any $q \in B_{\Theta(p)}(p)$:
\[|\delta_{\bp}(q)-\delta_{\bp}(p)|\le 1/10^3 \cdot \delta_{\bp}(p),\] in different terms
\[(1-1/10^3) \cdot \delta_{\bp}(p) \le \delta_{\bp}(q) \le (1+1/10^3) \cdot \delta_{\bp}(p).\]

From this and the triangle inequality we observe that for any $z \in arc_{\Lambda}(\E(a))$, $B_{\Theta(z)/2}(z) \subset \mathfrak{B}_{{\cal{A}},\Lambda,a}$. Also, from (i),
$\delta_{\bp}(z) \ge \alpha \cdot a/2$. Hence, $B_{\xi \cdot \alpha \cdot a/4}(z) \subset \mathfrak{B}_{{\cal{A}},\Lambda,a}$.\\

In turn, for any  $p \in \mathfrak{B}^H_{{\cal{A}},\Lambda,a}$ there is some $z \in arc_{\Lambda}(\E(a))$ so that
\[|\delta_{\bp}(z)-\delta_{\bp}(p)|\le L \cdot 4 \cdot \xi \cdot a \cdot \alpha < a \cdot \alpha/10.\]
and we get, either $\delta_{\bp}(p) \ge \delta_{\bp}(z) > a \cdot \alpha/2$, or $\delta_{\bp}(p) =  \delta_{\bp}(z) + \delta_{\bp}(p) - \delta_{\bp}(z)
= \delta_{\bp}(z) - |\delta_{\bp}(z)-\delta_{\bp}(p)| > a \cdot \alpha/4$.\\

For (iii), we argue similarly as in \ref{insu}. For any two points  $p,q \in \E(a) \subset H_\infty$ we find a sequence of $c$-skin uniform arcs in $k \cdot H$ which converge
some $c$-skin uniform arc in $H_\infty$ linking the two points. The set of arcs $\Lambda_\infty$ we get from this limit process is again a link space for $\I_{H_\infty}^c(\alpha \cdot a)$.\qed \\

\textbf{Skin Uniformity of $\mathfrak{B}_{{\cal{A}},\Lambda,a}$} \quad   The way we defined the $\mathfrak{B}$-hulls allows us to separate the verification
of their uniformity properties in two different cases.\\

In \ref{uni2}, the micro-uniformity we derived for QT-skin adapted covers in \ref{pu} gives us the existence of skin uniform curves
for pairs of points close to each other and close to the boundary. The chosen uniform arcs resemble hyperbolic geodesics in the Poincar\'{e} metric on the upper half plane.\\

In \ref{uni}, we treat the complementary case of pairs of points with bounded distance from the boundary or each other. The idea is to combine this
lower bound and with upper bounds for controlled paths from the given pair of points to pairs in $\I(a)$ which are hot-wired by skin uniform curves
in $\Lambda$.

\begin{lemma}\label{uni2}  \quad  For any $c$-skin uniform  $H \in \cal{H}$, $\xi \in (0, \min\{\alpha,1/(10^3 \cdot L_{\bp})\})$ and  any pair $p, q  \in \mathfrak{B}^H_{{\cal{A}},\Lambda,a}$ with
 \[\mm{ }\,
dist(p,\p \mathfrak{B}_{{\cal{A}},\Lambda,a}), dist(q,\p \mathfrak{B}_{{\cal{A}},\Lambda,a}) \le \rho_n \cdot \xi \cdot a/10 \mm{ and }
d_H(p,q) \le \rho_n \cdot \xi \cdot a/10,\]
can be linked by a curve $\gamma =
\gamma_{p,q} \subset \mathfrak{B}^H_{{\cal{A}},\Lambda,a}$, so that for any $z \in \gamma_{p,q}$:
\begin{equation} \label{eee} l(\gamma) \le c_n \cdot d_H(p,q) \,\mm{ and }\, l_{min}(\gamma_{p,q}(z)) \le c_n \cdot
\min\{dist(z,\p \mathfrak{B}_{{\cal{A}},\Lambda,a}),\delta_{\bp}(z)\},\end{equation}
where $c_n >0$ and $\rho_n \in (0,1)$ are the constants in \ref{pu}.
\end{lemma}

{\bf Proof} \quad Most of the assertions have already been proved in \ref{pu}. The additional factor $\xi <1$ only decreases the radius and, hence, the conclusions hold also for $\xi \cdot \rho_n$.\\

In particular we get from \ref{pu},(\ref{ain}), and \ref{cch}(ii), for any $z \in \gamma \subset \mathfrak{B}^H_{{\cal{A}},\Lambda,a}$:\\

$l(\gamma) \le c_n \cdot d_H(p,q) \le c_n \cdot \rho_n \cdot \xi \cdot a/10 \le c_n \cdot \alpha \cdot a/10$ and $\alpha \cdot a/4 \le \delta_{\bp}(z)$.
 \qed

\begin{lemma}\label{uni}  \quad  For any $c$-skin uniform  $H \in \cal{H}$ and $\xi >0$ small enough, there is some $d(c,L,\xi) >0$ so that any pair of points
$p, q  \in \mathfrak{B}^H_{{\cal{A}},\Lambda,a}$ satisfying \[\mm{ either } p, q  \in arc_\Lambda, \,\, \mm{ or }\,\, p, q  \in \mathfrak{B}^H_{{\cal{A}},\Lambda,a} \mm{ with }
d(p,q) \ge \rho_n \cdot \xi \cdot a/10,\] can be linked through an arc $\gamma = \gamma_{p,q}$, so that for any $z \in \gamma_{p,q}$: \[ l(\gamma) \le d \cdot d(p,q) \,\mm{ and }\,
l_{min}(\gamma_{p,q}(z)) \le d \cdot \min\{dist(z,\p \mathfrak{B}_{{\cal{A}},\Lambda,a}),\delta_{\bp}(z)\}.\]
\end{lemma}

{\bf Proof} \quad We follow the definition of $\mathfrak{B}^H_{{\cal{A}},\Lambda,a}$ and consider pairs of points first in $\E(a)$, then in $arc_{\Lambda}(\E(a))$ and
finally in $\mathfrak{B}_{{\cal{A}},\Lambda,a}$.\\

1. When $p,q \in \E(a)$, then the $c$-skin uniform arc $\gamma_{p,q}$ in $arc_{\Lambda}(\E(a))$ that joins $p,q$ and satisfies  for any $z \in \gamma_{p,q}$ \[
l(\gamma) \le c \cdot d(p,q) \,\mm{ and }\, l_{min}(\gamma_{p,q}(z)) \le c \cdot \delta_{\bp}(z),\] The remaining assertion is $ l_{min}(\gamma_{p,q}(z)) \le d \cdot  dist(z,\p
\mathfrak{B}_{{\cal{A}},\Lambda,a})$ for some suitable $d>0$.\\

For this, we combine the estimates  \[ l_{min}(\gamma_{p,q}(z)) \le c \cdot \delta_{\bp}(z), \,\,  \xi \cdot \alpha \cdot a/4 \le
dist(arc_{\Lambda}(\E(a)),\p \, \mathfrak{B}^H_{{\cal{A}},\Lambda,a})\]
 \[\mm{ and }\, |\delta_{\bp}(p)- \delta_{\bp}(z)|   \le L \cdot d(p,z):\] When we start from $p$ and
run along $\gamma_{p,q}$, the length $l_{min}(\gamma_{p,q}(z))$ gradually increases from $0$. Now we distinguish two subcases:\\

1a. We reach the midpoint of $\gamma_{p,q}$ before $l_{min}(\gamma_{p,q}(z))
\ge c \cdot a$. Then we have \[l_{min}(\gamma_{p,q}(z)) \le c \cdot a \le 4 \cdot c/( \xi \cdot \alpha) \cdot dist(arc_{\Lambda}(\E(a)),\p \, \mathfrak{B}^H_{{\cal{A}},\Lambda,a}),
\mm{ for any } z \in \gamma_{p,q}.\]

1b. There is a first point $z_0$ where $l_{min}(\gamma_{p,q}(z_0)) = c \cdot a = c \cdot \delta_{\bp}(z_0)$,  when we pass $z_0$. Now, the twisted double cone condition for $\gamma$
shows, that for the $z \in \gamma$ between $z_0$ and
the midpoint of $\gamma_{p,q}$, $\delta_{\bp}(z) >  \delta_{\bp}(z_0)$, that is, $z \notin \overline{\I_H(a)}$, and it also implies
\[l(\gamma_{p,q}[z_0,z]) \le c \cdot (\delta_{\bp}(z) - \delta_{\bp}(z_0)) \le c \cdot  L \cdot dist(z,\I_H(a))\]
Thus, we get:  {\small \begin{equation} \label{con2} l_{min}(\gamma_{p,q}(z)) \le c \cdot  L \cdot dist(z,\I_H(a)) +
4 \cdot c/( \xi \cdot \alpha) \cdot dist(arc_{\Lambda}(\E(a)),\p \mathfrak{B}^H_{{\cal{A}},\Lambda,a})\le \end{equation} \[c \cdot  L \cdot dist(z,\I_H(a)) +
4 \cdot c/( \xi \cdot \alpha) \cdot dist(\E(a),\p \mathfrak{B}^H_{{\cal{A}},\Lambda,a})  \le \big(c \cdot  L  + 4 \cdot c/( \xi \cdot \alpha)\big) \cdot dist(z,\p \mathfrak{B}^H_{{\cal{A}},\Lambda,a}).\]}

2.  For $p, q  \in arc_\Lambda$ we recall  that  $U_{\xi \cdot \alpha \cdot a/4}\left(arc_{\Lambda}(\E(a)) \right) \subset
        \mathfrak{B}_{{\cal{A}},\Lambda,a} \subset \E(\alpha \cdot a/4)$. Upon choosing $\xi$ small enough, we may assume $B_{\xi \cdot \alpha \cdot a/4}(z)$, $z  \in arc_\Lambda$,
is uniformly controlled so that after scaling by $(\xi \cdot \alpha \cdot a/4)^{-1}$ it is $\ve$-close to the unit ball in $\R^n$ in $C^3$-topology, for some tiny $\ve >0$.\\

2a. In the case where $d(p,q) \le \xi \cdot \alpha \cdot a/10$, we observe, that the shortest geodesic arc linking them satisfies the asserted inequalities for uniformly controlled constants.\\

2b. Otherwise, we may assume $d(p,q) \ge \xi \cdot \alpha \cdot a/10$ and  $p  \in arc_\Lambda  \setminus \E(a)$, in particular $\delta_{\bp}(p) \le a$,  belongs to some $c$-skin uniform arc
 $\gamma_{x,y} \subset \E(\alpha \cdot a)$ which joins two points $x(p),y(p) \in \E(a)$.\\

  Let $x$ be the one closer to $p$, in path length. Then we infer \,
 $l_{min}(\gamma_{x(p), y(p)}(p)) \le c \cdot \delta_{\bp}(p) \le c \cdot a$. We get a similar point $x(q)$ for $q$, the closer endpoint of another arc in $\Lambda$, and
 the same inequality for $q$, if  $q  \in arc_\Lambda  \setminus \E(a)$.
For $q  \in \E(a)$ we just take the constant path.\\

Now we define the curve $\gamma_{p,q}$ from $p$ to $q$ as follows:\begin{itemize}
    \item Starting from $p$  we follow some segment of some $c$-skin uniform arc $\in \Lambda$ of length $\le c  \cdot a$ until we reach $x(p) \in \E(a)$.
    \item Then we run from $x(p)$ to $x(q)$  along some $c$-skin uniform arc $\gamma_{x(p), x(q)} \in \Lambda$.
    \item Finally, we start from $x(q)$ and return to $q$ following the chosen short segment of another $c$-skin uniform arc $\in \Lambda$ of length $\le c  \cdot a$.
\end{itemize}
On the balance sheet, we get
\[l(\gamma_{p,q}) \le l(\gamma_{x(p), x(q)})+ 2 \cdot c \cdot a \le c \cdot d(x(p), x(q))+2 \cdot c  \cdot a \le c \cdot d(p,q)+
4 \cdot c \cdot a,\]
 with the lower bound  $d(p,q) \ge \xi \cdot \alpha \cdot a/10$, this shows\[l(\gamma_{p,q}) \le c \cdot d(p,q)+  4 \cdot c \cdot a \le (c + 40 \cdot c/(\xi \cdot \alpha)) \cdot d(p,q).\]
Also the estimate  $l_{min}(\gamma_{x(p),y(p)}(p)) \le  c \cdot  a$ and $arc_{\Lambda}(\E(a)) \subset \E(\alpha \cdot a)$ show that
\[ l_{min}(\gamma_{p,q}(z)) \le c \cdot a + c \cdot \delta_{\bp}(z) \le (c/\alpha + c )\cdot \delta_{\bp}(z)\]
Finally we append the estimate for the length of the connecting arcs in $\gamma_{x,y}$ from $p$ resp. $q$ to $arc_{\Lambda}(\E(a))$, that is $ c \cdot a$, to (\ref{con2})
and apply \ref{cch}(ii) to see
\begin{equation} \label{con3} l_{min}(\gamma_{p,q}(z)) \le \big(c \cdot  L  + 4 \cdot c/( \xi \cdot \alpha)\big) \cdot dist(z,\p \mathfrak{B}^H_{{\cal{A}},\Lambda,a})
+  c \cdot a \end{equation}
\[\le \big(c \cdot  L  + 4 \cdot c/( \xi \cdot \alpha) + 4 \cdot c/(\xi \cdot \alpha)\big) \cdot dist(z,\p \mathfrak{B}^H_{{\cal{A}},\Lambda,a}).\]

3. When $p, q \in \mathfrak{B}_{{\cal{A}},\Lambda,a}$ and $d(p,q) \ge \xi \cdot \alpha \cdot a/10$, we repeat the preceding arguments and additionally choose radial arcs from $p, q \in
\mathfrak{B}_{{\cal{A}},\Lambda,a}$ to points in $arc_\Lambda$. \qed

{\bf Proof of \ref{std}} \quad The properties collected in \ref{arc}, \ref{cch} and the results for two cases considered in \ref{uni} and \ref{uni2} above settle the proof of \ref{std} when we set
\begin{equation} \label{defco} \su(a):=\mathfrak{B}_{{\cal{A}},\Lambda,a},\end{equation}
for any link space $\Lambda$ of $\E(a)$ and some QT-skin adapted cover $\cal{A}.$\qed

\begin{remark} \textbf{(NTA Regularity)}\label{pu1}  \quad  A forerunner of the uniformity concept,  due to Jerison and Kenig [JK], is that of NTA-domains.
Being an NTA (=non-tangentially accessible) Euclidean domain $D$, can be described as being uniform with the additional exterior condition that for any boundary point $p \in \p D$ there is \emph{twisted exterior cones} in the complement of $D$, pointing to $p$ in the sense of \ref{jdo}. \\

At the time when NTA-domains were introduced the available techniques were not sufficient to treat potential theory and, in particular, Martin theory from interior conditions alone. Only later progress allowed to finally remove all exterior assumptions. This is, in particular, due to Bass and Burdzy [BB] and Aikawa [Ai1], [Ai2].\\

We use $\mathfrak{B}$-hulls only in the context of Martin theory. Hence, we can base our conclusions on the more recent results for uniform domains. However, the NTA-condition
has its own benefits and is incorporated in important other work. To open optional applications of $\mathfrak{B}$-hulls also in this context, we  outline here a way to improve the
uniformity of the $\su$, stated in \ref{std}, to that of being
NTA-domains.\\

The main point is that QT-properties of our QT-skin adapted cover $\cal{A}$
can be refined from pairwise quantitatively transversal intersections to its counterpart for finite intersections of balls in $\cal{A}$. The result now reads:\\

There is a transversality constant $\tau_n \in (0,1)$, depending only on the dimension $n$,
so that a cover ${\cal{A}}$ of $H \setminus \Sigma$, can be chosen so that, in addition to the properties in \ref{besi} we have the following intersection estimates
for any $m \le c(n)$ balls $B_{2 \cdot \Theta(p_1)}(p_1), ..., B_{2 \cdot\Theta(p_m)}(p_m) \in {\cal{A}}$:
{\small \begin{description}
\item[(QT1)*]  $\,B_{2 \cdot \Theta(p_1)}(p_1) \cap ...B_{2 \cdot\Theta(p_m)}(p_m) \n \, \Rightarrow \, B_{(2-\tau_n) \cdot \Theta(p_1)}(p_1) \cap ... B_{(2-\tau_n) \cdot\Theta(p_m)}(p_m) \n$,
\item[(QT2)*] $\,B_{2 \cdot \Theta(p_1)}(p_1) \cap ...B_{2 \cdot\Theta(p_m)}(p_m) \v \, \Rightarrow \, B_{(2+\tau_n) \cdot \Theta(p_1)}(p_1) \cap ... B_{(2+\tau_n) \cdot\Theta(p_m)}(p_m) \v$.
\end{description}}

The proof follows the lines of \ref{besi2} and inductively takes intersection of $m \le c(n)$ balls into account, with \ref{besi2} settling the case $m=2$. The additional twist is the formulation of the
shell conditions (S1) and (S2): in step $m+1$,  we take shells around the intersections of the $m$-balls and also in each step we get a counterpart to (\ref{ang}) now for solid angles to make sure
that these intersections do not move to far away when the radii of the involved balls are changed.\\

The benefit of this extension to (QT1)* and (QT2)*  is a refinement of the micro-uniformity to an NTA-property near the boundary. Then the skin uniformity can be used as before to
globalize the NTA-property to $\su$. In turn, note that for $H \setminus \Sigma$ there is no exterior part left and, thus, one could equally name it a skin NTA-space. \qed
\end{remark}

\setcounter{section}{5}
\renewcommand{\thesubsection}{\thesection}
\subsection{Appendix: Oriented Boundaries and Rectifiable Currents}\label{app}

\bigskip

In this appendix we collect some pieces of the geometric measure theory for area minimizing hypersurfaces we used and referred to in the main body of this paper but also in the two follow-ups [L1] and [L2].\\

\textbf{I.} \, \textbf{Existence of Area Minimizers} \quad A transparent approach to geometric measure theory is that of oriented boundaries we describe here. Good
sources for further reading are [AFP], [Gi], [M] and [MM].\\

The more
technical but also more powerful approach using currents is explained in part III below. \\

We start with some usually bounded open set $\Omega \subset \R^m$ and define $\int_\Omega|D f|$, for a function $f \in L^1(\Omega,\R)$ as
\[\int_\Omega|D f| := \sup \{\int_\Omega f \cdot \mbox{div} g  \, d\mu \;| \; g \in C_0^1(\Omega,\R^m), |g|_{C^0} \le 1 \}\] and $f$ is called of a \emph{function of bounded
variation} or just \emph{BV-function} if $\int_\Omega|D f| < \infty$. The \emph{BV-norm} of a BV-function $f$ is defined as
$$|f|_{BV(\Omega)} := |f|_{L^1(\Omega)} + \int_\Omega|D f|.$$ For the characteristic function $\chi_E$  of some Borel set $E \subset \R^m$ one calls $\int_\Omega |D \chi_E|$: the
\emph{perimeter} $P(E,\Omega)$ of $E$ in $\Omega$. The explanation for this name also shows us the relation of the $BV$-norm to geometry: $P(E,\Omega)$  equals the
$(m-1)$ - dimensional Hausdorff-measure ${\cal{H}}_{m-1}(\p E \cap \Omega)$ of the part of the boundary $\p E$ within $\Omega$, as soon as
$\p E$ is sufficiently regular, e.g. $C^2$-smooth, cf. [Gi], 1.4. \\

Thus  to find area minimizing hypersurfaces one starts with a perimeter minimizing sequence (of  characteristic functions) of Caccioppoli sets $E_j \subset \R^m$, i.e. Borel
sets with locally finite perimeter in  $\Omega$, subject to the condition that $E_j \setminus \Omega \equiv L \setminus \Omega$ for some given Caccioppoli sets $L$
considered as boundary values.  For the purposes of this paper, we refer to it as the \textbf{Plateau problem with boundary data} $\p L \cap \p \Omega$ within $\Omega$, cf.
[GMS],Vol.II,p.565
and [Ag] for the varifold reformulation.\\

We presume $\p \Omega$ is smooth enough so that Rellich compactness theorem holds on $\Omega$. That is, the embedding $$BV_{loc}(\Omega) \hookrightarrow
L^1_{loc}(\Omega) \mm{
is compact},$$ cf. [AFP], 3.23. Here we wrote the localized versions to include the case where $\Omega$ is non-compact, which happens to the case when we consider cones.\\

A sufficient condition for this compactness to hold is that $\Omega$ is an \emph{extension
domain}, cf. [AFP], 3.20 - 3.23 and 3.49, which includes the cases of Lipschitz regular boundaries and also the more general uniform domains cf. C.I.\\

 Then there exists a minimizing Caccioppoli set $E$ with $E \setminus \Omega \equiv L \setminus \Omega$:
\[P(E,\Omega) \le P(F,\Omega), \mbox{ for all Caccioppoli sets } F \mbox{ with } F \setminus \Omega \equiv L \setminus \Omega\] An \textbf{outline of the argument} reads as follows: Rellich
compactness shows that sets of functions uniformly bounded in BV-norm on $\Omega$ are relatively compact in $L_{loc}^1$ (cf. [Gi], 1.19). And actually the BV-norm of the
$\int_\Omega|D \cdot |$-minimizing sequence $\chi_{E_j}$ is uniformly bounded since $$0 \le \int_\Omega|D \chi_{E_j}|\mm{ and } 0 \le  |\chi_{E_j}|_{L^1(\Omega)} \le
Vol(\Omega).$$ Thus there is an $L_{loc}^1$-converging subsequence of the $\chi_{E_j}$ and since this implies convergence almost everywhere we may assume the
$L_{loc}^1$-limit is again a characteristic function $\chi_{E}$ for some $E \subset \R^m$. Now one uses the semicontinuity of BV-norms (cf. [Gi], 1.9), that is, for a sequence
of functions $f_j$ in $BV(\Omega)$ which converges in $L_{loc}^1$ to a function $f$ we have: $$\int_\Omega|D f| \le \liminf_{j \ra \infty} \int_\Omega|D f_j|.$$

Thus, we observe that $E$ is a Caccioppoli set with $$P(E,\Omega) \le \liminf_{j \ra \infty} P(E_j,\Omega)$$ and hence we may think of $\p E$ as an area minimizing
hypersurface within $\Omega$. This same method also proves the \textbf{compactness theorem} for the space $$\mathfrak{C}_\Omega\{ E_\alpha \mbox{ Caccioppoli set  }
\, | \, E_\alpha \setminus \Omega \equiv L \setminus \Omega \}:$$ any $\int_\Omega|D \cdot |$-bounded sequence, that is with uniformly bounded perimeter, in
this space has a $L_{loc}^1$-converging subsequence.\\

\textbf{II.} \, \textbf{Regularity Theory for  Almost Minimizers} \quad While the previously stated compactness results for $BV$-functions provide the existence of
some minimizer, the yet proven $L_{loc}^1$-regularity for such a minimizer is rather weak. However, De Giorgi and others developed a partial regularity theory for such
minimizers and proved the following classical result.

\begin{proposition}\emph{\textbf{(Partial Regularity of Minimizers)}} \label{dg} \, Let $E \subset \R^m$ be a minimizing Caccioppoli set, then $\p E$ can be written as an analytic hypersurface except for some singular set
$\Sigma$ of  Hausdorff-dimension $\le m-8$.
\end{proposition}

To formulate this result we actually merged two stages in the development of the theory in \ref{dg}. The first one is the genuine regularity theory and secondly, much lighter in weight, the
estimate for the size of $\Sigma$. \\

\textbf{Stage 1} \, De Giorgi realized that \emph{partial} regularity of $\p E$ is the achievable goal to go for and proved central regularity results. Miranda refined his methods
and the final outcome was that $\p E$ can be written as an analytic hypersurface in $\R^m$ except for some singular set of Hausdorff-dimension $\le m-1$. Actually, the proof
gives us a description of $\p E$ as a \emph{countable disjoint union} of pieces of \emph{$C^1$-hypersurfaces} and some remaining set $\Sigma$
of $m-1$-measure zero.\\

 Note that the non-constructive methods to detect the $C^1$-pieces renders $\Sigma$ as a terra incognita without a proper internal structure.\\

To further explain the nature of the regularity result we consider an already more general case of \emph{almost minimizers} due to Tamanini [T1],
[T2], [MM], Bombieri [Bo] and Allard [A]. The following result is the basic content of [T1],Theorem 1.
\begin{proposition}\label{arm} \, Let $\Omega \subset \R^m$ be open, $E \subset \R^m$ be a Caccioppoli set, which is almost minimizing in $\Omega$ in the sense that
the following estimate holds for some $K > 0$, $\alpha \in (0,1)$ {\small $$\psi(E, B_\rho(x)) := \int_{B_\rho(x)}|D \chi_{E}| - \inf \left\{\int_{B_\rho(x)}|D \chi_{F}|\, \Big| \, F
\Delta E \subset \! \subset B_\rho(x) \right\} \le K \cdot \rho^{m-1+2 \cdot \alpha}$$}
for any $x \in \Omega$, $\rho \in (0,R)$, for some $R > 0$ (where $F \Delta E := F \setminus E \cup  E \setminus F$).\\

Then $\p E \cap \Omega$ can be written as a $C^{1, \alpha}$- hypersurface except for some singular set of Hausdorff-dimension $\le m-8$.
\end{proposition}

\begin{remark} \label{riem} \, 1. This includes the area minimizing case where $\psi \equiv 0$. The regularity statement means that we can locally  find a hyperplane $\subset \R^m$  so that $\p E$ is
the graph of some $C^{1, \alpha}$-function.\\

2. One readily checks that the minimal hypersurfaces $S_C=\p B_1(0) \cap C$, for any area minimizing cone $C$, are almost minimizers. The $S_C$ are critical points of the area functional, but they are not area minimizing and they cannot even be stable, since $Ric_{S^n} >0$.\\

3. Images of almost minimizers under diffeomorphisms of the ambient space remain to be almost minimizers. Since this is a statement on the behavior in arbitrarily small
balls around points in the image of the almost minimizer we only need to consider the case of a linear isomorphism $\phi$ between Euclidean spaces.\\

Now in this case, there is some $k \ge 1$ so that $1/k  \cdot |v| \le D\phi(v) \le  k  \cdot |v|$, for any $v \in \R^n$ and therefore {\footnotesize
 $$\int_{B_\rho(\phi(x))}|D \chi_{\phi(E)}| - \inf
\left\{\int_{B_\rho(\phi(x))}|D \chi_{F}|\, \Big| \, F \Delta \phi(E) \subset \! \subset B_\rho(\phi(x)) \right\} \le $$
$$k^{n-1} \cdot \int_{B_{k \cdot \rho}(x)}|D \chi_{E}| - \inf \left\{\int_{B_{k \cdot \rho}(x)}|D \chi_{G}|\, \Big| \, G
\Delta E \subset \! \subset B_{k \cdot \rho}(x) \right\}$$}

It is noteworthy that this does \textbf{not} imply that bounds on the mean curvature towards singular points would be preserved under local diffeomorphisms.\\

4. With these results for almost area minimizers we also get the regularity theory also for area minimizers in Riemannian manifolds since we can now use local
coordinate charts. Under these maps area minimizers in a manifold become an almost minimizer in $\R^m$. Then we apply \ref{arm} and transfer the regularity result back
to the manifold.\\

Further improvements of the smoothness result are owing to standard elliptic theory. For instance, this shows that when the ambient manifold is analytic the regular portion of an
area minimizing hypersurface is also analytic cf. [Mo], 5.7. \qed
\end{remark}

Two important consequences of \ref{arm} (part of [T1],Theorem 1) say that a sequence of almost minimizers $E_i$ converging to some limit $E_\infty$ will eventually become
smooth near smooth limit points in $\p E_\infty$. And, secondly, the $L^1_{loc}$-convergence implies a $C^1$-convergence when the limit is known to be $C^{1,
\alpha}$-smooth. The latter may not be evident from [T1] but could be derived from its sources. Explicitly this can be located in Allard's work, cf. [A] or  [Si2], 23.1.\\

 Technically,
this reaches to the heart of the regularity theory and it exploits the averaged oscillation of the normal vector, the \emph{excess}. Since the control and effect of thresholds for the
excess are the clue to the following result
\begin{corollary}\label{carm} \quad Let $E_i$, $i \ge 0$ be Caccioppoli sets satisfying the condition $$\psi(E_i, B_\rho(x)) \le K \cdot \rho^{m-1+2 \cdot \alpha} \mbox{ for fixed } K > 0 \mbox{ and }
 \alpha \in (0,1)$$ in some open $\Omega \subset \R^m$.
 \begin{enumerate}
\item  Assume the $E_i$ converge in $L^1_{loc}$ to  $E_\infty$, there are points $p_i \in \p E_i$ which converge to some $p_\infty \in \p E_\infty$ and that $p_\infty$ is
    a smooth point in $\p E_\infty$. Then, for large $i$, the $p_i$ are also smooth points $x_i \in \p E_i$.
\item  If the limit $E_\infty$ in (i) has a $C^{1, \alpha}$-boundary in $\Omega$, then $\p E_i$ converges in $C^1$-topology to $\p E_\infty$.
\end{enumerate}
\end{corollary}

Note that $E_\infty$ also satisfies  the condition $\psi(E_\infty, B_\rho(x)) \le K \cdot \rho^{m-1+2 \cdot \alpha}$ which can be proved as in [Gi],9.1. Thus, again, this result
may be used in the Riemannian world and, also,  elliptic theory allows us to improve the $C^1$-convergence to $C^k$-convergence for any $k \ge 0$ when the ambient
Riemannian manifold is of class $C^\infty$. In other words, \emph{a flat norm converging sequence of area minimizers will be locally $C^k$-converging around smooth
points of the limit surface}, cf. [Gi],11.4 for more details.\qed

\textbf{Stage 2} \, In \ref{dg} we also asserted an estimate for the dimension of the singular set. It is based on a remarkable use of \ref{carm}.\\

Federer, building on work of Bombieri, Giusti, De Giorgi, Simons, Almgren and others, founded a way to estimate the dimension of the singular set more efficiently. He used
approximations by tangent cones to derive estimates for their size by a dimensional induction starting from the fact that all minimizing cones in dimension $\le 7$ are smooth
hyperplanes and  thus singularities are isolated points in dimension $8$. The result is
\begin{proposition}\emph{\textbf{(Partial Regularity of Almost Minimizers)}} \label{dg8} \, Let $E \subset \R^m$ be a almost minimizing Caccioppoli set in the sense of \ref{arm}, then the singular set
$\Sigma^{m-8} \subset \p E$ is a closed subset of (potentially locally varying) Hausdorff-dimension $\le m-8$.
\end{proposition}
The transfer to Riemannian manifolds is not a problem since the inductive cone reduction - which is the main new ingredient at this stage - considers an infinitely scaled
ambient geometry, that is, the same Euclidean geometry which is used for area minimizers in $\R^m$.\qed

\textbf{Standard Applications} \, This regularity theory is typical for hypersurfaces and does not apply area minimizers in higher codimension.
For reference in the main body of the paper, we explicitly state some folklore applications we could not properly localize in the literature. They are valid, again,
only in the hypersurface case.

\begin{corollary}\label{tt} \quad Let $D \subset H$ an open domain in an oriented minimal boundary $H$. Then we the following statements are equivalent
\begin{itemize}
\item All points in $D$ are regular (manifold points).
\item Around any point $p \in D$, $|A| \le c$, for some $c(p) >0$.
\item For any $p \in D$, the hyperplane is a tangent cone.
\end{itemize}
\end{corollary}

Next we consider a \emph{non-extinction} result for oriented minimal boundaries in $\R^{n+1}$, we need to properly apply compactness results.
It shows that opposing sheets of sequences of such minimizers cannot approach each others too closely and annihilate in the limit. Implicitly this is contained in the estimate [Gi], Prop.5.14. For the sake of completeness we include an argument.

\begin{lemma}  \label{nex} \quad   Let $H_i \subset \R^{n+1}$ be a sequence of oriented minimal boundaries with $0 \in H_i$ and $|A| \le 1$ on $B_2(0) \subset H_i$. Then, for any
 compactly converging subsequence $H_{i_k}$, we have
 \begin{itemize}
   \item the limit hypersurface $H_\infty$  is an oriented minimal boundary,
   \item $0 \in H_\infty$ and $|A| \le 1$ on $B_1(0) \subset H_\infty$,
   \item $B_1(0) \subset H_{i_k}$ converge smoothly to  $B_1(0) \subset H_\infty$, in the sense of $\D$-maps.
 \end{itemize}
\end{lemma}

{\bf Proof} \quad We show that, in $\R^{n+1}$, the ball $B_1(0) \subset H_i$  is not approached from $O_i:=H_i \setminus B_2(0)$
when $i \ra \infty$. That is there is a lower positive distance bound between these balls and the $O_i$, independent of $i$.\\

Otherwise, we get for $i \gg 1$, a subset of $O_i$ which can be written as a smooth graph $G_i$ over $B_1(0) \subset H_i$ arbitrarily close to $B_1(0)$ in $C^3$-norm.
This is a standard consequence of DeGiorgi-Allard regularity theory, cf. [Si],24.2 and the Harnack inequality cf.[So],p.73.\\

Since the $H_i$ bound open sets $U_{H_i}  \subset \R^{n+1}$ we may assume that $G_i$ and $B_i$ are oppositely oriented.  For $i \gg 1$ we consider  $G_i \cup B_i$ and  join
$\p G_i$ and $\p B_i$ linearly through some hypersurface $F_i$. Then we add the bounded open set $V_i  \subset \R^{n+1}$ with\[\p V_i = G_i \cup B_i \cup F_i \]  to
$U_{H_i}$ and observe that this open set $U_{H_i} \cup V_i$ has an oriented boundary \[\p (U_{H_i} \cup V_i) \equiv H_i  \mm{ relative } \R^{n+1} \setminus \overline{V_i}.\]
And relative $\overline{V_i}$, it has smaller area than $H_i$. For the latter assertion, we note that the $B_i$ have uniformly bounded geometry since $|A_{H_i}|\big|_{B_2(0)} \le 1$.\\

For large $i$,  $\p (U_{H_i} \cup V_i)$ is an admissible compactly supported perturbation of $H_i$. This contradicts the area minimizing property of $H_i$ and shows
that $O_i$ remains in a positively lower bounded distance of $B_i$, for all $i$. The remaining assertions follow from the standard regularity theory.\qed

Also, we note some weak type of a Harnack inequality or a quantitative identity theorem for $|A|$.

\begin{lemma}  \label{haa}  For any $\lambda \in (0,1]$ and given radius $R_0 >0$, there is a constant $c(\lambda, n, R_0) >0$, so that
for any oriented minimal boundary  $H \subset \R^{n+1}$ and any $p \in H$ with \[ \sup \{|A|(x) \, | \, x \in B_{R_0}(p) \cap H\} \ge 1.\] we have the positive lower
estimate
 \[ \sup \{|A|(x) \, | \, x \in B_{\lambda  \cdot R_0}(p) \cap H\} \ge c.\]

\end{lemma}

 {\bf Proof} \quad It suffices to check the case where $R_0 =1$ and we may assume $p=0 \in \R^{n+1}$. Assume there is no such constant $c > 0$. Then there is some $\lambda \in (0,1)$ and a sequence of such
hypersurfaces $H_k$ and points $0 \in H_k$ so that
$\sup \{|A|(x) \, | \, x \in B_\lambda(0) \cap H_k\} \le 1/k.$\\

Due to the minimality of the $H_i$  a subsequence of the $H_k$ which converges compactly on $\R^{n+1}$ in flat norm to some limit hypersurface $H_\infty$. As in the preceding result we may assume from $\sup \{|A|(x) \, | \, x \in B_\lambda(0) \cap H_k\} \le 1/k$ that, in $B_\lambda(0)$, this is a $C^k$-convergence, for some $k \ge 5$.\\

Then the analytic
minimizer $H_\infty$ is a hyperplane, since the limit of the $B_\lambda(0) \cap H_k$ in  $H_\infty$ must be flat. Note that $H \setminus \Sigma$ is connected, for any oriented boundary $H$, cf.\ref{coh}.\\

But then regularity theory transforms the flat convergence to $C^k$-convergence, also outside $B_\lambda(0)$. This shows that $\sup \{|A|(x) \, | \, x \in B_R(0) \cap H_k\} \ra 0$, for any given $R >0$, which contradicts the assumption.\qed

\textbf{III.} \, \textbf{Currents} \quad Currents can be viewed as a generalization of submanifolds, or better to say chains and cochains. The basic concepts, due de Rham and
Whitney [R],[Wh2] were introduced in a differential topological context. Later they became the basis for geometric measure theory that formally includes the approach of
oriented boundaries we outlined above, cf. [F], [GMS] or [Si2] for comprehensive presentations.\\

The idea is to enlarge the space of submanifolds to the space of $m$-currents $\mathcal{D}_m(U)$, with $U \subset \R^n$ open, which is the dual of $\mathcal{D}^m(U)$,
the space of smooth $m$-forms compactly supported in $U$.
Integrating the forms over a submanifold $N$ we can be interpret $N$ as a current denoted by $\llb F \rrb$.\\

Towards a topological use of currents, we define the boundary $\p T  \in \mathcal{D}_{m}(U)$ of any current $T \in \mathcal{D}_{m+1}(U)$ through \[\p T(\omega):= T(d \omega), \mm{ for any compactly supported }
(m-1)\mm{-form }\omega.\]

The support $supp \: T $ of a current $T$ is the complement of the union of all open sets $W$ with $T (\omega) = 0$ for  $\omega \in \mathcal{D}^n(U)$ with
$supp \; \omega \subset W.$ For any open $W \subset U$ and $T  \in \mathcal{D}_{m}(U)$ we write $T \llcorner W$ for the current in $\mathcal{D}_{m}(W)$ we get from restricting $T$ to $W$. \\

For any compactly supported current  $T  \in \mathcal{D}_{m}(U)$ we define its push-forward $f_\sharp T$ by
\[f_\sharp T(\omega):= T(f^*\omega), \mm{ for the pull-back } f^* \omega \mm{ of any }m\mm{-form } \omega.\]
The weighted area of   $T  \in \mathcal{D}_{m}(U)$, its \emph{mass} ${\bf{M}}_U (T)$, is defined by
\[{\bf{M}}_U (T) = \sup_{|\omega| \le 1, \mbox{{\tiny{supp}}} \omega \subset U}T(\omega).\]
We define the flat metric topology on $\mathcal{D}_{m}(U)$: For any open subsets of $\R^n$: $W \subset \overline{W} \subset U
\subset \mathbb{R}^n$. Roughly, the \emph{flat (pseudo)metric} measures the volume between two currents $C_1, C_2 \subset \mathcal{D}_{n-1}(U)$:
\begin{equation}\label{fll}
\db_W(C_1,C_2) :=
\end{equation}
\[ \mbox{inf}\{{\bf{M}}_W (S) +
{\bf{M}}_W (R)\, |\, C_1 - C_2 = S + \p R, S \in \mathcal{D}_{n-1}(U), R \in \mathcal{D}_n(U)) \}.\]
The family of these $\db_W$ generate the flat metric topology. When $W$ is the entire ambient space, that is, $\R^n$ or the given manifold, we drop the index $W$ and merely write $\db$.\\

The geometrically most relevant spaces of currents are the spaces $\mathcal{R}_m(U) \subset \mathcal{D}_m(U)$ of integer multiplicity \emph{rectifiable currents} and
$\mathcal{I}_m(U) \subset \mathcal{D}_m(U)$ of  \emph{integral currents}. Integral currents are those rectifiable currents with rectifiable boundary.\\

$T \in \mathcal{D}_m(U)$ is an integer multiplicity rectifiable current, if for any $\ve >0$ and any compact set $K \subset U$ there is a compactly supported $m$-dimensional
polyhedral chain with integer coefficients $P=P(K, T,\ve) \subset \R^k$ of oriented simplices and a Lipschitz map $f:\R^k \ra \R^n$ so that $supp \, f_\sharp P \subset K$ and
\[{\bf{M}}_U (T-f_\sharp P) < \ve.\]

As in the BV-approach there are compactness results for integral currents now expressed in flat metric topology in place of $L^1$-topology.\\

These concepts and results extend to any compact manifold $M^n$, via local charts. The following existence result follows from compactness results for subsets of
$\mathcal{I}_{n-1}(M)$ under mass bounds. For the existence we refer to [F1],4.2.17, 4.4.5 and 5.1.6, or [GMS], 5.4.1,Cor.1. For the regularity theory, cf. [F1],5.3 and for the dimension of the singular set [F2].

\begin{proposition}\emph{\textbf{(Homological Minimizers)}}\label{hm} \, For any $\alpha \in H_{n}(M^{n+1}, \Z)$ there is a mass minimizing integral current $X^n \in \alpha$ whose support is a smooth
hypersurface except for some singular set $\Sigma_X$ of codim $\ge 8$ in $M^{n+1}$.
\end{proposition}

The statement refers to the homology of integral currents defined using the boundary operator $\p$ on the space of integral currents.\\

Singular homology and integral current homology are isomorphic when the ambient space is sufficiently regular, for instance a smooth manifold,
cf. [Dp] and [H] for details.\\

\textbf{IV.} \, \textbf{Decomposition of Rectifiable Currents} \quad In a rather concrete sense, rectifiable currents can be understood from the case of oriented
boundaries of measurable sets. Namely, there is a decomposition theorem for these currents into oriented boundaries, cf.[F],4.5.17 and also [Si],Ch.37 or
[GMS], I.4.3.1,Th.7.

\begin{proposition}\label{dic} \quad  For any $R \in \mathcal{R}_n(\R^{n+1})$, with $\p R \v$, there exist measurable sets $A_i \subset \R^{n+1}$, $i \in \Z$, $A_i \subset A_{i+1}$  such that for any bounded open $W \subset \R^{n+1}$:
  \[R= \sum_{i \in \Z} \p \llb A_i \rrb \,\mm{ and } \, {\bf{M}}_W (R)= \sum_{i \in \Z} {\bf{M}}_W (\p \llb A_i \rrb).\]
\end{proposition}

In the case of a local mass minimizing $R$ one may assume that the sets $A_i$ are open and the $ \p \llb A_i \rrb$ are oriented boundaries,
each of them minimizes the perimeter in the BV-sense.\\

There is also a localized version of \ref{dic} for currents in a manifold $M^{n+1}$ which can be directly be derived from \ref{dic}.\\

When $U$ is a proper ball in $M$ and $R \in \mathcal{R}_n(U)$, we define some diffeomorphism $f: U \ra \R^{n+1}$,
apply \ref{dic} to $f_\sharp R \in \mathcal{R}_n(\R^{n+1})$ and consider the pull-back of the resulting decomposition on $U$ with the retransformed masses.\\

More generally,  we can use the long exact homology sequence for the homology of integral currents to infer that for an open subset $U \subset M$ with $H_n(M, M\setminus U)
= 0$ there is a decomposition relative $U$.\\

The following result is a version of this local decomposition adapted to the case of area minimizers.

\begin{proposition} \emph{\textbf{(Local Decompositions)}} \label{dicc} \,
Any locally mass minimizing current $T \in \mathcal{R}_n(U)$, with $\p T =0$, in an open subset $U \subset M$ with $H_n(M, M\setminus U)
= 0$, can be decomposed in $U$ into oriented boundaries $\p M_i$, for open sets $M_i \subset
U$, $i \in \Z$, $M_i \subset M_{i+1}$, each locally area minimizing in $U$
 so that for any open $W \subset \! \subset U$
\begin{equation}\label{ld}
 T\llcorner U= \sum_{i \in \Z} \p \llb M_i \rrb \llcorner U \,\mm{ and } \, {\bf{M}}_W (T\llcorner U)= \sum_{i \in \Z} {\bf{M}}_W (\p \llb M_i \rrb \llcorner U).
\end{equation}
\end{proposition}

In general, the latter sums may contain infinitely many terms. But we observe that for a compact manifold $M$ and $T$ a mass minimizing current, that represents a given homology class $\alpha \in H_{n}(M^{n+1}, \Z)$, this sum is finite:\\

To see this, we choose a small ball $B_{5 \cdot r}(p) \subset U \subset M$, so that $(5 \cdot r)^{-1} \cdot B_{5 \cdot r}(p)$ is nearly isometric to the unit ball in $\R^{n+1}$.
Then the minimality of each $\p M_i$ which intersects $B_r(p)$ gives the estimate [Gi],Ineq.(5.16):
${\bf{M}}_{B_{2 \cdot \rho}(p)} (\p \llb M_i \rrb \llcorner U \ge c_n \cdot r^n$, for some constant $c_n >0$ depending only on $n$, for each of its components.
Thus the finiteness of the total mass of $T$ shows that there are only finitely many $\p M_i$ involved.\\

Finally, we note, from the indicated arguments, that the term \textbf{local} refers to the choice of a suitable set $U$ in the ambient manifold, \emph{independent} of the given current.
This allows us to use the arguments, within the fixed set $U$, also when we consider converging sequences of such currents.\\

An instructive example how such decompositions can be used to reduce problems to oriented boundaries can be found in Simon's proof of the strict maximum principle for
minimizing integral currents in [Si3],Ch.2. In turn, from this maximum principle, the oriented boundaries in the sum (\ref{ld}) are locally either disjoint or, in cases of currents of higher multiplicities, they are equal.\\


\small
\begin{thebibliography}{LLLL}
\bibitem[Ai1]{Ai1} Aikawa, H.: Boundary Harnack principle and Martin boundary for a
uniform domain, J. Math. Soc. Japan 53(1) (2001), 119-145.
\bibitem[Ai2]{Ai2} Aikawa, H.: Potential analysis on nonsmooth domains - Martin boundary and boundary Harnack principle,
in Complex Analysis and Potential Theory, Vol. 55 of CRM, Proc. Lecture Notes, AMS (2012), 235-253
\bibitem[Ai3]{Ai3} Aikawa, H.: Potential-theoretic characterizations of nonsmooth domains, Bull. London Math. Soc. 36 (2004), no. 4, 469-482
\bibitem[A]{A} Allard, W.K.: On the first variation of a varifold, Ann. of Math. 95 (1972), 417-491
\bibitem[Ag]{Ag}  Almgren, F.: Plateau's problem: an invitation to varifold geometry, revised edition, AMS (2001)
\bibitem[AG]{AG} Armitage,D. and Gardiner,S.: Classical Potential Theory, Springer-Verlag (2001)
\bibitem[AFP]{AFP} Ambrosio, L., Fusco, N. and Pallara, D.: Functions of Bounded Variation and Free Discontinuity Problems, Oxford Math. Monographs, Clarendon Press (2000)
\bibitem[An1]{An1} Ancona, A.:  Th\'{e}orie du potentiel sur les graphes et les vari\'{e}t\'{e}s, in:  Ecole d'\'{e}t\'{e} de Prob. de Saint-Flour XVIII-1988, LNM 1427, Springer (1990), 1-112
\bibitem[An2]{An2} Ancona, A.: Negatively curved manifolds, elliptic operators, and the Martin boundary, Ann of Math 125 (1987), 495-536
\bibitem[An3]{An3} Ancona, A.: On strong barriers and an inequality of Hardy for domains in $\R^n$, J. London Math. Soc. 34 (1986), 274-290
\bibitem[BEL]{BEL} Balinsky, A., Evans, W.D. and Lewis, R.: The Analysis and Geometry of Hardy's Inequality, Springer (2015)
\bibitem[BB]{BB} Bass, R. F.  and Burdzy, K.: A boundary Harnack principle in twisted Hölder domains, Ann. of Math. 134( (1991), 253-276
\bibitem[Bo]{Bo} Bombieri, E:  Regularity theory for almost minimal currents, Arch. Rational Mech. Anal. 78 (1982) 99-130
\bibitem[BG]{BG} Bombieri, E.  and Giusti,E. : Harnack's inequality for elliptic differential equations on minimal surfaces, Invent. Math. 15 (1972) 24-46
\bibitem[BHK]{BHK} Bonk, M., Heinonen,  J. and Koskela, P.: Uniformizing Gromov hyperbolic spaces, Ast\'{e}risque 270 (2001)
\bibitem[BJ]{BJ} Borel, A. and Ji,L.: Compactifications of Symmetric and Locally Symmetric Spaces, Birkhauser (2006)
\bibitem[BM]{BM} Brezis, H. and Marcus, M.: Hardy's Inequalities Revisited, Ann. Scuola Norm. Sup. Pisa Cl. Sci. Vol. XXV (1997), 217-237
\bibitem[Ch]{Ch} Chavel, I.: Eigenvalues in Riemannian Geometry, Academic Press (1984)
\bibitem[CHS]{CHS} Caffarelli, L. A.  Hardt, R. and Simon, L.: Minimal surfaces with isolated singularities, Manuscripta Math. 48 (1984),1-18
\bibitem[CN]{CN} Cheeger, J.  and Naber, A.: Quantitative Stratification and the Regularity of Harmonic Maps and Minimal Currents,
Comm. Pure and Appl. Math.  66 (2013), 965-990
\bibitem[D]{D} De Giorgi, E.: Frontiere orientate di misura minima, Sem. Mat. Sc. Norm. Pisa (1961), 1-56
\bibitem[DS]{DS} David, G. and  Semmes, S.: Quasiminimal surfaces of codimension 1 and John domains. Pacific J. Math. 183 (1998), 213-277
\bibitem[Di]{Di} DiBenedetti, E.: Real Analysis, Birkhaeuser (2002)
\bibitem[Do]{Do} Doob,J.: Classical Potential Theory and its Probabilistic Counterpart. Springer, New York (1984)
\bibitem[Dp]{Dp} De Pauw, T.: Comparing homologies: Cech's theory, singular chains, integral flat chains and integral currents, Rev. Mat. Iberoam. 23 (2007), 143-189
\bibitem[F1]{F1} Federer, H.: Geometric Measure Theory, Spinger Verlag, Berlin  (1969)
\bibitem[F2]{F2} Federer, H.: The singular set of area minimizing rectifiable currents with codimension one and of area minimizing chains modulo two with
arbitrary codimension, Bull. AMS 76 (1970), 767-771
\bibitem[Fa]{Fa} Falconer, K.: Fractal Geometry, Mathematical Foundations and Applications, 3rd ed., Wiley (2014)
\bibitem[GMS]{GMS} Giaquinta, M., Modica, G. and Sourek, J.: Cartesian Currents in the Calculus of Variations, Vol. I and II, Springer (1998)
\bibitem[GT]{GT} Gilbarg, D. and Trudinger,N.: Elliptic Partial Differential Equations of Second Order,2nd edition, Springer Verlag, Berlin (1983)
\bibitem[GN]{GN} Ghoussoub, N. and Moradifam, A.: Functional Inequalities: New Perspectives and New Applications, Math. Surveys and Monographs 187,  AMS (2013)
\bibitem[G]{G} Gromov, M.: Plateau-Stein manifolds, Central Europ. Journal Math. 12, (2014), 923-951
\bibitem[GL]{GL} Gromov, M. and Lawson, B.: Positive scalar curvature and the Dirac operator on complete Riemannian manifolds, Publ. Math. IHES 58 (1983), 295-408
\bibitem[GO]{GO} Gehring, F.W. and Osgood,B.G, Uniform domains and the quasihyperbolic metric, Journal d'Analyse Math. 36 (1979), 50-74
\bibitem[Gi]{Gi} Giusti, E.: Minimal Surfaces and functions of bounded variations, Birkhaeuser Verlag, Basel (1984)
\bibitem[H]{H} Hardt, R.: Plateau  Problems  in  Metric  Spaces  and  Related  Homology, in Geometry and Topology of Submanifolds and Currents, Contemp.Math.646, AMS (2015),1-17
\bibitem[HS]{HS} Hardt, R. and Simon, L.: Area minimizing hypersurfaces with isolated singularities, Crelle Journal, 362 (1985), 10-129
\bibitem[Hz]{Hz} Heinz, E.: \"Uber die L\"osungen der Minimalfl\"achengleichung; Nachr. Akad. Wiss. G\"ottingen. Math.-Phys. Kl. (1952), 51-56
\bibitem[He]{He} Herron, D.. Uniform spaces and Gromov hyperbolicity,  Quasiconformal Mappings and their Applications, Narosa (2007), 79-115
\bibitem[JK]{JK} Jersion, D. and Kenig, C.: Boundary behaviour of harmonic functions in nontangentially accessible domains, Adv. in Math. 46 ,No.1(1982), 80-147
\bibitem[K]{K} Koskela, P.: Old and New on the Quasihyperbolic Metric, in, Quasiconformal mappings and analysis,
Springer-Verlag,  205-219 (1998)
\bibitem[L1]{L1} Lohkamp, J.: Hyperbolic Geometry and Potential Theory on Minimal Hypersurfaces, Preprint
\bibitem[L2]{L2} Lohkamp, J.: Skin Structures in Scalar Curvature Geometry, Preprint
\bibitem[L3]{L3} Lohkamp, J.: Positive scalar curvature in dim $\ge 8$, C. R. Acad. Sci. Paris, Ser. I 343 (2006), 585-588
\bibitem[M]{M} Maggi, F.: Sets of Finite Perimeter and Geometric Variational Problems, Cambridge Studies in Adv. Math. (2012)
\bibitem[MM]{MM} Massari, U. and Miranda, M.: Minimal Surfaces of Codimension One, Math. Studies 91, North-Holland (1984)
\bibitem[Mo]{Mo} Morrey, C.: Multiple integrals in the calculus of variations, Springer 1966
\bibitem[R]{R} de Rham, G.: Differentiable manifolds: Forms, Currents, Harmonic Forms, Spinger (1984)
\bibitem[SG]{SG} Shilov, G.E. and Gurevich, B.L.: Integral, Measure and Derivative, Dover Publications, New York (1977)
\bibitem[SU]{SU} Schoen, R.  and Uhlenbeck, K.: A regularity theory for harmonic maps, J.  Diff. Geom. 17 (1982), 307-335
\bibitem[SY]{SY} Schoen, R. and Yau, S.T.: Existence of
incompressible minimal surfaces and the topology of three dimensional manifolds with non-negative scalar curvature, Ann. of Math. 110 (1979), 127-142
\bibitem [Si1]{Si1} Simon,L.: Rectifiability of the singular set of energy minimizing maps, Calc.of Var.P.D.E.3(1995),1-65
\bibitem[Si2]{Si2} Simon, L.:  Lectures on Geometric Measure Theory, Proceedings of the Centre for Mathematical Analysis, Australian
National University, Canberra, 1983
\bibitem[Si3]{Si3} Simon, L.: A strict maximum principle for area minimizing hypersurfaces, J. Diff. Geom. 26, 327-335, (1987)
\bibitem[So]{So}  Solomon, B.: On foliations of $\R^{n+1}$ by minimal hypersurfaces. Comment. Math. Helv. 61 (1986), 67-83
\bibitem[St]{St} Stein, E.: Singular Integrals and Differentiability Properties of Functions, Princeton Math. Series 30, (1971)
\bibitem[T1]{T1} Tamanini, I.: Boundaries of Cacciopoli sets with H\"older continuous normal vector, J. Reine Angew. Math., 334 (1982), 27-39
\bibitem[T2]{T2} Tamanini, I.: Regularity Results for Almost Minimal Oriented Hypersurfaces in $\R^n$, Quaderni del Dipartimento di Matematica dell
Universit\`{a} di Lecce 1 (1984)
\bibitem[V]{V} V\"ais\"al\"a J.: Relatively and inner uniform domains. Conform. Geom. Dyn. 2, (1998), 56-88
\bibitem[W]{W} White, B.: Stratification of minimal surfaces, mean curvature flows, and harmonic maps. J. Reine Angew. Math. 488 (1997),1-35
\bibitem[Wh1]{Wh1} Whitney, H.: Analytic extensions of differentiable functions defined in closed sets, Trans. AMS, 36 (1934), 63-89
\bibitem[Wh2]{Wh2} Whitney, H.: Geometric Integration Theory, Princeton (1957)
\bibitem[YC]{YC} Yomdin, Y. and Comte, G.: Tame Geometry with Application in Smooth Analysis,
LNM 1834, Springer (2004)

\end{thebibliography}
\end{document}